%% file: paper-peristaltic-flow-2024-R.tex
\documentclass[preprint,12pt]{elsarticle_arxiv}

\usepackage{graphicx}
\usepackage{placeins}
\usepackage{subcaption}

\usepackage{amssymb}
\usepackage{amsmath}
\usepackage{amsthm}
\usepackage{amsfonts}

\usepackage{hyperref}
\usepackage[usenames,dvipsnames]{xcolor}
\usepackage{multirow}
\usepackage{xspace}
\usepackage{esint}
\usepackage{color}

\usepackage{changes}
\definechangesauthor[name={VL}, color=orange]{vl}
\definechangesauthor[name={ER}, color=orange]{er}

\usepackage[mathlines]{lineno}

\newcommand\revE[1]{#1}
\newcommand\newE[1]{#1}
\newcommand\chE[1]{#1}
\newcommand\revVa[2]{#1}{}  
\newcommand\revVr[2]{#1}{}  

\graphicspath{.}
\def\pathONE{.}

\newtheorem{remark}{Remark}[section]

\def\eq#1{(\ref{#1})}

\def\bmi#1{\textbf{\textit{#1}}}
\def\ol#1{\overline{#1}}

\def\Dvg{\mbox{\rm Div}}

\def\zerobf{{\bf 0}}
\def\onebm{{\bf{1}}}

\def\sigmabf{{\boldsymbol{\sigma}}}
\def\taubf{{\boldsymbol{\tau}}}
\def\omegabf{{\boldsymbol{\omega}}}
\def\etabf{{\boldsymbol{\eta}}}
\def\vthetabf{{\boldsymbol{\vartheta}}}
\def\xibf{{\boldsymbol{\xi}}}
\def\chibf{{\boldsymbol{\chi}}}

\def\Dlt{\Delta}
\def\dlt{\delta}

\def\veps{\varepsilon}
\def\vepsz{{\varepsilon_0}}
\def\Hveps{{\homlabel,\varepsilon_0}}

\def\vphi{\varphi}
\def\vtheta{\vartheta}
\def\pd{\partial}
\def\dd{\,{\rm{d}}}

\def\Om{\Omega}
\def\OmR{{\Omega_R}}

\def\hOmSe{{\hat\Omega_S^\veps}}
\def\hOmFe{{\hat\Omega_F^\veps}}
\def\hOmMe{{\hat\Omega_m^\veps}}

\def\hGmSFe{{\hat\Gamma_{\it{SF}}^\veps}}

\def\om{\omega}
\def\RR{{\mathbb{R}}}
\def\intY{\fint}

\def\dV{\mathrm{\,dV}}  %
\def\dVx{\mathrm{\,dV}_{\kern-0.1em x}}  %
\def\dVy{\mathrm{\,dV}_{\kern-0.1em y}}  %
\def\dVz{\mathrm{\,dV}_{\kern-0.1em z}}  %
\def\dVxy{\mathrm{\,dV}_{\kern-0.1em xy}}  %
\def\dVxz{\mathrm{\,dV}_{\kern-0.1em xz}}  %
\def\dY{\mathrm{\,dV}_{\kern-0.1em y}}  %
\def\dS{\mathrm{\,dS}}  %
\def\dSx{\mathrm{\,dS}_{x}}  %
\def\dSy{\mathrm{\,dS}_{y}}  %
\def\mx{{[m]}}
\def\mxe{{[m*]}}
\def\fx{{[f]}}
\def\rhs{{r.h.s.{~}}}
\def\ext{{\rm{ext}}}
\def\mic{{\rm{mic}}}

\def\aYms#1#2{a_{Y}^{m*} \left ({#1},\,{#2}\right )}

\def\gYm#1#2{g_{Y}^m \left ({#1},\,{#2}\right )}
\def\dYm#1#2{d_{Y}^m \left ({#1},\,{#2}\right )}
\def\aYf#1#2{a_f \left ({#1},\,{#2}\right )}
\def\ipYf#1#2{\left \langle{#1},\,{#2}\right \rangle_{Y_f}}

\newcommand\wrt{{w.r.t.{~}}}
\newcommand\fsi{{\it{fs}}}
\newcommand\ie{{\it{i.e.~}}}
\newcommand\eg{{\it{e.g.~}}}
\newcommand\cf{{{cf.~}}}
\newcommand\wtilde[1]{\widetilde{#1}}
\newcommand\ul[1]{\underline{#1}}
\newcommand\sym{{\rm{sym}}}

\def\what#1{\widehat{#1}}

\def\Vcal{\mathcal{V}}
\def\Wcal{\mathcal{W}}
\def\Mcal{\mathcal{M}}
\def\Ical{\mathcal{I}}

\def\Qcal{\mathcal{Q}}
\def\Lcal{\mathcal{L}}

\def\Ucalbf{\boldsymbol{\mathcal{U}}}
\def\Wcalbf{\boldsymbol{\mathcal{W}}}

\def\IIcalbf{{\mbox{\boldmath$\mathcal{I}$\unboldmath \kern-0.5em\boldmath$\mathcal{I}$\unboldmath}}}
\def\DDcalbf{{\boldsymbol{\mathcal{I}} \kern-0.5em\boldsymbol{\mathcal{D}}}}

\newcommand\GrxyS[2]{\left(\eebx{#1} + \eeby{#2}\right)}

\newcommand{\dist}[2]{{\rm{dist}}\left({#1},{#2}\right)}
\newcommand{\jump}[1]{\left[{#1}\right]}

\newcommand{\Tuf}[1]{{\mathcal{T}}_\veps{\left ({#1}\right )}}
\newcommand\dltsh{\dlt_\tau}

\newcommand\fbav{{\what{\bmi{f}}}}
\newcommand\bbav{{\what{\bmi{b}}}}

\def\ww{{w}}

\def\kkk{{\alpha}}

\def\bb{{\bmi{b}}}
\def\ub{{\bmi{u}}}
\def\vb{{\bmi{v}}}
\def\wb{{\bmi{w}}}
\def\fb{{\bmi{f}}}
\def\nb{{\bmi{n}}}
\def\db{{\bmi{d}}}
\def\eb{{\bmi{e}}}
\def\ssb{{\bmi{s}}}
\def\gb{{\bmi{g}}}
\def\Bb{{\bmi{B}}}
\def\Db{{\bmi{D}}}
\def\Cb{{\bmi{C}}}
\def\Ub{{\bmi{U}}}
\def\Wb{{\bmi{W}}}
\def\Gb{{\bmi{G}}}
\def\Kb{{\bmi{K}}}
\def\Sb{{\bmi{S}}}
\def\Xb{{\bmi{X}}}
 \def\Fb{{\bmi{F}}}
\def\Hb{{\bmi{H}}}
\def\Qb{{\bmi{Q}}}
\def\Ib{{\bmi{I}}}

\def\Pibf{\boldsymbol{\Pi}}
\def\Xibf{\boldsymbol{\Xi}}
\def\Hdb{{\bf{H}}^1}
\def\Hpdb{{\bf{H}}_\#^1}
\def\HpdbO{{\bf{H}}_{\#0}^1}
\def\Hdb{{\bf{H}}^1}

\def\Hpdbav{\wtilde{\bf{H}}_\#^1}

\def\Lb{{\bf{L}}}

\def\Dop{{{\rm I} \kern-0.2em{\rm D}}}
\def\Hop{{{\rm I} \kern-0.2em{\rm H}}}
\def\Fop{{{\rm I} \kern-0.2em{\rm F}}}
\def\Aop{{{\rm A} \kern-0.6em{\rm A}}}%
\def\Iop{{{\rm I} \kern-0.2em{\rm I}}}%

\def\parg{\ul{\gb}}

\def\thetabf{\boldsymbol{\theta}}

\def\eeb#1{\eb({#1})}
\def\eeby#1{\eb_y({#1})}
\def\eebx#1{\eb_x({#1})}
\def\eebz#1{\eb_z({#1})}
\def\R{\hbox{\rm I\kern-0.2em R}}
\def\Z{\hbox{\rm Z\kern-0.3em Z}}

\def\sx{{[s]}}
\def\fx{{[f]}}

\def\circparp#1{{\langle{#1}\rangle_p}}
\def\circpare#1{{\langle{#1}\rangle_\eb}}

\def\MeanYs#1#2{\left \langle{#1}\right \rangle_{#2}}

\def\aYs#1#2{a_{Y}^{s} \left ({#1},\,{#2}\right )}
\def\bYs#1#2{b_{Y}^{s} \left ({#1},\,{#2}\right )}

\def\homlabel{HOM}
\def\reflabel{REF}
\def\hmodel{\homlabel-model\xspace}
\def\rmodel{\reflabel-model\xspace}
\def\prevstep{^{\triangleleft}}

\def\QP{QP}
\def\epot{\bar\varphi}

\begin{document}

\begin{frontmatter}

\title{Homogenized model of peristaltic deformation driven flows in piezoelectric porous media}

\author[NTIS]{E.~Rohan\corref{cor1}}
\ead{rohan@kme.zcu.cz}
\cortext[cor1]{Corresponding author}
\author[NTIS]{V.~Luke\v{s}}
\ead{vlukes@kme.zcu.cz}

\address[NTIS]{Department of Mechanics \&
  NTIS New Technologies for Information Society, Faculty of Applied Sciences, University of West Bohemia in Pilsen, \\
Univerzitn\'\i~22, 30100 Plze\v{n}, Czech Republic}

\begin{abstract}
The paper presents a new type of weakly nonlinear two-scale model of
controllable periodic porous piezoelectric structures saturated by Newtonian
fluids. The flow is propelled by peristaltic deformation of microchannels which
is induced due to piezoelectric segments embedded in the microstructure and
locally actuated by voltage waves. The homogenization is employed to derive a
macroscopic model of the poroelastic medium with effective parameters modified
by piezoelectric properties of the skeleton. To capture the peristaltic
pumping, the nonlinearity associated with deforming configuration must be
respected. In the macroscopic model, \revE{this nonlinearity is introduced through
 homogenized coefficients depending on the deforming micro-configurations}. For this, linear expansions
based on the sensitivity analysis of the homogenized coefficients with respect
to deformation induced by the macroscopic quantities are employed. This enables
to avoid the two-scale tight coupling of the macro- and microproblems otherwise
needed in nonlinear problems. The derived \revE{reduced-order model} is implemented and verified
using direct numerical simulations of the periodic heterogeneous medium.
Numerical results demonstrate the peristaltic driven fluid propulsion in
response to the electric actuation and the efficiency of the proposed treatment
of the nonlinearity. The paper shows new perspectives in homogenization-based
computationally efficient modelling of weakly nonlinear problems where
continuum microstructures are perturbed by coupled fields.
\end{abstract}


\end{frontmatter}



\section{Introduction}

The peristaltic flow is induced by deforming walls of a channel. The study of
this phenomenon is of a great importance in physiology and biomechanics
\cite{Maiti-Mirsa2011,Carew-Pedley1997,Pozrikidis1987,Fung1968}, however, as a
driving mechanism of fluid transport, it presents an important and challenging
issue in the design of smart ``bio-inspired'' materials. We consider locally
periodic porous structures saturated by a Newtonian fluid controllable by
locally distributed piezoelectric segments which can induce peristaltic
deformation wave of the microchannels, in contrast with biological tissues
where other physiological mechanisms, namely muscular contraction provide the
desired actuation.

The intention of the paper is to contribute \revE{to the mathematical modelling and computational methods for newly emerging and perspective
area of controllable and adaptive structures. Since this is a broad and challenging research area, we focus on periodic structures}  equipped with embedded actuators
and sensors connected to electric circuits, which provide a higher level of
multi-functionality such as the capability to harvest energy or to suppress
undesired vibrations. The controllable metamaterials bring many challenging
problems for the field of continuum mechanics and multiscale modelling. Besides
the passive nonlinear mechanical behaviour (geometric nonlinearity due to large
deformations, or displacements, or material nonlinearity due to multiphysics
coupling), the nonlinearity emerges also due to the electric circuit control.
Moreover, by the consequence, such metamaterials can exhibit non-local
properties which, in a sense, break the axiom of neighbourhood. Controlling and
tuning of nonlinearities arising from the fluid-structure interaction coupled
in space and time, being connected with a ``smart'' on-line/real-time control
system provides wide range of applicability in the design of robotic structures
and adaptive shape morphing applications.


Effective material properties of composites involving piezoelectric components
have been studied \revE{and even optimized since a long time, see \eg \cite{NelliSilva1998-cmame}. The research of electro-mechanical transduction in porous piezoelectric structures was motivated by biological tissues \cite{telega1991piezoelectricity,Telega2000}, further  motivated the applications in
activated bio-piezo porous implants \cite{Sikavitsas-etal-2001,Wiesmann-etal-2001}.} In the context of the multiscale modelling, the design of a new type of bio-materials assisting in bone healing
and regeneration processes motivated the homogenization of the periodic
composites consisting of piezoelectric matrix and elastic anisotropic
inclusions \cite{Miara-piezo}, \cf \cite{Rohan-Miara-MAMS-2006} where the
sensitivity of the effective material properties \wrt the inclusion shape was
developed, \revE{\cf \cite{koutsawa_belouettar_makradi_nasser_2010}}. \revE{Piezoelectric fiber composites were investigated \eg in \cite{Medeiros2015,Iyer-IJSS-2014}.}
An asymptotic analysis has been applied to derive higher order
models of piezoelectric rods and beams, starting from the 3D piezoelectricity
problem \cite{Viano-pz-beams-IJSS2016}. Besides the classical periodic
homogenization \cite{Sanchez1980Book,Cioranescu1999book,Cioranescu2008a}, the
Suquet method of homogenization has been used to obtain analytic models of
particulate and fibrous piezoelectric composites \cite{Iyer-IJSS-2014}. The
micromechanics approaches including the Mori-Tanaka and self-consistent
upscaling schemes were employed in \cite{Ayuso-etal-IJSS2017-pz-homog}. Beyond
the linear piezoelectric effect, the issue of composite materials providing an
apparent flexoelectric effect, related to nonuniform deformation, has been
studied using the higher order continua approaches
\cite{Yvonnet_2020,Chen_2021}, \cf \cite{Mawassy_2021}. A theory of controlled
dynamic actuation of flexoelectric membranes interacting with viscous flows was
introduced in \cite{Herrera-Valencia_2014}, a work inspired by biological
systems. \revE{The porous piezoelectric medium treated in the present paper is featured by embedded electrodes enabling to impose locally electric electric field. To compensate this effect in the homogenization limit, the piezoelectric material coefficients scaling must be introduced. By the consequence of the induced ``strong heterogeneity'', one cannot reduce the obtained  homogenized model to the standard homogenization without any scaling which leads to an extended piezoelectric effective constitutive law, as discussed in detail in \cite{Rohan-Lukes-PZ-porel}.}


\paragraph{Novelty and paper contribution}
The paper contributes to the state of the art in \revE{three} aspects:
\textbf{a)} a new type two-scale model of controllable periodic porous
piezoelectric media saturated by a Newtonian fluid is derived using
the homogenization of locally periodic structures; \textbf{b)} we show new perspectives in the homogenization-based computationally efficient modelling of weakly nonlinear problems using the linearization approach based on the sensitivity analysis of the homogenized coefficients \wrt perturbed microconfigurations in response to the local macroscopic fields. \revE{In particular, we demonstrate the efficiency of this modelling approach in the context of the \emph{Model Order Reduction} for two-scale problems within the range of linear kinematics, but non-negligible deformation, such that the equilibrium conditions are imposed in the ``true'' configuration.}
\revE{Moreover, \textbf{c)}} the topic of peristaltic deformation induced in a
fluid-saturated composite medium in response to an electro-mechanical
actuation has not been reported in the literature, namely in the
context of the multiscale modelling.

The paper builds on our recent
works devoted to the three main issues which are synthesized in the
proposed homogenization-based model of the ``smart'' fluid-saturated
porous media describing the peristalsis-driven flow. To account for
the connected porosity, the piezo-poroelastic model
\cite{Rohan-Lukes-PZ-porel} is extended by the Darcy flow model which
is obtained by the Stokes flow homogenization. For studying the
fluid-structure interaction in porous flexoelectric media under large
deformation in the microstructure, the homogenization approach seems
to be perspective \cite{Sandstrom_2016,Lukes_Rohan_2022}, \revE{\cf \cite{Brown2014,Collis2017,Miller2021}}, but
inevitably leads to very complex two-scale computations (the well known ``finite element-square'' (FE$^2$)
approach) requiring reiterated solving the microproblems at all integration nodes of the discretized macroproblem. 
However,  under the restriction to small deformation, the
nonlinearity induced by the deformation-dependent effective model
parameters can be approximated by virtue of the first order expansion formulae derived using the sensitivity analysis
of the homogenized coefficients, following the approach proposed in \cite{Rohan-Lukes-nlBiot2015}. Due to this modelling strategy, the two-scale model is implemented using a computationally efficient scheme avoiding the ``finite element-square'' complexity. The numerical studies reported in the paper demonstrate the necessity of capturing the nonlinearity in the model of the electro-active
porous medium generating the fluid transport.

 \paragraph{Paper organization}
The problem of the fluid interacting with the piezoelectric solid
composite at the heterogeneity level of the periodic microstructure is
introduced in Section~\ref{sec-problem}. The homogenization procedure is
described in Section~\ref{sec-homog}, \revE{where the macroscopic model is derived for the linear problem.}
It involves the homogenized coefficients expressed in terms of the so-called
characteristic responses of microstructures. These are given in terms of
solutions to the micro-problems. The macroscopic model nonlinearity is
introduced in Section~\ref{sec-varHC}. \revE{For this, homogenization of the linearized incremental formulation is explained in \ref{app-incr} using a simplified problem for disconnected porosity. This justifies the use of the perturbation analysis of the
homogenized coefficients depending on the macroscopic state variables for which some technical details are postponed in the \ref{sec-sa}.}
Numerical illustrations of the two-scale modelling of the ``piezo-controlled'' porous
microstructures inducing the peristalsis-driven flow are reported in
Section~\ref{sec:numex}. The pumping effect of the homogenized medium and the
key role of the model nonlinearity are demonstrated using examples of 3D
microstructures with 1D and 2D fluid channel connectivity.
\revE{
  Technical parts related to the model derivations are postponed in the Appendix which consists of four parts: the reference problem is introduced in \ref{sec-DNS}; the modelling framework based on the incremental formulation and approximation using the initial configuration is explained in \ref{app-incr}. Some details on the differentiation of the homogenized coefficients \wrt the microstructure deformation is reported in \ref{sec-sa}.
Finally, in~\ref{sec-1D}, the fluid transport driven by the peristaltic deformation in porous structures is illustrated using a reduced 1D model of the homogenized
medium, showing the importance of the deformation-dependent permeability for the pumping effect.}


\paragraph{Notation}
We employ the following notation. Since we deal with a two-scale problem, we
distinguish the ``macroscopic'' and ``microscopic'' coordinates, $x$ and $y$,
respectively. We use $\nabla_x = (\pd_i^x)$ and $\nabla_y = (\pd_i^y)$ when
differentiation \wrt coordinate $x$ and $y$ is used, respectively, whereby
$\nabla \equiv \nabla_x$. By $\eeb{\ub} = 1/2[(\nabla\ub)^T + \nabla\ub]$ we
denote the strain of a vectorial function $\ub$, where the transpose operator
is indicated by the superscript ${}^T$. The Lebesgue spaces of 2nd-power
integrable functions on an open bounded domain $D\subset \RR^3$ is denoted by
$L^2(D)$, the Sobolev space $\Wb^{1,2}(D)$ of the square integrable
vector-valued functions on $D$ including the 1st order generalized derivative,
is abbreviated by $\Hdb(D)$. Further $\Hpdb(Y)$ is the Sobolev space of
vector-valued Y-periodic functions (the subscript $\#$). \revE{The surface integrals are closed by $\dS,\dSx$, or $\dSy$, while volume integrals are written without an ``integration element'', like $\dd V$.}


\input{aux-Model-pz-porel3.tex}


\input{pz_flow_simulations.tex}


\section{Conclusion}\label{sec-conclusion}

The presented work summarizes our research efforts aimed towards developing
efficient computational tools to simulate behaviour of ``smart porous''
materials controllable by electric field which are designed to transport fluids
using the principle of the peristaltic deformation. In particular, using the
homogenization of the linearized fluid-structure interaction problem, we
derived a model of electroactive porous material. The considered periodic
composite structure involves piezoelectric actuators which provide a
time-and-space control handle to induce a desired deformation wave. In this
respect, the homogenized model (\hmodel) should be used for a given scale
$\vepsz >0$, \ie a given size of the representative periodic cell $\vepsz Y$,
since it is featured by the strong heterogeneity associated with the scaled
material parameters of the micromodel. The \hmodel was validated by means of
the DNS of the heterogeneous structure response taken as ``the reference''
(\rmodel). As expected, the response of the homogenized medium provides a good
approximation for the deformation field, while the approximation of the
pressure is worse, in general. However, it should be noted, that, because of
different treatment of the electrode (equipotential) segments, the prescribed
control -- the electric potential wave, is not quite the same for both the
reference and the \hmodel, so that the comparison is rather affected by the
control wave length. In our validations, to see the fluid transport effect, we
had to keep the wave rather short compared to the number of the considered
microstructure periods.

As confirmed by numerical studies, to achieve the desired pumping effect of the
homogenized continuum, it is necessary to account for the nonlinearity
associated with deformation-dependent microconfigurations and, hence,
respecting deformation-dependent effective properties of the homogenized
material. For this, the sensitivity analysis approach has been applied which
leads to a computationally efficient numerical scheme for solving the nonlinear
problem.

\revE{The main advantage of the derived homogenized model is the reduction of
computational complexity of solving the nonlinear problem. This can be shown by comparison of three possible modelling approaches: a) the DNS (the reference model), b) the FE$^2$ approach -- the homogenization based on updating the spatial local reference microstructures (explained in \ref{app-incr}), and c) the proposed reduced order model based on the homogenized coefficient approximation using the sensitivity of the responses of the undeformed configuration (MOR/SA). In the case a), the reference computation with the FE mesh
constituted by 20 periodic cells leads to a linear system with more than 300000
degrees of freedom, which must be solved in several iterations in each time
step. The solution of each time step takes approximately 2.4 minutes on a personal
computer (10 core Intel i9 processor, 64\,GiB RAM). In the case b), the two-scale calculation
includes solving several microscopic subproblems, evaluating homogenized
coefficients, and calculations at the macroscopic scale using a time stepping
iterative algorithm. With a FE mesh at the microscopic level, consisting of $\approx 16000$ tetrahedral elements and $\approx 3000$ nodes. For illustration, solving a single micro problem without sensitivity analysis takes approximately 8s, which for 160 macroscopic integration points means the solution time of approximately 1280s for the microproblems at each macroscopic iteration step.
Finally, in the case c), the ``MOR/SA'' approach requires solving the ``cell problem'' and the sensitivity analysis only once, prior to solving the macroscopic problem (independently on the number of iterations and time steps). This about 35s and the evaluation of the approximated coefficients is performed in tens of milliseconds (point-wise at the macroscopic integration points). One may conclude that the computational reduction of the MOR/SA method compared to the FE$^2$ method is by the factor of 1/30 to 1/40. Expressed in the wall-clock time, the numerical simulations reported in Section 5.3 using the proposed MOR/SA approach is completed in about 5 minutes, while the solution times for both the DNS and FE$^2$ method for the same problems require  more than 2 hours. }



The two-scale piezoelectric model respecting the geometrical nonlinearity seems
to be an effective tool for simulating peristaltic flows in porous
piezoelectric structures. Further studies and model refinement are envisaged
towards using the acoustic wave energy for the fluid transport, in analogy with
the ``acoustic streaming'' phenomenon. For this, the influence of the flow
dynamics and respecting of all the inertia related phenomena will be important
for higher frequencies of the voltage actuation through the piezoelectric
actuators.

\noindent \paragraph{Acknowledgement}
The research has been supported by the grant projects GA~21-16406S and GA~24-12291S
 of the Czech Science Foundation.

\input{aux-Appendix-pz-porel-short.tex} 

\input{aux-1Dmodel-R.tex}

\bibliographystyle{elsarticle-harv}
\bibliography{Biblio-PZ-flow}

\end{document}

%% file: aux-Model-pz-porel3.tex
\section{Micromodel and the fluid-structure interaction}\label{sec-problem}

In the geometrical framework describing the fluid saturated porous medium at
the heterogeneity level, we introduce the fluid-structure interaction (FSI) problem
for which the two-scale problem is obtained by the homogenization method in
Section~\ref{sec-homog}.
\chE{Although in this paper we consider small deformations of the solid and the linear constitutive law \eq{eq-1}, as will be shown, to capture the ``pumping effect'' of the peristaltic deformation, it is necessary to consider the nonlinearity of the FSI on the interface due to the deforming pores, thereby modifying the hydraulic permeability. A brief illustration of this phenomenon using a 1D reduced model is presented in \ref{sec-1D}.}

\subsection{\newE{Linearization based on the Eulerian frame}}
The FSI problem can be formulated in the spatial configuration $\Om\subset \RR^3$ and further approximated using a reference configuration $\Om_R\subset \RR^3$ which enables to introduce an incremental formulation involving the stress increment defined in terms of the Lie derivative  $\dlt_R^\Lcal \taubf_R$. The equilibrium for both the solid and fluid can be expressed by the virtual power law where the virtual (test) displacement field $\tilde \vb$ is defined in $\Om_R = \Om_R^s \cup \Om_R^f$, being continuous on the interface $\Gamma_R^\fsi = \pd\Om_R^s \cup \pd\Om_R^f$,
\begin{equation}\label{eq-fsi1}
  \begin{split}
    \int_{\Om_R} \dlt_R^\Lcal \taubf_R:\eb_R(\tilde \vb) J_R^{-1} \dd x + \int_{\Om_R}  \taubf_R:\dlt\etabf_R(\tilde \vb) J_R^{-1} \dd x \\
    = \int_{\pd_\sigma \Om_R^f} \gb_R(\bar p)\cdot \tilde \vb \dd S - \int_{\Om_R}\taubf_R:\eb_R(\tilde \vb) J_R^{-1} \dd x\;,
  \end{split}
\end{equation}
to be satisfied for all test displacements $\tilde \vb$ which are considered to be smoothly extended to the whole domain $\Om_R$. 

Further we neglect the variation the nonlinear part $\dlt\etabf_R(\tilde \vb)$ of the Green Lagrange strain associated with the displacement increment $\dlt\ub$. Correspondingly, we shall assume $J_R \approx 1$ by virtue of the small deformation (so that we do not distinguish between the  Cauchy and the Kirchhoff stresses), hence \eq{eq-fsi1} becomes
\begin{equation}\label{eq-fsi1a}
  \begin{split}
    \int_{\Om_R} \dlt_R^\Lcal \taubf_R:\eb_R(\tilde \vb) \dd x + \int_{\Om_R}\taubf_R:\eb_R(\tilde \vb)  \dd x = \int_{\pd_\sigma \Om_R^f} \gb_R(\bar p + \dlt \bar p)\cdot \tilde \vb \dd S \;.
  \end{split}
\end{equation}
By the consequence, assuming the equilibrium is attained in $\Om_R$ for the stress $\taubf_R$ and external loads due to $\bar p$ at $\pd_\sigma \Om_R^f$, \eq{eq-fsi1a} would further be simplified, yielding
\begin{equation}\label{eq-fsi1b}
  \begin{split}
    \int_{\Om_R} \dlt_R^\Lcal \taubf_R:\eb_R(\tilde \vb) \dd x 
    = \int_{\pd_\sigma \Om_R^f} \gb_R(\dlt \bar p)\cdot \tilde \vb \dd S \;.
  \end{split}
\end{equation}
Equation \eq{eq-fsi1a} can be presented in the form of \eq{eq-fsi3}. More comments on solving the nonlinear problem using the are postponed in \ref{sec-DNS}. In what follows we use the displacement $\ub$ as being extended in the whole reference configuration $\Om_R$, so being defined also in the fluid $\Om_R^f$.

\begin{figure}
\centering 
  \includegraphics[width=0.7\linewidth]{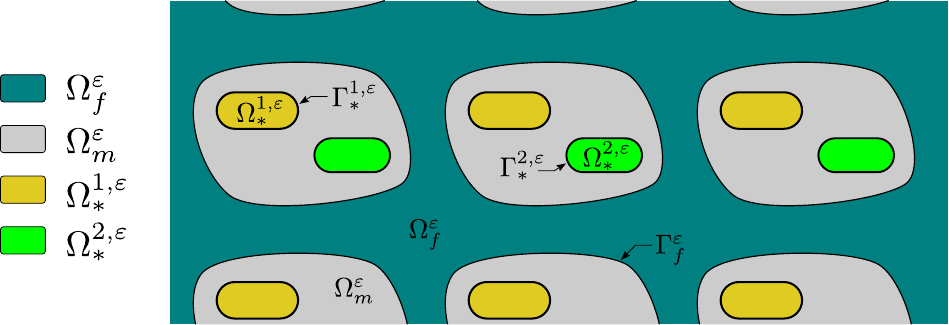}
\caption{\revE{The periodic microstructure and its decomposition into subdomains, as described in Eqs.~\eq{eq-pzm1} and \eq{eq-pzm1b}.}}
  \label{fig-micro-scheme}
\end{figure}

\subsection{Porous piezoelectric solid saturated by static fluid}\label{sec:fspm}  

We consider a quasi-static loading of a piezoelectric skeleton interacting with
a viscous fluid saturating pores in the skeleton.
The piezo-poroelastic medium occupies an open bounded domain $\Om \subset \RR^3$
with Lipschitz boundary $\pd \Om$.
The following decomposition of $\Om $ into the piezoelectric matrix,
$\Om_m^\veps$, elastic conductive inclusions, $\Om_*^\veps$, and 
the fluid occupying the channels, $\Om_f^\veps$, is considered:
\begin{equation}\label{eq-pzm1}
\begin{split}
\Om = \Om_f^\veps \cup \Om_m^\veps \cup \Om_*^\veps\cup \Gamma^\veps\;,
\quad \Om_f^\veps \cap \Om_m^\veps \cap \Om_*^\veps = \emptyset\;,\\
\Om_s \equiv \Om_{m*} = \Om_m^\veps \cup \Om_*^\veps \cup \Gamma_*^{\veps}\;,\quad
\mbox{ where } \Om_*^\veps = \bigcup_\kkk \Om_*^{\kkk,\veps}\;,
\end{split}
\end{equation}
where $\Gamma^\veps$ is the union of all interfaces separating the main three parts.
Using the notation $\Om_s$  we often refer to the solid part consisting of the piezoelectric matrix and
the conductor subparts which are separated by $\Gamma_*^{\veps}$; this interface consists of its
 segments $\Gamma_*^{\kkk,\veps}$ introduced, as follows:
\begin{equation}\label{eq-pzm1b}
  \Gamma_*^{\veps} = \bigcup_\kkk\Gamma_*^{\kkk,\veps} \quad  \mbox{ and }\quad \Gamma_*^{\kkk,\veps} = \pd \Om_*^{\kkk,\veps} \cap \ol{\Om_m^\veps}\;.
\end{equation}

In \eq{eq-pzm1}, each $\Om_*^{\kkk,\veps}$ is assumed to be a union of segments with a diameter $\approx \veps$; this restriction prevents simply connected conducting parts. Instead,  theses small segments are employed as the control handles allowing for a local control of the electric field due to (locally) prescribed voltage, as will be specified below.
%
The decomposition of the periodic microstructure is displayed in Fig.~\ref{fig-micro-scheme}.
By $\Gamma_f^\veps$ we denote the solid-fluid interface, $\Gamma_f^\veps  =
\ol{\Om_m^\veps \cup \Om_*^\veps} \cap \ol{\Om_f^\veps}$. The interface between
the piezoelectric matrix and the conductors $\Gamma_*^{\veps}$ 
Further we denote by $\pd_\ext\Om_m^\veps = \ol{\Om_m^\veps} \cap \pd\Om$ the
exterior boundaries of $\Om_m^\veps$. In analogy, we define
$\pd_\ext\Om_f^\veps$ as the exterior parts on the
boundaries of the fluid $\Om_f^\veps$.   We assume that both the matrix $\Om_m^\veps$ and the fluid-filled channels $\Om_f^\veps$ are connected domains. 

 The boundary conditions for both the solid and fluid phases are prescribed on
the external boundary segments $\pd_\ext\Om_{m*} = \pd\Om_{m*}\cap\pd\Om$ and $\pd_\ext\Om_f = \pd\Om_f\cap\pd\Om$.
The following two splits are defined,
$\pd_\ext\Om_{m*}^\veps = \Gamma_{\ub}^\veps \cup \Gamma_\sigma^\veps$ and
$\pd_\ext\Om_f^\veps = \Gamma_p^\veps \cup \Gamma_w^\veps$ 
such that the ``Dirichlet'' type  boundary conditions are prescribed on $\Gamma_{\ub}^\veps$ and $\Gamma_p^\veps$, whereas the surface traction forces and fluid fluxes are prescribed on $\Gamma_\sigma^\veps$ and $\Gamma_w^\veps$, respectively, where
\begin{equation}\label{eq-pzm1a}
\Gamma_\sigma^\veps = \pd_\ext\Om_{m*}^\veps\setminus\Gamma_{\ub}^\veps\;, \quad 
\Gamma_w = \pd_\ext\Om_f^\veps\setminus\Gamma_p^\veps\;.
\end{equation}

To respect spatial fluctuations of the material parameters, by virtue of the
scale parameter introduced above, all material coefficients and unknown
functions involved in the mathematical model which depend on the scale will be
labelled by superscript $\veps$. In the piezoelectric solid, the Cauchy stress
tensor $\sigmabf^\veps$ and the electric displacement $\vec D^\veps$ depend on
the strain tensor $\eeb{\ub^\veps} = (\nabla\ub^\veps + (\nabla\ub^\veps)^T)/2$
defined in terms of the displacement field $\ub^\veps = (u_i^\veps)$, and on the
electric field $\vec E(\vphi^\veps) = \nabla \vphi^\veps$ defined in terms of the
electric potential $\vphi^\veps$. The following constitutive equations
characterize the piezoelectric solid in $\Om_m^\veps$,
\begin{equation}\label{eq-1}
\begin{split}
\sigma_{ij}^\veps(\ub^\veps,\vphi^\veps) & = A_{ijkl}^\veps e_{kl}^\veps(\ub^\veps) - g_{kij}^\veps E_k^\veps(\vphi^\veps)\;,\\
D_k^\veps(\ub^\veps,\vphi^\veps) & =  g_{kij}^\veps e_{ij}^\veps(\ub^\veps) + d_{kl}^\veps E_l^\veps(\vphi^\veps)\;,
\end{split}
\end{equation}
where $\Aop^\veps=(A_{ijkl}^\veps)$ is the elasticity fourth-order symmetric
positive definite tensor of the solid, $A_{ijkl} = A_{klij} = A_{jilk}$. The deformation is coupled with the electric field through the 3rd
order tensor $\parg^\veps = (g_{kij}^\veps)$, $g_{kij}^\veps = g_{kji}^\veps$ and
$\db = (d_{kl})$ is the permittivity tensor. The conductive solid is described by
its elasticity $\Aop^\veps$ only. The permittivity in a conductor is
infinitely large, so that we assume a constant potential $\vphi^\veps = \bar\vphi^\kkk$ in any $\Om_*^{\kkk,\veps}$, whereby
\eq{eq-1} reduces to the elasticity law.

In this paper we shall use a compact (global) notation, such that \eq{eq-1} can
be written (we drop the superscript $^\veps$ for the moment),
$\sigmabf_s = \Aop\eeb{\ub} - \parg^T\cdot\vec E(\vphi)$, and
$\vec D = \parg:\eeb{\ub} + \db\vec E$, where $(\parg:\eb)_k = g_{kij} e_{ij}$. We shall not distinguish piezoelectric and dielectric (non-piezo electric) materials which both are governed by \eq{eq-1} with the only difference of vanishing $\parg$ in the latter case. In the conductor, \eq{eq-1}$_2$ is not considered.

Concerning the fluid flow described by the pressure $p^\veps$ and the velocity $\vb^{f,\veps}$ subject to the nonslip  boundary conditions on the pore walls, the fluid viscosity must be scaled. 
We consider Newtonian fluids, such that the viscous stress rheology is given by
$D_{klij}^{f,\veps} = \mu^\veps (\dlt_{ik}\dlt_{jl} + \dlt_{il}\dlt_{jk} - (2/3)\dlt_{ij}\dlt_{kl})$, where  $\mu^\veps$ is the dynamic viscosity. 
Hence, the stress in fluid is given by
\begin{equation}\label{eq-S*2}
    \sigmabf_f^\veps = -p^\veps \Ib + \Dop^{f,\veps} \eeb{\vb^{f,\veps}}\;, \quad
    \mbox{ where }\quad \Dop^{f,\veps} = 2\mu^\veps(\Iop - \frac{1}{3}\Ib\otimes\Ib)\;.
\end{equation}

\subsection{Problem formulation for the heterogeneous medium}\label{sec:prob}
\revE{By virtue of the incremental formulation introduced above and in \ref{sec-DNS}, in this section, we concern the linear subproblem imposed in a fixed configuration. The unknown fields $\ub^\veps,\vb^{f,\veps}, p^\veps,\vphi^\veps$ can be considered as the perturbations, in the context of solving the nonlinear problem, see \ref{app-incr}.}


\subsubsection{Strong formulation (with inertia terms)}
We shall now consider all inertia effects which are relevant to state the problem. Later on, due to restrictions to quasi-static events only, these inertia effects will be disregarded to derive the homogenized model.
The micromodel involves the following differential equations governing the fluid-solid interaction and the electric field coupled with the deformation through the piezoelectric constitutive law \eq{eq-1},
\begin{equation}\label{eq-pzfl0}
\begin{split}
\rho_s\ddot\ub^\veps -\nabla\cdot \sigmabf_s^\veps(\ub^\veps,\vphi^\veps) & = \fb^{s,\veps}\;, \quad \mbox{ in }  \Om_{m*}^\veps\;, \\
-\nabla\cdot \vec D^\veps(\ub^\veps,\vphi^\veps) & = q_E^\veps\;,\quad \mbox{ in }  \Om_m^\veps\;,\\
\mu^\veps \nabla^2 \vb^{f,\veps}  - \rho_f \left(\dot\vb^{f,\veps}
+ 
\wb^\veps \cdot\nabla)\vb^{f,\veps}
\right) & = \nabla p^\veps - \fb^{f,\veps}\;,\quad \mbox{ in } \Om_f^\veps\;,\\
\quad \gamma \dot p + \nabla\cdot\vb^{f,\veps} & = 0\;,\quad \mbox{ in } \Om_f^\veps\;,
\end{split}
\end{equation}
where $\gamma$ is the fluid compressibility.

In this respect, it is convenient  to use the fluid velocity decomposition 
\begin{equation}\label{eq-FS10}
  \begin{split}
    \wb^\veps = \vb^{f,\veps} - \dot{\wtilde\ub}^\veps\;,\\
  \end{split}
\end{equation}
where $\dot{\wtilde\ub}^\veps$ is the solid phase velocity field extended from $\Om_m^\veps$ to pores $\Om_f^\veps$, and $\wb^\veps = \zerobf$ on the pore walls $\Gamma_\fsi^\veps$.
\revE{Note that $\dot\ub$ is the solid velocity (``dot'' means the material derivative) and $\tilde\ub$ denotes a smooth extension of $\ub$ from the solid to the pores. 
}

\revE{In general, the boundary conditions prescribed for the solid involve the Dirichlet conditions concerning the fields $\ub^\veps,\vb^{f,\veps}$ and $\vphi^\veps$, whereas the Neumann type conditions comprise the boundary tractions $\sigmabf_s^\veps\cdot \nb = \bb^{s,\veps}$ acting on $\Gamma_\sigma$ of the solid surface, the fluid pressure $p^\veps$ representing the stress in the fluid acting on $\Gamma_p^\veps$ (assuming the dissipative part of $\sigmabf_f^\veps\cdot \nb$ vanishes) and the electric insulation, such that}
\revE{
\begin{equation}
  \begin{split}\label{eq-BC1}
    \ub^\veps  = \ub_\pd \quad \mbox{ on } \Gamma_u^\veps\;,& \quad \sigmabf_s^\veps\cdot \nb = \bb^{s,\veps}\mbox{ on } \Gamma_\sigma\;, \\
       \vb^{f,\veps}  = 0  \quad \mbox{ on } \Gamma_w^\veps\;,& \quad  p^\veps  = p_\pd \quad \mbox{ on } \Gamma_p^\veps\;, \\
    \vphi^\veps  = \bar \vphi_\veps^\kkk\;,\quad \mbox{ on } \Gamma_*^{\kkk,\veps}\;,\;\kkk = 1,\dots,\bar\kkk\;,& \quad \vec D^\veps\cdot\nb = 0\quad \mbox{ on } \pd\Om_m^\veps\setminus \Gamma_*^\veps\;.
    \end{split}
\end{equation}
}
It should be noted that the control electric potential $\bar \vphi_\veps^\kkk(t,x)$ can be defined for all $x \in \Om$ attaining a constant value over each simply connected domain $\Om_*^{\kkk,\veps}$ at a given time $t$. We assume $\bar \vphi_\veps^\kkk(t,\cdot)\in C^0(\Om)$, is such that it converges strongly 
to some  $\bar \vphi^\kkk(x,t)$.

The following interface conditions express the continuity of the traction stress, velocity and electric potential,
\begin{equation}
\begin{split}\label{eq-BC2}
  \nb \cdot \jump{\sigmabf^\veps}  = 0\;,\quad 
\jump{\ub^\veps}  = 0\;,\quad \mbox{ on } \Gamma_*^\veps\;,\quad \\
\vb^{f,\veps}  = \dot\ub^\veps\;,\quad 
    \nb\cdot\sigmabf_f^\veps  = \nb\cdot\sigmabf_s^\veps\;,  \mbox{ and }\quad 
    \nb\cdot \vec D^\veps  = 0\quad\mbox{ on } \Gamma_\fsi^\veps\;.
\end{split}
\end{equation}

In the weak formulations introduced below, the following admissibility sets and spaces will be employed, being defined using the boundary conditions \eq{eq-BC1},
\begin{equation}\label{eq-5a}
\begin{split}
\Ucalbf^\veps(\Om_{m*}^\veps) & = \{\vb \in \Hdb(\Om_{m*}^\veps)|\; \vb = \ub_\pd \mbox{ on }  \Gamma_{\ub}^\veps\}\;,\\
\Vcal_*(\Om_m^\veps,\{\Gamma_*^{\kkk,\veps}\}_\kkk) & = \{\vphi \in H^1(\Om_m^\veps)|\,\vphi = \bar\vphi_\veps^\kkk
\mbox{ on } \Gamma_*^{\kkk,\veps}\}\;,\\
\Qcal^\veps(\Om_f^\veps)  & = \{q \in H^1(\Om_f^\veps) |\; q = p_\pd \mbox{ on }  \Gamma_p^\veps\}\;,\\
 \Wcalbf_0^\veps(\Om_f^\veps,\Gamma_\fsi^\veps) & = \{\vb \in \Hdb(\Om_f^\veps)|\;\vb = \zerobf \mbox{ on } \Gamma_\fsi^\veps\cup \Gamma_{w}^\veps\}\;,
\end{split}
\end{equation}
where the space $\Wcalbf_0^\veps(\Om_f^\veps,\Gamma_\fsi^\veps)$ will be employed for the seepage velocity $\wb^\veps$.
%
By $\Ucalbf_0^\veps(\Om_{m*}^\veps)$ and $\Qcal_0^\veps(\Om_f^\veps)$ we denote the spaces associated with the affine sets  $\Ucalbf^\veps(\Om_{m*}^\veps)$ and $\Qcal^\veps(\Om_f^\veps)$ defined above. In analogy, $\Vcal_0^\veps(\Om_m^\veps,\{\Gamma_*^{\kkk,\veps}\}_\kkk)$ is derived from $\Vcal_*$ by taking there $\bar\vphi_\veps^\kkk := 0$.


\subsubsection{Weak formulation for the quasistatic flow (inertia effects neglected)}
\revE{Further in this paper, we disregard the inertia effects in \eq{eq-pzfl0}, while restricting the model to describe events characterized by slow flows and deformation. By the consequence, the homogenization using the periodic unfolding method is applied to the weak formulation of the quasistatic fluid-solid interaction which is obtained  upon multiplying the respective equations by test functions and integrating by parts.} Then,  upon substituting  the constitutive relationships, the quasistatic restriction yields the following problem for
$(\ub^\veps(t,\cdot),\vphi^\veps(t,\cdot),\wb^\veps(t,\cdot),\bar p^\veps(t,\cdot)) \in \Ucalbf^\veps(\Om_{m*}^\veps) \times\Vcal_*(\Om_m^\veps,\{\Gamma_*^{\kkk,\veps}\}_\kkk)\times \Wcalbf_0^\veps(\Om_f^\veps,\Gamma_\fsi^\veps)\times\Qcal(\Om_f^\veps)$ satisfying at $t>0$
\begin{equation}\label{eq-pzfl2}
\begin{split}
  \int_{\Om_{m*}^\veps} \left(\Aop^\veps\eeb{\ub^\veps} -\parg^\veps \nabla \vphi^\veps
  \right): \eeb{\vb^\veps} \dV& \\
  - \int_{\Gamma_\fsi^\veps} \nb^\veps\otimes\vb^\veps: \sigmabf_f^\veps(\wb^\veps,\dot{\tilde\ub},p^\veps) \dS & =\int_{\Om_{m*}^\veps} \fb^{s,\veps}\cdot\vb^\veps \dV + 
  \int_{\pd_\sigma\Om_{m*}^\veps} \bb^{s,\veps}\cdot\vb^\veps \dV\;,\\
    \int_{\Om_{m}^\veps} \left((\parg^\veps)^T:\eeb{\ub^\veps} \dV
     + \db^\veps \nabla \vphi^\veps\right)\cdot \nabla \psi^\veps  \dV
      & =  \int_{\Om_{m}^\veps} \rho_E^\veps \psi^\veps \dV\;,\\
  \int_{\Om_{f}^\veps} \left(2\mu^\veps\eeb{\wb^\veps + \dot{\tilde\ub}^\veps} \dV
   - p^\veps\Ib\right) : \eeb{\thetabf^\veps} \dV
   & = \int_{\Om_{f}^\veps} \fb^{f,\veps}\cdot\thetabf^\veps \dV \;,\\
  \int_{\Om_{f}^\veps} q^\veps\left( \gamma \dot p^\veps + \nabla\cdot(\wb^\veps + \dot{\tilde\ub})\right) \dV & = 0\;,\\
\end{split}
\end{equation}
for all  $(\vb^\veps,\psi^\veps,\thetabf^\veps,q^\veps) \in  \Ucalbf_0^\veps(\Om_{m*}^\veps) \times\Vcal_0(\Om_m^\veps,\{\Gamma_*^{\kkk,\veps}\}_\kkk)\times 
\Wcalbf_0^\veps(\Om_f^\veps,\Gamma_\fsi^\veps)\times\Qcal_0^\veps(\Om_f^\veps)$


\revE{
\begin{remark}\label{rem-bi-fsi}
Formulation \eq{eq-fsi3} presented in \ref{sec-DNS} concerning the direct numerical simulations in the heterogeneous structures is obtained from \eq{eq-pzfl2} by virtue of the interaction integral which can be rewritten using the volume integral (note the normal orientation $\nb^\sx = -\nb^\fx$), so that (we drop ${}^\veps$, assuming a given scale $\veps_0$, $\Om_R^f = \Om^f$ in the context of the reference configuration of the iterative algorithm)
\begin{equation}\label{eq-dns5}
  \begin{split}
    -\int_{\Gamma_\fsi}\sigmabf_f:\nb^\sx\otimes\wtilde\vb \dS & = 
     \int_{\Om_f} \sigmabf_f: \nabla\wtilde\vb \dV -\int_{\pd_\sigma\Om}\sigmabf_f:\nb^\fx\otimes\wtilde\vb\dV\;,
  \end{split}
\end{equation}
where the equilibrium identity $\nabla\cdot\sigmabf_f = 0$ in $\Om_f$ has been employed.
  \end{remark}
}

\section{Homogenization of the fluid-structure interaction problem}\label{sec-homog}
\chE{Although we consider upscaling of a linear problem, in \ref{app-incr}, we show that it can be
interpreted in the context of solving a nonlinear problem using a sequence of
linear subproblems associated with an incremental formulation.}

We build on the results obtained in \cite{Rohan-Lukes-PZ-porel}, where disconnected porosity formed by the union of fluid-filled inclusions of size $\sim \veps$ and  a static, or quasi-static loading without any flow effect in these separated inclusions was considered. 
Even for the connected porosity allowing for the nonstationary Stokes flow, it can be shown that
the homogenized parameters of the effective piezo-poroelastic constitutive law
can be obtained as in the static, or quasistatic case independently on the
time. In \cite{Rohan-zamp2020}, this property of decoupled homogenization of
the solid and the fluid was obtained even in the fully dynamic case when the
inertia effects are respected in terms of an acoustic linearization. This applies also in our present situation involving the piezoelectric skeleton -- the homogenized flow obeys a Darcy law  described in terms of a permeability tensor being computed independently of the deformation and the electric field.

\subsection{Periodic geometry and the representative cell}\label{sec:per}

\begin{figure}
  \centering
  \includegraphics[width=0.9\linewidth]{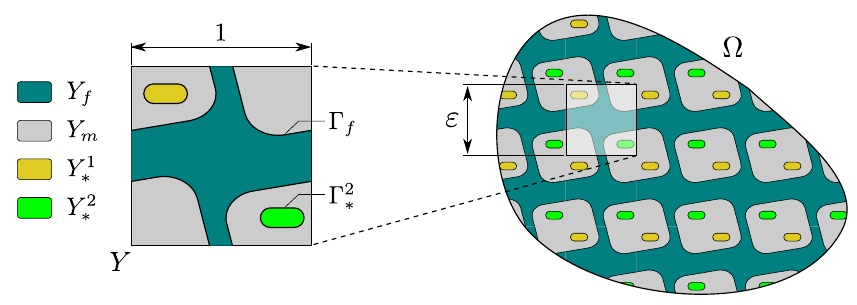}
  \caption{The scheme of the representative periodic cell decomposition and the
  generated periodic structure. The porous microstructure
  with conducting parts $Y_*^\alpha$ embedded in the solid piezoelectric skeleton
  $Y_m$ and fluid part $Y_f$. \revE{$\Om_k^\veps$ is generated as a periodic lattice by $\veps Y_k$,  $k = m,*,f$. (Note the volume $|\veps Y| \sim \veps^3$ in 3D.)}}\label{fig-CFDF}
\end{figure}

\revE{The fluid saturated porous medium is generated as a periodic lattice by a reference periodic cell $Y$
decomposed into three non-overlapping subdomains $Y_m$, $Y_f$ and $Y_*$, see
Fig.~\ref{fig-CFDF}, whereby $|Y|\approx 1$. Specifically, $\Om_k^\veps$ is generated by repeating the rescaled cell subdomain $\veps Y_k$ within the lattice, where $k = m,*,f$.
Thus, in accordance with the decomposition \eq{eq-pzm1} }
\begin{equation}\label{eq-6}
\begin{split}
  Y = Y_m \cup Y_f \cup Y_* \cup \Gamma_Y\;, \quad
   Y_i \cap Y_j = \emptyset\, \mbox{ for } j\not = i \mbox{ with } i,j \in \{f,m,*\}\;,\\
   Y_* =\bigcup_\alpha Y_*^\alpha\;,\quad Y_*^\alpha\cap Y_*^{\beta} = \emptyset\, \mbox{ for } \alpha\not = \beta\;, \quad
   Y_s \equiv Y_{m*} = Y_m\cap Y_*\;,
\end{split}
\end{equation}
where \revE{$\dist{Y_*^\alpha}{Y_*^\beta} \geq s^*(1-\delta_{\alpha\beta})$ with $s^*>0$ being the minimum distance between different conductive parts.} $\Gamma_Y$,
representing the union of all the interfaces, splits in two disjoint parts (in
the sense of the surface measure), assuming $\pd Y_f \cap \pd Y_* = \emptyset$ (the conductors separated from the fluid)
\begin{equation}\label{eq-6a}
\Gamma_Y = \Gamma_\fsi \cup \Gamma_{m*}\;, \quad \Gamma_\fsi  = \ol{Y_s} \cap \ol{Y_f}\;,\quad 
\Gamma_{m*}^\alpha = \ol{Y_*^\alpha} \cap \ol{Y_m}\;,\quad
\Gamma_{m*} = \bigcup_\alpha \Gamma_{m*}^\alpha \;.
\end{equation}
We also assume well separated conductor parts, \ie $\ol{Y_*^\alpha}\cap \pd Y = \emptyset$, such that the periodic lattice $\Om_*^\alpha$ is constituted by mutually disconnected inclusions with the perimeter $\approx \veps$.

\subsection{Material parameters and scaling}\label{sec:mat}

Specific features of the micromodel can be retained in the limit
$\veps\rightarrow 0$, if some of the involved material parameters are
considered to depend on $\veps$. To allow for microstructures with strongly
controlled electric field, in formulation \eq{eq-pzfl2}, potentials
$\bar\vphi^\alpha$ are given for each simply connected domain
$\Om_*^{\alpha,\veps}$ occupied by the perfect conductor and represented by
$Y_*^\alpha$ within the cell $Y$. The following assumptions are made:
%
\begin{quote}
\begin{tabular}{lll}
  \textbf{a)}  Strongly controlled field: & $\bar\vphi^{\veps} = \bar\vphi^{\alpha}$ & in $\Om_*^{\alpha,\veps}$\;, \\
  \textbf{b)} Weakly piezoelectric material: & $\parg^\veps = \veps \bar\parg\;,\quad  \db^\veps  = \veps^2 \bar\db$ & in $\Om_m^\veps$\;,\\
  \textbf{c)} Nonslip condition on pore walls: & $\mu^\veps  = \veps^2 \bar\mu$ & in $\Om_f^\veps$\;.
\end{tabular}
\end{quote}

Consistently with the assumption \textbf{a)}, since $\bar\vphi^\alpha$ does not
vanish with $\veps\rightarrow 0$, steep gradients on the electric potential are
assumed for small $\veps$, which must be compensated by the scaling
\textbf{b)}, \revE{as discussed in \cite{Rohan-Lukes-PZ-porel}}. In the fluid part, the scaling \textbf{c)} of the viscosity by
factor $\veps^2$ is coherent with the quasistatic flow model with the nonslip
conditions for the fluid velocity at $\Gamma_f^\veps$, see \eg
\cite{Allaire-nonsteady-NS-homog,Hornung1997book,Clopeau-Mikelic-etal-MCM2001}.
This holds for problem \eq{eq-pzfl2}, when the advection acceleration can be
neglected. However, if fluid inertia is concerned, $\mu^\veps = \veps^\beta
\bar\mu$ with a constant $\bar\mu>0$, whereby exponent $\beta > 0$ is to be
determined, \cf \cite{Chen-homog-NS-Forchheimer2001}.

\subsection{Asymptotic expansions}

The following truncated expansions (the so-called recovery sequences) of
displacements, electric potential, pressure and the fluid seepage velocity can be introduced formally such that the limit equations can be derived:
\begin{equation}\label{eq-FS14a}
\begin{split}
\Tuf{\ub^{\veps}} \approx \ub^{R\veps}(x,y) & := \ub^{0\veps}(x) + \veps \ub^{1\veps}(x,y)\;,  \\
\Tuf{\vphi^\veps}\approx\vphi^{R\veps}(x,y) & := \hat\vphi^{0\veps}(x,y) \;,\\
\Tuf{p^\veps}\approx p^{R\veps}(x,y) & := p^{0\veps}(x) + \veps p^{1\veps}(x,y) \;,\\
\Tuf{\wb^\veps} \approx \wb^{R\veps}(x,y) & := \hat\wb(x,y)\;.
\end{split}
\end{equation}
The test functions 
are considered in the analogous form, thus
\begin{equation}\label{eq-FS14c}
\begin{split}
  \Tuf{\vb^\veps(x)} & = \vb^0(x) + \veps\vb^1(x,y)\;,\\
\Tuf{\psi^\veps}& = \hat\psi^{0\veps}(x,y) \;,\\
  \Tuf{q^\veps(x)} & = q^0(x) + \veps q^1(x,y)\;,\\
    \Tuf{\thetabf^\veps(x)} & = \hat\thetabf(x,y)\;.
\end{split}
\end{equation}
All the two-scale functions are $Y$-periodic in the second variable $y$ and for almost all $t>0$
\begin{equation}\label{eq-FS14d}
\begin{split}
  \ub^1(\cdot,t), \vb^1 & \in L^2(\Om; \Hpdb(Y_s))\;,\quad \ub^0(\cdot,t) \in \Hdb(\Om)\;,\\
\hat\wb(\cdot,t),\hat\thetabf & \in L^2(\Om;\HpdbO(Y_f))\;,\\
p^1(\cdot,t),q^1 & \in L^2(\Om; H_\#^1(Y_f))\;,\quad p^0(\cdot,t) \in H^1(\Om)\;, \\
\hat\vphi^0(x,\cdot,t) & \in L^2(\Om; H_{\#0*}^1(Y_m)) + \sum_\alpha \Phi^\alpha \bar\vphi^\alpha(\cdot,t)\;,
\end{split}
\end{equation}
where 
\begin{equation}\label{eq-Hsp}
  \begin{split}
    \HpdbO(Y_f) & = \{\vb \in \Hpdbav(Y_f)|\; \vb = \bmi{0}  \mbox{ on } \Gamma_{fs}\}\;,\\
    H_{\#0*}^1(Y_m) & = \{\psi \in H_\#^1(Y_m)|\;\psi = 0 \mbox{ on } \Gamma_*\}\;,\\
    H_{\#0,\alpha}^1(Y_m) & =  H_{\#0*}^1(Y_m) + \Phi^\alpha \;, \quad
    \Phi^\alpha = \delta_{\alpha\beta} \mbox{ on } \Gamma_*^\beta\;,
\end{split}
\end{equation}
with $\Phi^\alpha \in H_\#^1(Y_m)$ otherwise given arbitrarily.   
Note that the nonslip condition of $\hat\wb$ on  $\Gamma_\fsi$ is imposed due to the  space $\HpdbO(Y_f)$.
As the result, $\int_{Y_f} \nabla_y\cdot \hat\thetabf = 0$ for any $\hat\thetabf \in \HpdbO(Y_f)$.

Further we shall employ the volume fractions $\phi_f:= |Y_f|/|Y|$ and $\phi_s = 1-\phi_f$. The notion of the surface phase fractions $\bar\phi_s$ and $\bar\phi_f$ is needed as well, being introduced by virtue of the following convergences,
\begin{equation}\label{eq-FS20a}
    \int_{\Gamma_\sigma^\veps}\vtheta\dS \rightarrow \int_{\pd_\sigma\Om} \bar\phi_s  \vtheta\dS
    = \int_{\pd\Om} \bar\phi_s  \vtheta\dS\;, \quad
    \int_{\Gamma_p^\veps}\vtheta\dS \rightarrow \int_{\pd_p\Om} \bar\phi_f  \vtheta\dS
    = \int_{\pd\Om} \bar\phi_f  \vtheta\dS\;, 
\end{equation}
for all $\vtheta \in H^1(\Om)$, where the equalities hold due to defining $\bar\phi_s := 0$ on $\pd\Om\setminus {\pd_\sigma\Om}$ and $\bar\phi_f := 0$ on $\pd\Om\setminus {\pd_p\Om}$. This will enable to introduce the effective traction stress, see below \eq{eq-H1*S1a}.

\subsection{Limit two-scale equations -- Fluid response in $\Om \times Y_f$}

We recall the separability of the fluid and solid responses of the
microstructure as discussed \eg in \cite{Rohan-zamp2020}, \cf
\cite{Clopeau-Mikelic-etal-MCM2001,Mik-Wheel-2012}. It enables to derive the
limit equations of the fluid flow in the microstructure formally independently
of the solid part. Obviously, the fluid and solid responses are coupled at the
macroscale. We omit details on the quasistatic flow homogenization -- the
result can be adopted from \cite{Rohan-zamp2020} dealing the dynamic
fluid-structure interaction in the porous media.

Upon substituting the asymptotic expansions \eq{eq-FS14a}-\eq{eq-FS14c} into \eq{eq-pzfl2}$_{3,4}$, 
 the following limit equations are obtained,
\begin{equation}\label{eq-FS15}
  \begin{split}
& \int_\Om \intY_{Y_f}\Dop\eeby{\hat\wb}:\eeby{\hat\thetabf} \dVy\dVx  -\int_\Om \intY_{Y_f}p^1\nabla_y\cdot\hat\thetabf \dVy\dVx \\
& \quad = \int_\Om(\fbav^f +  p^0 \nabla_x)\dVx \cdot\intY_{Y_f}\hat\thetabf\dVy\;,\\
& \int_\Om \left( \nabla_x\cdot\left(\phi_f\dot\ub^0 +  \intY_{Y_f}\hat\wb\dVy\right) + \intY_{Y_f}\nabla_y\cdot \dot{\wtilde{\ub}}^1\dVy
 + \gamma p^0\right) q^0 \dVx = 0\;,\\
& \int_\Om\intY_{Y_f}q^1\nabla_y\cdot\hat\wb \dVy\dVx = 0\;,
  \end{split}
\end{equation}
for all $\hat\thetabf \in L^2(\Om;\HpdbO(Y_f))$, $q^0 \in L^1(\Om)$
and $q^1 \in L^2(Y_2)$. Above, $\fbav^f$ is the averaged volume force satisfying
$\int_{Y_f}\fb^f\cdot\hat\thetabf = \fbav^f\cdot\int_{Y_f}\hat\thetabf$.

Now the autonomous characteristic responses $(\wb^k, \pi^k)$, $k=1,2,3$ can be defined,
such that (the summation in $k$ applies),
\begin{equation}\label{eq-A10a}
\begin{split}
  \hat\wb & = \wb^k(\what f_k - \pd_k^x p^0)\;,\quad
 p^1  = \pi^k(\what f_k - \pd_k^x p^0)\;.
\end{split}
\end{equation}
The characteristic responses  $\wb^{k}(y)$ are solutions of the following problem: Find $\wb^{k}\in\HpdbO(Y_f)$ and $\pi^{k}\in L^2(Y_f)$, such that
\begin{equation}\label{eq-A10c}
  \begin{split}
    \aYf{{\wb}^k}{\vb} -\ipYf{{\pi}^k}{\nabla_y\cdot\vb}  & = \ipYf{\onebm_k}{\vb}\;, \\ 
    \ipYf{q}{\nabla_y\cdot{\wb}^k} & = 0\;, 
\end{split}
\end{equation}
for all $\vb \in \HpdbO(Y_f)$ and $q \in L^2(Y_f)$.
Using the macroscopic seepage velocity defined by
\begin{equation}\label{eq-FS17a}
\phi_f\wb^0 := \intY_{Y_f} \hat \wb\dVy\;,
\end{equation}
mass conservation \eq{eq-FS15}$_2$ {yields} the macroscopic flow equation,
\begin{equation}\label{eq-FS17}
  \begin{split}
\gamma \dot p^0 +  \nabla_x\cdot[\phi_f(\wb^0 + \dot\ub^0)] +  \intY_{Y_f} \nabla_y\cdot\dot{\wtilde\ub}^1 \dVy& = 0\;,
  \end{split}
\end{equation}
where $\wtilde\ub^1$ is to be substituted using the scale-splitting formula \eq{eq-L2b}, defined below.

\subsection{Limit two-scale equations -- Solid response in $\Om \times Y_s$}

To obtain the limit of the solid phase equilibrium, the interaction term involving the traction of the stress in the fluid expressed in terms of the fluid velocity $\wb^\veps + \dot{\ub}^\veps$ and the pressure,
\begin{equation}\label{eq-FS18}
  \begin{split}
    \Ical^\veps(\sigmabf_f^\veps,\vb^\veps) & = \int_{\Gamma_\fsi^\veps} \nb^\sx\cdot\sigmabf_f^\veps\cdot\vb^\veps\dSy\;,\quad\mbox{ where }\quad
   \sigmabf_f^\veps  = -{p}^\veps \Ib + 2\mu^\veps \eeb{\wb^\veps + \dot{\wtilde\ub}^\veps}\;.
 \end{split}
\end{equation}
Recall that in \cite{Rohan-Lukes-PZ-porel}, only the disconnected inclusions were considered (no flow), so that the limit term depends naturally on $p^0$. For the quasistatic flow, according to \cite{Rohan-zamp2020}, the limit interaction term yields
\begin{equation}\label{eq-FS20}
  \begin{split}
   \Ical^\veps(\sigmabf_f^\veps,\vb^\veps) & \rightarrow
     \int_\Om \phi_f\left( p^0\nabla + \fbav^f\right)\cdot\vb^0\dVx\\ & - \int_\Om p^0 \intY_{\Gamma_\fsi} \vb^1\cdot \nb^\sx\dSy\dVx 
    +  \int_{\pd_p\Om} \bar\phi_f  \bar\bb^f \cdot\vb^0\dSx\;,
\end{split}
\end{equation}
where $\vb^0 \in \Ucalbf_0(\Om)$, $\vb^1(x,\cdot) \in \Hpdb(Ys)$, and
$\bar\phi_f$ is the boundary porosity. Note that on $\Gamma_p^\veps$, the
prescribed pressure $p_\pd$ represents the boundary traction stress acting on
the fluid part. Hence, by virtue of \eq{eq-FS20a} the condition prescribing the
averaged traction stress $\bar\phi_f\bar\bb^f = -\nb p_\pd$ on $\pd_p\Om$
emerges in \eq{eq-FS20} due to the convergence of $\bb^{f,\veps}$ acting on
$\Gamma_p^\veps$. Further issues concerning the boundary conditions are
discussed in Section~\ref{sec:macroproblem}.

The homogenization procedure explained in \cite{Rohan-Lukes-PZ-porel} is pursued in analogy to derive the resulting limit equations. 
Thus, we obtain the limit of \eq{eq-pzfl2} using \eq{eq-FS20}. We report here separately the local and the global subproblems; the local ones are obtained for vanishing $\vb^0$, such that
\begin{equation}\label{eq-L*S1}
\begin{split}
& \intY_{\Om \times Y_s}
\eeby{\vb^1} : \Aop \GrxyS{\ub^0}{\ub^1}\dVxy
- \intY_{\Om \times Y_m }\eeby{\vb^1} :\bar\parg^T\nabla_y\hat\vphi^0 \dVxy\\
& =  \int_\Om p^0 \intY_{\Gamma_f} \vb^1 \cdot \nb^\fx \dSy \dVx  \;,\\
& \intY_{\Om \times Y_m}\nabla_y\hat\psi^0 \cdot [\bar\parg:\GrxyS{\ub^0}{\ub^1} + \bar\db\nabla_y\hat\vphi^0]\dVxy = \intY_{\Om \times Y_m}\bar\rho_E \hat\psi^0\dVxy\;,
\end{split}
\end{equation}
must hold for all $\vb^1 \in \Lb^2(\Om;\Hpdb(Y_m))$ and $\hat\psi\in L^2(\Om;W_{\#*}(Y_m)$.
Due to the linearity, the standard scale decoupling applies by introducing the correctors functions (the characteristic responses) such that the two-scale functions can be expressed, as follows:
\begin{equation}\label{eq-L2b}
\begin{split}
\ub^1(x,y) & = \omegabf^{ij}e_{ij}^x(\ub^0)  - p^0 \omegabf^P + \omegabf^\rho \rho_E + \sum_\alpha\hat\omegabf^\alpha\bar\vphi^\alpha\;,\\
\vphi^0(x,y) & = \hat\eta^{ij}e_{ij}^x(\ub^0) - p^0 \hat\eta^P + \hat\eta^\rho \rho_E + \sum_\alpha\hat\vphi^\alpha\bar\vphi^\alpha\;.
\end{split}
\end{equation}
All $\omegabf,\hat\eta$ and $\hat\vphi$ are Y-periodic, representing the  displacements in the
entire solid part, $Y_{m*} = Y_m \cup Y_*$ and the electric potential in the
matrix part $Y_m$.

\subsubsection{Local characteristic responses}\label{sec-charpb}
The characteristic microproblems of the piezo-poroelastic medium are expressed using  the  bilinear forms:
\begin{equation}\label{eq-L3}
\begin{split}
\aYms{\ub}{\vb} & = \intY_{Y_{m*}} [\Dop \eeby{\ub}]:\eeby{\vb}\dY\;,\\
\gYm{\ub}{\psi} & = \intY_{Y_m} g_{kij} e_{ij}^y(\ub) \pd_k^y \psi\dY\;,\\
\dYm{\vphi}{\psi} & = \intY_{Y_m} [\db\nabla_y \vphi]\cdot\nabla_y \psi\dY\;.
\end{split}
\end{equation}
The local characteristic responses are solutions of the following autonomous problems (independent of the macroscopic responses):
\begin{itemize}
\item Find
$(\omegabf^{ij},\hat\eta^{ij})\in \Hpdb(Y_{m*})\times H_{\#0*}^1(Y_m)$ for any $i,j = 1,2,3$
satisfying
\begin{equation}\label{eq-mip1}
\begin{split}
\aYms{\omegabf^{ij} + \Pibf^{ij}}{\vb} - \gYm{\vb}{\hat\eta^{ij}}& = 0\;, \quad \forall \vb \in  \Hpdb(Y_{m*})\;,\\
\gYm{\omegabf^{ij} + \Pibf^{ij}}{\psi} + \dYm{\hat\eta^{ij}}{\psi}& = 0\;, \quad \forall \psi \in  H_{\#0*}^1(Y_m)\;,
\end{split}
\end{equation}
where  $\Pibf^{ij} = (\Pi_k^{ij})$, $i,j,k   = 1,2,3$ with components $\Pi_k^{ij} = y_j\delta_{ik}$ represents the homogeneous displacement field due to the principal macroscopic strain modes.
\item Find
$(\omegabf^P,\hat\eta^P)\in \Hpdb(Y_{m*})\times H_{\#0*}^1(Y_m)$
satisfying
\begin{equation}\label{eq-mip2}
\begin{split}
\aYms{\omegabf^P}{\vb} - \gYm{\vb}{\hat\eta^P}& = -\intY_{\Gamma_{f}} \vb\cdot \nb^\fx \dSy\;, \quad
\forall \vb \in  \Hpdb(Y_m) \;,\\
\gYm{\omegabf^P}{\psi} + \dYm{\hat\eta^P}{\psi}& = 0\;, \quad \forall \psi \in  H_{\#0*}^1(Y_m)\;.
\end{split}
\end{equation}
\item Find
$(\omegabf^\rho,\hat\eta^\rho)\in \Hpdb(Y_{m*})\times H_{\#0*}^1(Y_m)$
satisfying
\begin{equation}\label{eq-mip3}
\begin{split}
\aYms{\omegabf^\rho}{\vb} - \gYm{\vb}{\hat\eta^\rho}& = 0\;, \quad
\forall \vb \in  \Hpdb(Y_m) \;,\\
\gYm{\omegabf^\rho}{\psi} + \dYm{\hat\eta^\rho}{\psi}& = \intY_{\Gamma_{mc}} \psi \dSy\;, \quad \forall \psi \in  H_{\#0*}^1(Y_m)\;.
\end{split}
\end{equation}
\item Find
$(\hat\omegabf^\alpha,\hat\vphi^\alpha)\in \Hpdb(Y_{m*})\times H_{\#0,\alpha}^1(Y_m)$
satisfying, for $\alpha = 1,2,\dots,\alpha^*$,
\begin{equation}\label{eq-mip4}
\begin{split}
\aYms{\hat\omegabf^\alpha}{\vb} - \gYm{\vb}{\hat\vphi^\alpha}& = 0 \;, \quad
\forall \vb \in  \Hpdb(Y_m) \;,\\
\gYm{\hat\omegabf^\alpha}{\psi} + \dYm{\hat\vphi^\alpha}{\psi}& = 0\;, \quad \forall \psi \in  H_{\#0*}^1(Y_m)\;.
\end{split}
\end{equation}
\end{itemize}

In \eq{eq-mip2}, the interface integral can be rewritten using the  volume integral: for any $\vb \in \Hpdb(Y_{m*})$, 
\begin{equation}\label{eq-mip0a}
\begin{split}
\intY_{\Gamma_f} \vb\cdot\nb^\fx\dSy = - \intY_{\Gamma_f} \vb\cdot\nb^\mx \dSy= - \intY_{Y_{m*}}\nabla_y\cdot\vb\dVy\;.
\end{split}
\end{equation}

\subsection{Macroscopic problem -- linear model of quasistatic flow and deformation}\label{sec:macroproblem}

At the macroscopic level, the pore flow is coupled with the solid deformation in terms of the macroscopic quasistatic  equilibrium arising from the solid part, and the mass conservation arising from the fluid part. The macroscopic state of the piezo-poroelastic medium is given by the couple $(\ub^0,p^0)$ whereas the macroscopic electric potentials $\{\bar\vphi^\alpha\}$ are involved as the space-time control fields. We introduce the admissibility sets, 

\begin{equation}\label{eq-UPspaces}
\begin{split}
\Ucalbf(\Om) = \{\ub \in \Hdb(\Om)|\; \ub = \ub_\pd \mbox{ on }\pd_u\Om\}\;,\\
\Qcal(\Om) = \{ p \in H^1(\Om)|\; p = p_\pd \mbox{ on }\pd_p\Om\}\;.
\end{split}
\end{equation}
Further we employ the test function spaces $\Ucalbf_0(\Om)$ and $\Qcal_0(\Om)$, associated with $\Ucalbf(\Om)$ and $\Qcal(\Om)$, respectively, by virtue of prescribing the zero Dirichlet boundary conditions.

The limit boundary conditions of the Neumann type must be specified. For the solid phase, the traction stress is acting on $\Gamma_\sigma^\veps$ represented by $\pd_\sigma \Om$ in the limit. The prescribed fluid pressure on $\Gamma_p^\veps$ presents the Dirichlet-type  boundary conditions, however, in the limit model, it acts also as the traction stress $-\nb p_\pd$ in the macroscopic equilibrium equation for the mixture. 

\subsubsection{Homogenized coefficients of the piezo-poroelastic medium}

We first recall the limit equations \eq{eq-FS17} and \eq{eq-FS17a} arising from the flow subproblem. Substituting \eq{eq-A10a} in \eq{eq-FS17} yields the hydraulic permeability $\Kb = (K_{ij})$,
\begin{equation}\label{eq-K}
K_{ij} = \intY_{Y_f} w_i^j\dVy\;,\quad K_{ij} = K_{ji}\;.
\end{equation}
Its symmetry and positive semidefiniteness follows by \eq{eq-A10c} upon substituting there $\vb:=\wb^j$. Recall that the rank deficiency of $K_{ij}$ appears if the porosity is not a simply connected domain.

%
The quasistatic equilibrium equation is obtained in the standard way due to the limit interaction integral \eq{eq-FS20} yielding the traction vector $\bar\phi_f\bar\bb^f$ on $\pd_p \Om$. In general, the effective volume and traction forces are involved,
\begin{equation}\label{eq-H1*S1a}
\begin{split}
  \fbav = \phi_s \fbav^s + \phi_f \fbav^f\;,\quad \bbav = \bar\phi_s \bar\bb^s + \bar\phi_f \bar\bb^f\;,
  \quad
  \phi_d\fbav^d = \intY_{Y_d} \fb^d\dVy\;,\quad d = s,f\;,
\end{split}
\end{equation}
where $\fbav^d$ are the average volume forces.
The effective traction stress $\bbav$ is defined in terms of surface fractions $\bar\phi_d$, $d=f,s$ of the two phases and using the mean traction stresses $\bar\bb^f$ and $\bar\bb^s$ loading the fluid and solid, respectively.


Macroscopic functions $\ub^0 \in \Ucalbf(\Om)$ and $p^0\in L^2(\Om)$ satisfy 
\begin{equation}\label{eq-H1*S1}
\begin{split}
\int_{\Om} &
e_{ij}^x(\vb^0)\left[\aYms{\ub^1-\Pibf^{kl}e_{kl}^x(\ub^0)}{\Pibf^{ij}} - \gYm{\Pibf^{ij}}{\hat\vphi^0 }  \right] \dVx \\
& -\int_{\Om} p^0 \phi_f \nabla_x\cdot \vb^0 \, \dVx  =
\int_{{\Om}}\fbav \cdot \vb^0 \, \dVx + \int_{\pd {\Om}}\bbav \cdot \vb^0\dSx\;,
\end{split}
\end{equation}
and the fluid mass conservation, see \eq{eq-FS17}, 
\begin{equation}\label{eq-H1*S2}
\begin{split}
\int_{\Om}&  q^0\left(\phi_f\nabla_x\cdot \ub^0
- \intY_{\Gamma_f} \ub^{1} \cdot \nb^\mxe \dSy \right ) \, \dVx + \int_{\Om}\phi_f\nabla_x\cdot\wb^0 \, \dVx
+\gamma \int_{\Om}\phi_f  p^0 q^0\, \dVx = 0\;,
\end{split}
\end{equation}
for all $(\vb^0,q^0) \in \Ucalbf_0(\Om)\times \Qcal_0(\Om)$. It should be noted that, the electric charge conservation is satisfied at the microlevel.
In \eq{eq-H1*S1} and  \eq{eq-H1*S2}, the two-scale functions are substituted using the split \eq{eq-L2b}. Upon collecting the terms matching different macroscopic function, the expressions for the homogenized coefficients are identified, see \cite{Rohan-Lukes-PZ-porel},
\begin{equation}\label{eq-HC1}
\begin{split}
A_{klij} & = \aYms{\omegabf^{ij} + \Pibf^{ij}}{\omegabf^{kl}+\Pibf^{kl}} + \dYm{\hat\eta^{kl}}{\hat\eta^{ij}}\;, \\
B_{ij} & = \aYms{\omegabf^P}{\Pibf^{ij}} - \gYm{\Pibf^{ij}}{\hat\eta^P} + \phi_f\delta_{ij} =
-\intY_{Y_m} \nabla_y\cdot \omegabf^{ij}\dVy + \phi_f\delta_{ij}
\;,\\
M & 
= \aYms{\omegabf^P}{\omegabf^P} + \dYm{\hat\eta^P}{\hat\eta^P}
+ \phi_f\gamma\;,\\
H_{ij}^\alpha & = \aYms{\hat\omegabf^\alpha}{\Pibf^{ij}} -  \gYm{\Pibf^{ij}}{\hat\vphi^\alpha}\;,\\
S_{ij} & = \aYms{\omegabf^\rho}{\Pibf^{ij}} -  \gYm{\Pibf^{ij}}{\hat\eta^\rho}\;,\\
R & = - \intY_{\Gamma_{f}} \omegabf^\rho\cdot \nb^\fx\dSy\;,\\
Z^\alpha & = - \intY_{\Gamma_f}\hat\omegabf^\alpha\cdot \nb^\fx \dSy\;.
\end{split}
\end{equation}
Above the symmetric expressions were derived using the equations of the characteristic problems \eq{eq-mip1}-\eq{eq-mip4}.
Clearly, $A_{klij} = A_{ijkl} = A_{klji}$, $B_{ij}=B_{ji}$, as expected for the Biot-type continuum and also $S_{ij}=S_{ji}$, $H_{ij}^\alpha =H_{ji}^\alpha$.
We shall employ the ``boldface'' notation $\Aop = (A_{ijkl})$, $\Bb = (B_{ij})$, $\Sb = (S_{ij})$ and $\Kb=(K_{ij})$, whereas
$\ul{\Hb} = (H_{ij}^\alpha)$ and $\ul{Z} = (Z^\alpha)$ matches with $\ul{\vphi} = (\bar\vphi^\alpha)$, so that
$\sum_\alpha \Hb^\alpha \bar\vphi^\alpha = \ul{\Hb}\cdot\ul{\vphi}$.

\subsubsection{Macroscopic problem -- linear model of coupled flow and deformation}

The limit macroscopic equations are obtained from the two-scale limit equations
with the only non-vanishing macroscopic test functions $\vb^0$ and $q^0$.
 The weak formulation reads: 
For any time $t > 0$, find a couple $(\ub^0(t,\cdot),p^0(t,\cdot)) \in \Ub(\Om)\times P(\Om)$ which satisfies
\begin{equation}\label{eq-S25}
\begin{split}
&\int_{\Om} \left (\Aop  \eeb{\ub^0}
-  p^0  \Bb + \ul{\Hb}\cdot \ul{\vphi^0} \right ):\eeb{\vb^0}\dVx\\
& \quad = \int_{\Om} \fbav \cdot \vb^0 \dVx
+ \int_{\pd \Om}\bbav \cdot \vb^0\dSx\;, \quad \forall \, \vb^0 \in
\Ucalbf_0(\Om) \;,\\
&\int_{\Om} q^0\left (\Bb :\eeb{\dot\ub^0} +  \dot p^0 M -\ul{Z} \cdot \dot{\ul{\vphi^0}}\right)\dVx + \int_{\Om} \frac{\Kb}{\bar\mu}\left(\nabla_x p -\what{\fb}^f\right) \cdot \nabla_x q^0 \dVx\\
& \quad = 0\;,\quad \forall q^0 \in \Qcal_0(\Om).
\end{split}
\end{equation}


Upon integrating by parts in \eq{eq-S25}, the following system of equations
governing the fluid flow in deforming piezo-poroelastic medium in domain $\Om
\subset \RR^3$ is obtained,
\begin{equation}\label{eq-M1}
\begin{split}
-\nabla\cdot \sigmabf^H(\ub^0,p^0) = \hat\fb \;,\quad \mbox{ in }  \Om\;,\\
\mbox{ where } \quad \sigmabf^H(\ub^0,p^0) = {\Aop}\eeb{\ub^0} - p^0  {\Bb}
+ \sum_\alpha  {\Hb^\alpha} \bar\vphi^\alpha + {\Sb \rho_E} \;,
\end{split}
\end{equation}
and
\begin{equation}\label{eq-M2}
\begin{split}
    {\Bb}:\eeb{\dot\ub^0} +  {M} \dot p^0 + \nabla\cdot\wb^0 = \sum_\alpha  {Z^\alpha} \dot{\bar\vphi}^\alpha + {R  \dot\rho_E}\;,\\
 \mbox{ where } \quad \wb^0 = -\frac{1}{\bar\mu}{\Kb}(\nabla p^0 - \fb^f)\;.
\end{split}
\end{equation}
\revE{For vanishing piezoelectric coefficients $\parg \equiv \zerobf$, the Biot model of the poroelasticity is recovered, since, in \eq{eq-mip1}-\eq{eq-mip3}, all $\eta^* \equiv 0$, such that expressions \eq{eq-HC1} yield the standard effective poroelasticity coefficients $\Aop,\Bb$ and $M$, whereas 
${\Hb^\alpha},R$ and ${Z^\alpha}$ vanish.}

In the numerical examples reported in Section~\ref{sec:numex},  $\rho_E \equiv 0$, as well as the volume forces are not considered.

We shall employ boundary conditions 
  prescribed on $\pd \Om$ being decomposed into disjoint parts $\Gamma_i$, $i = 0,1,2$, such that
$\pd \Om = \Gamma_0 \cup \Gamma_1 \cup \Gamma_2$.
 Without loss of generality, we consider $p_\pd = \bar P^i$ on $\Gamma_i$, $i = 1,2$, and $\ub_\pd = \zerobf$ on $\Gamma_u \equiv \Gamma_1$  
\begin{equation}\label{eq-H*S2}
\begin{split}
  p^0 = \bar P_i \mbox{ on } \Gamma_i \;,i = 1,2\;,\\
  \nb\cdot\wb^0 = 0 \mbox{ on } \Gamma_0 \;,\\
  \ub^0 = \zerobf\mbox{ on } \Gamma_1\;,\\
  \nb\cdot\sigmabf^H= \zerobf\mbox{ on } \Gamma_0\;,\\
   \nb\cdot\sigmabf^H= -\bar P_2\nb \mbox{ on } \Gamma_2\;.
\end{split}
\end{equation}


\paragraph{Problem formulation}
Let us assume zero initial conditions, vanishing $(\ub^0,p^0)$ at $t = 0$ in $\Om$
corresponding to unloaded initial configuration. For given $P^k(t)$ functions
of $t>0$, whereby $P^k(0) = 0$ for $k=1,2$ and given electric actuation by
voltages $\{\bar\vphi^\alpha(t,x)\}_\alpha$ defined for $(t,x) \in [0,T]\times\Om$, with
$T>0$, whereby $\bar\vphi^\alpha(0,\cdot) = 0$ in $\Om$, find a solution
$(\ub^0,p^0)(t,x)$ satisfying \eq{eq-M1}, \eq{eq-M2} and \eq{eq-H*S2}.

\section{Modified formulation for nonlinear response}\label{sec-varHC}

Although the homogenized model \eq{eq-M1}-\eq{eq-M2} was derived for a fixed
reference configuration, as explained in \ref{app-incr}, it can be modified \revE{to capture the nonlinear effects associated with flow in the deforming microstructure. This is important to simulate the peristalsis-driven flow, see \ref{sec-1D} where the 1D reduced model is employed to show the importance of respecting the permeability dependence on the deformation.} We adhere the
approach suggested in \cite{Rohan-Lukes-nlBiot2015} which enables to
approximate the ``deformation-dependent'' homogenized coefficients (DDHC). \revE{The notion of DDHC actually reflects the dependence of the microstructure configuration on the local deformation which, in the homogenization limit, besides the macroscopic strain depends also on all macroscopic variables involved in the expansion of displacement $\ub^1$ and, thereby, on $\ub^\mic$. In \ref{app-incr}, we show the link between the sensitivity analysis employed below to get such linear expansions for the DDHC and linearized transformations between the initial and deformed configurations.}

To lighten the notation in this section, we drop the superscript $(\cdot)^0$,
which was previously used to denote macroscopic quantities.

\subsection{Approximation of deformation-dependent homogenized coefficients}\label{sec:appHC}

The approximation is based on the 1st order Taylor expansion which uses the
sensitivity of homogenized coefficients \wrt the macroscopic variables. We
shall consider all tensors and vectors labelled by $\tilde{}$ to be defined in
the perturbed microconfiguration $\tilde\Mcal({\ub^\mic}(x,\cdot),Y)$ which is
characterized by the deformed ``local'' representative periodic cell $Z(x) \equiv \tilde
Y(x) = Y + \{{\ub^\mic}(x,Y)\}$ for $x \in \Om$. By ${\ub^\mic}$ we refer to the
displacement field reconstructed in the microstructures,
\begin{equation}\label{eq-S22}
\begin{split}
  {\ub^\mic}(x,y) = \Pibf^{ij}(y) e_{ij}^x(\ub(x)) + \ub^1(x,y)\;,
\end{split}
\end{equation}
where $\ub^1(x,y)$ is defined in \eq{eq-L2b}. 
 Since
$\tilde Y(x)$ depends on the macroscopic coordinate, the
microstructure is perturbed from its periodic structure. However, the
homogenization procedure can still be applied and the homogenized
coefficients denoted by $\Hop(\tilde\Mcal({\ub^\mic},Y))$, in a generic sense, can be computed for the perturbed geometry which is determined by the field ${\ub^\mic}$.
Due to the
sensitivity analysis, as explained in \cite{Rohan-Lukes-nlBiot2015} and presented below, and by virtue of the split formulae
\eq{eq-L2b}, the perturbed coefficients
$\Hop(\tilde\Mcal({\ub^\mic},Y)) \approx \tilde \Hop(\eeb{\ub},p,\{\bar\vphi^\alpha\})$ can be approximated using the first
order expansion formulae which have the generic form
\begin{equation}\label{eq-S24}
\begin{split}
\tilde \Hop(\eeb{\ub},p,\{\bar\vphi^\alpha\}_\alpha) & =  \Hop^0 + \delta_{\eb} \Hop^0 : \eeb{\ub} + \delta_p \Hop^0 p + \sum_\alpha\delta_{\vphi,\alpha} \Hop^0 \bar\vphi^\alpha\;,\\
(\delta_{\eb}H^0)_{ij} & := \big(\pd_{\eb} (\delta \Hop^0\circ{\ub^\mic})\big)_{ij} = \delta \Hop^0 \circ {(\omegabf^{ij} + \Pibf^{ij})}\;,\\
\delta_p \Hop^0 & := \pd_p (\delta \Hop^0\circ{\ub^\mic}) = \delta \Hop^0 \circ (-\omegabf^P)\;,\\
\delta_{\vphi,\alpha} \Hop^0  & := \pd _{\vphi,\alpha}(\delta \Hop^0\circ{\ub^\mic}) = \delta \Hop^0 \circ \hat\omegabf^\alpha\;.
\end{split}
\end{equation}
Coefficients $\Hop^0$ are computed using \eq{eq-HC1} and \eq{eq-K} for the
reference ``initial'' configuration $\Mcal(Y)$ and $\delta \Hop^0$ are
sensitivities which are derived in \ref{sec-sa} using the ``design velocity field'' method according to the shape sensitivity terminology. Thus, for any $\Hop$ we can obtain a sensitivity expression $\dlt\Hop(\vec\Vcal) = \dlt\Hop\circ\vec\Vcal$ which is linear in the vector field $\vec\Vcal$, such that for perturbed microconfiguration $z_i(y,\tau) = y_i + \tau\Vcal_i(y)$, $y \in Y$, $i=1,2,3$, $\tilde\Hop = \Hop^0 + \tau \dlt\Hop\circ\vec\Vcal$. Clearly, the derivatives \eq{eq-S24} are obtained upon substituting $\vec\Vcal$ by the specific vector field.

It is worth of noting that the potentials, 
$\{\bar\vphi^\alpha(x,t)\}_\alpha$ are considered as the (given) control functions which, however, they influence directly the homogenized model parameters due to the expansion \eq{eq-L2b}.

\subsection{Semilinear formulation of the macroscopic problem}\label{sec:nonlin}

The variational formulation \eq{eq-S25} can be modified to respect the
deformation-dependent homogenized coefficients listed in \eq{eq-HC1} and
\eq{eq-K} which are further referred to by $\Hop^0$. For this, all $\Hop^0$ are
substituted using $\tilde \Hop$, in the generic sense, as defined in
\eq{eq-S24}. In the same context, by $\ol{\Hop}$ the coefficients of the
linearized subproblems are denoted.


It should be stressed out, that the microproblems \eq{eq-mip1}-\eq{eq-mip4} and \eq{eq-A10c} are solved  only once independently of the macroscopic problem. As well, the associated sensitivity analysis is performed for the reference non-perturbed configuration.

\subsubsection{\newE{Incremental formulation}}
We consider an incremental formulation based on the following modelling steps:

\paragraph{Step 1: Solid equilibrium in the deformed configuration.} Based on the smallness of the deformation assumption, we consider a perturbation of the spatial reference configuration $\hOmSe$, assuming the quasistatic equilibrium of the stresses $\tilde\sigmabf_s^\veps$ was attained for loads by volume forces $\tilde\fb^\veps$ and surface tractions $\tilde\bb^\veps$, including the fluid stress traction $\nb^\sx\tilde\sigmabf_f^\veps$ acting on $\hGmSFe$. The solid state displacement and electric potential increments, $\dlt\ub^\veps$ and $\dlt\vphi^\veps$, respectively, should satisfy
\begin{equation}\label{eq-E1}
  \begin{split}
    \int_\hOmSe [\tilde\sigmabf_s^\veps + \dlt\sigmabf_s^\veps(\dlt\ub^\veps,\dlt\vphi^\veps)] : \eeb{\vb} \dVx
     - \int_\hGmSFe\vb\cdot(\tilde\sigmabf_f^\veps+ \dlt \sigmabf_f^\veps)\cdot \nb^\sx \dSx \\
    =\int_\hOmSe(\tilde\fb^\veps + \dlt\fb^\veps)\cdot\vb \dVx
     + \int_{\pd_\sigma\hOmSe} (\tilde\bb^\veps + \dlt\bb^\veps)\cdot\vb\dSx 
\end{split}
\end{equation}
for all $\vb \in \Ucalbf_0^\veps(\hOmSe)$. 
Should the equilibrium be attained ``precisely'' in $\hOmSe$  (in the context of the finite increment associated with the ``$\dlt$'' perturbation used in the past increments), the actual perturbation pair  $(\dlt\ub^\veps,\dlt\vphi^\veps)$ must satisfy
\begin{equation}\label{eq-E2}
  \begin{split}
    \int_\hOmSe \dlt\sigmabf_s^\veps(\dlt\ub^\veps,\dlt\vphi^\veps) : \eeb{\vb}\dVx
     - \int_\hGmSFe\vb\cdot\dlt \sigmabf_f^\veps\cdot \nb^\sx \dSx \\
    =\int_\hOmSe \dlt\fb^\veps\cdot\vb\dVx + \int_{\pd_\sigma\hOmSe} \dlt\bb^\veps\cdot\vb\dSx \;,\\
  \int_\hOmMe \vec D^\veps(\dlt\ub^\veps,\dlt\vphi^\veps)\cdot\nabla\psi \dVx = 0\;,
\end{split}
\end{equation}
for all $\vb \in \Ucalbf_0^\veps(\hOmSe)$ and $\psi \in \Vcal_*(\hat\Om_m^\veps,\{\Gamma_*^{\alpha,\veps}\}_\alpha)$. In \eq{eq-E1}, we do not present the electric charge conservation for the sake of brevity.

\paragraph{Step 2: Homogenization of the linear subproblems.}
The fluid-structure interaction problem is posed in the reference spatial configuration; the solid response in $\hOmSe$ is governed by \eq{eq-E2}, which formally is represented by \eq{eq-pzfl2}, the 1st and the 2nd equation, where ${\Om_{m*}^\veps}$ is considered as  $\hOmSe$ and ${\ub^\veps,\vphi^\veps}$ is replaced by  $(\dlt\ub^\veps,\dlt\vphi^\veps)$.
In the fluid channel $\hOmFe$, the Stokes flow model is considered, as described in \eq{eq-pzfl2}, the 3rd and the 4th equation, where $\dot{\tilde\ub}^\veps$ is replaced by  $\dlt\dot{\tilde\ub}^\veps$ by virtue of the velocity definition. The homogenization procedure applies, as described above, leading to the macroscopic problem formally identical with \eq{eq-S25}, to be satisfied by the increments $(\dlt\ub^0,\dlt p^0)$ and involving the control increments $\dlt\ul{\vphi}^0$. Since it is obtained by homogenizing the spatial heterogeneous configuration, in \eq{eq-S25}, the domain $\Om$ is the actual reference spatial macroscopic configuration and the effective poroelastic coefficients and the permeability describe the effective constitutive corresponding to the locally deformed microconfigurations.

\paragraph{Step 3: Modified poroelastic coefficients, incremental constitutive laws.} Besides the permeability $\Kb$, the other effective material parameters describe the tangent moduli of the macroscopic stress $\ol{\sigmabf}_S$. These all effective parameters can be computed in the local deformed representative cells $\tilde Y(x)$, however, we employ their approximation based on the sensitivity analysis (SA), as described above. In \ref{app-incr}, a brief derivation of the incremental macroscopic constitutive law is presented to justify relevance of the  SA based approximation. The stress is given by a hypoelastic constitutive law which can be presented, denoting by $\ssb(t,x) = (\ub,p,\ul{\vphi})(t,x)$ the state of the (homogenized) poroelastic medium, for $x \in \Om$,
\begin{equation}\label{eq-E2b}
  \begin{split}
    \sigmabf(\ssb(t,x)) & = \int_0^t \dlt_{\ssb} \sigmabf({\ssb}(\tau,x))\circ\hat\ssb(x)\dd \tau\;,\quad \ssb(t,x) = \hat\ssb(x) t\;,\\
    \dlt_{\ssb} \sigmabf\circ\dlt\ssb(\tau,x) & = \tilde\Aop  \eeb{\dlt\ub}
- \dlt p  \tilde\Bb + \tilde{\ul{\Hb}}\cdot \ul{\dlt\vphi}\;,
\end{split}
\end{equation}
where the homogeneity $\dlt\ssb(\tau,x) = \hat\ssb(x) \dd \tau$ applies due to elasticity/piezoelectricity of the solid.
By virtue of the incremental algorithm introduced below, in Section~\ref{sec:nonlin-macro}, the simple updating  $\sigmabf := \bar\sigmabf + \dlt\sigmabf$ can be employed in each iteration.

\paragraph{Semi-discretized macroscopic problem}

Within the time discretization with levels $t_k = k \Delta t$, $k=0,1,2,\dots$
we consider two consecutive steps and the abbreviated notation $t:=t_k$, $t\prevstep :=
t_{k-1}$. By $(\ub,p)\approx (\ub(t,\cdot),p(t,\cdot))$ we shall refer to the
unknown fields at time $t = t_k \in\{t_i\}_i$, whereas $(\ub\prevstep,p\prevstep) \approx
(\ub(t\prevstep,\cdot),p(t\prevstep,\cdot))$, thus being associated with time $t - \Delta t$.
Using the implicit approximation $\dot\ub \approx (\ub - \ub\prevstep)/\Delta t$,
$\dot p \approx (p - p\prevstep)/\Delta t$ and $\dot \vphi \approx (\vphi -
\vphi\prevstep)/\Delta t$,

Since the nonlinear problem is solved using a Newton-type iterative method, it
is convenient to introduce the residual-based incremental formulation.
Let us introduce the residual function, 
\begin{equation}\label{eq-rif1}
  \begin{split}
  & \Psi^t((\ub,p),\Fop;(\vb,q))  =
\int_{\Om} \left (\tilde\Aop  \eeb{\ub}
-  p  \tilde\Bb + \tilde{\ul{\Hb}}\cdot \ul{\vphi} \right ):\eeb{\vb}\dVx \\
& +\int_{\Om} q\left (\tilde\Bb :\eeb{\ub-\ub\prevstep} +  (p - p\prevstep) \tilde M - \tilde{\ul{Z}}\cdot(\ul{\vphi} - \ul{\vphi}\prevstep)\right)\dVx \\
& + \int_{\Om} \frac{\Delta t}{\bar\mu}\tilde\Kb\left(\nabla_x p -\fbav^f\right) \cdot \nabla_x q \dVx - \left( \int_{\Om} \fbav \cdot \vb \dVx
+ \int_{\pd_\sigma \Om}\bbav \cdot \vb\dSx \right)\;,
\end{split}
\end{equation}
evaluated at time $t$ (hence the superscript in $\Psi^t$), where the ``tilde'' notation means the application of $\tilde \Hop$ in the generic sense of \eq{eq-S24} for all the homogenized coefficients, and $\Fop$ represents all  loads, or control given in  time and space domain $[0,T]\times\Om$. In particular, $\Fop$ comprises $\{\bar\vphi^\alpha(x,t)\}_\alpha$ and all the volume forces and surface tractions.

Now the time-discretized problem \eq{eq-S25} is reformulated, as follows:

For a given solution $(\ub\prevstep,p\prevstep)$, find a couple $(\ub,p) \in \Ucalbf(\Om)\times \Qcal(\Om)$ which satisfies
 \begin{equation}\label{eq-S25TD}
 \Psi^t((\ub,p),\Fop;(\vb,q)) = 0\quad\mbox{  for any }\quad (\vb,q)  \in \Ucalbf_0(\Om)\times \Qcal_0(\Om)\;.
 \end{equation}
It can be solved by successive iterations based on the
straightforward decomposition introduced in terms of
the recent approximation $(\bar \ub,\bar p)$ of the solutions $(\ub,p)\in \Ucalbf(\Om)\times \Qcal(\Om)$ and the correction $(\delta \ub, \delta p)\in  \Ucalbf_0(\Om)\times \Qcal_0(\Om))$,
\begin{equation}\label{eq-S50}
\begin{split}
\ub = \bar \ub + \delta\ub\;,\quad p = \bar p + \delta p\;.
\end{split}
\end{equation}
The first order expansion in $(\delta \ub, \delta p)$ leads to the obvious approximation 
\begin{equation}\label{eq-rif2}
  \begin{split}
    0  & = \Psi^t((\ub,p),\Fop;(\vb,q)) \approx \Psi^t((\bar\ub,\bar p),\bar\Fop;(\vb,q)) \\
    & + \delta_{(\ub,p)}  \Psi^t((\bar\ub,\bar p),\bar\Fop;(\vb,q))\circ(\dlt\ub,\dlt p)
    + \delta_{\Fop} \Psi^t((\bar\ub,\bar p),\bar\Fop;(\vb,q))\circ\dlt\Fop\;,
\end{split}
\end{equation}
to hold for any $(\vb,q) \in  \Ucalbf(\Om)\times \Qcal(\Om)$,
where $(\bar\ub,\bar p)$ is the recent iteration.

 The total variation $\delta\Psi^t$ can be expressed in terms of tangential incremental coefficients $\ol{\Hop}$ defined below, in \eq{eq-S40}-\eq{eq-S41}, see Section~\ref{sec-tan},
\begin{equation}\label{eq-sts6}
\begin{split}
&\delta_{(\ub,p)}\Psi^t((\bar \ub,\bar p),\bar\Fop;( \vb,q)) \circ (\delta\ub,\delta p)
 = \int_\Om \left(\ol{\Aop}\eeb{\delta\ub} - \delta p \ol{\Bb}\right):\eeb{\vb}\dVx\\
& \quad
  - \int_\Om \big(\pd_\eb{\fbav}\circ\eb(\delta\ub)+\pd_p{\fbav}\circ\delta p\big)\cdot \vb\dVx
-\int_{\pd\Om}\big(\pd_\eb{\bbav}\circ\eb(\delta\ub)+\pd_p{\bbav}\circ\delta p\big)\cdot \vb\dVx \\
& \quad
 +  \frac{\Delta t}{\bar\mu}\int_\Om \nabla q\cdot
\left(\ol{\Kb}\nabla \delta p + \ol{\Gb}:\eeb{\delta\ub} + \ol{\Qb}\delta p
\right)\dVx
 + \int_\Om q\big(\ol{\Db}:\eeb{\delta\ub}+\ol{M}\delta p\big)\dVx
\;.
\end{split}
\end{equation}
Variation $\delta_{\Fop} \Psi^t((\bar\ub,\bar p),\bar\Fop;(\vb,q))\circ\dlt\Fop$ can be expressed in analogy, however, since $\Fop$ is considered to be a given ``control'' at each time level and fixed during the iterations, \ie $\delta\Fop \equiv 0$,   $\delta_{\Fop} \Psi^t$ is not needed.

\subsubsection{Algorithm of time increment step}\label{sec:nonlin-macro}

Knowing the response at time $t\prevstep = t-\Delta t$, the one at time level $t$  is computed by solving \eq{eq-S25TD} iteratively using a modified Newton method based on the linearization due to \eq{eq-rif2}
explained above.
Since the boundary data, namely $p_\pd(t,\cdot)$, see \eq{eq-H*S2}), can depend
on time, we use the notation $\Ucalbf^t(\Om)$ and $\Qcal^t(\Om)$, for the
admissibility sets. For simplicity, however, we consider the split of $\pd\Om$
independent of time, thus, we keep using $\Ucalbf_0(\Om)$ and $\Qcal_0(\Om)$.

Given $(\ub\prevstep,p\prevstep) \approx \big(\ub(t\prevstep,\cdot),p(t\prevstep,\cdot)\big)$, compute $(\ub,p) \approx \big(\ub(t,\cdot),p(t,\cdot)\big)$:
\begin{list}{$\bullet$}{}
\item step 0:
  set the time-dependent control $\Fop^t:=\Fop(t)$ comprising the electric
  potential $\ul{\vphi}(t,\cdot$ given in $\Om$ and the boundary conditions,
  thereby $U^t(\Om)$ and $P^t(\Om)$ and define an approximation $(\bar \ub,\bar
  p) \approx (\ub\prevstep,p\prevstep)\in \Ucalbf^t(\Om)\times\Qcal^t(\Om)$,
  where the perturbed homogenized coefficients $\tilde \Hop$ (in the sense of the
  generic notation) are evaluated using $(\ub\prevstep,p\prevstep)$ inserted in
  \eq{eq-S24}$_1$ where $\Hop^0$ and the gradients $\dlt_{()}\Hop^0$ are
  given for the unperturbed microconfiguration (by off-line computing).

\item step 1: use $(\bar \ub,\bar p)$ to define the perturbed
  homogenized coefficients $\tilde \Hop$ and \eq{eq-rif1} to establish the residual
  $r(\vb,q) := \Psi^t((\bar \ub,\bar p),\Fop^t,(\vb,q))$, where the test
  functions $(\vb,q) \in \Ucalbf_0(\Om)\times\Qcal_0(\Om)$.

  If the desired tolerance $\epsilon > 0$ is achieved, i.e. $|r(\vb,q)| < \epsilon \|(\vb,q)\|$ for any $(\vb,q) \in \Ucalbf_0(\Om)\times\Qcal_0(\Om)$, put
$\ub\prevstep := \ub = \bar \ub$, $p\prevstep := p  = \bar p$,
and proceed to solving the response at the next time level $t\prevstep := t$. GOTO step 0.

\item step 2: solve the following linear problem for
  $(\delta\ub,\delta p) \in \Ucalbf_0(\Om)\times\Qcal_0(\Om)$:
\begin{equation}\label{eq-rif1a}
 \delta_{(\ub,p)}  \Psi^t((\bar\ub,\bar p),\Fop^t;(\vb,q))\circ(\dlt\ub,\dlt p)  = - \Psi^t((\bar\ub,\bar p),\Fop^t;(\vb,q)), 
\end{equation}
for any $(\vb,q) \in \Ucalbf_0(\Om)\times\Qcal_0(\Om)$,  
then update $\bar\ub := \ub = \bar \ub + \delta \ub$, $\bar p := p = \bar p + \delta p$,
and GOTO step 1.
\end{list}
\paragraph{Comments} In the space-discretized version, the standard residual
vector in the sense of the Galerkin finite element method replaces the
functional $r(\cdot,\cdot)$.
 
\subsection{Tangential incremental coefficients}\label{sec-tan}

The perturbations $\dlt\ub,\dlt p,\dlt \vphi$ employed in step 2 of the
Algorithm~\ref{sec:nonlin-macro}, see \eq{eq-rif1a}, substituted in the
integrals involved in $\Psi^t$, see \eq{eq-rif1}, yield approximate linear
expressions, as shown for the term related to the permeability,
\begin{equation}\label{eq-S35}
\begin{split}
& \int_\Om \nabla q \cdot \tilde\Kb(\eeb{\ub},p) \nabla p \dVx\\ 
& \quad = \int_\Om \nabla q \cdot (\Kb^0 + \pd_\eb \Kb^0\circ\eeb{\ub} + \pd_p \Kb^0\circ p + \pd_\vphi\Kb^0\circ \vphi) \nabla p\dVx\\
& \quad \approx 
\int_\Om \nabla q \cdot \left(\Kb^0 + \pd_\eb \Kb^0\circ(\eeb{\bar\ub } + \eeb{\delta \ub}) \right)\nabla {\bar{p}}\dVx\\
& \qquad + \int_\Om \nabla q \cdot \left(\pd_p \Kb^0\circ (\bar p + \delta p) + \pd_\vphi\Kb^0\circ (\bar\vphi + \delta \vphi) \right)\nabla {\bar{p}}\dVx\\
& \qquad + \int_\Om \nabla q\cdot \left(\Kb^0 + \pd_\eb \Kb^0\circ\eeb{\bar\ub } + \pd_p \Kb^0\circ \bar p + \pd_\vphi\Kb^0\circ \bar\vphi\right) \nabla {\delta{p}}\dVx\;,
\end{split}
\end{equation}
neglecting the ``$o(\dlt^2)$'' terms. In analogy, we proceed with all other integrals in  \eq{eq-rif1}.


In what follows we shall need some further notation employed in the linearization scheme: let $\Wb = \pd_\eb(\Xb:\bar\eb)\circ \circpare{~}$, then $\Wb\delta\eb =  \pd_\eb(\Xb:\bar\eb)\circ \delta\eb$; in analogy, let 
$X = \pd_p(\Xb\bar\eb)\circ \circparp{~}$, then $X\delta p =
\pd_p(\Xb\bar\eb)\circ \delta p$. Also the following abbreviations
will be used:
\begin{equation}\label{eq-S40}
  \begin{split}
    \ol{\pd}\Hop^0 &= \pd_\eb \Hop^0\circ\eeb{\bar\ub } + \pd_p \Hop^0\circ \bar p + \pd_\vphi \Hop^0\circ \bar \vphi\;,\\
\ol{\pd}\Kb^0 &= \pd_\eb \Kb^0\circ\eeb{\bar\ub } + \pd_p \Kb^0\circ \bar p + \pd_\vphi \Kb^0\circ \bar \vphi\;,\\
\end{split}
\end{equation}
in the generic sense (the example shown for the permeability). 
With this notation in hand we can introduce coefficients $\ol{\Hop}(\bar \ub, \bar p)$ involved in \eq{eq-rif1}, depending on
$(\bar \ub, \bar p)$, and on the response $(\ub\prevstep, p\prevstep)$ at the
previous time step $t\prevstep$, (recalling that the actual time level is $t = t\prevstep +
\Delta t$)
\begin{equation}\label{eq-S41}
\begin{split}
\ol{\Aop}(\bar\ub,\bar p) & = \Aop^0 + \ol{\pd}\Aop^0 + \pd_\eb (\Aop^0\eeb{\bar\ub })\circ \circpare{~} - \pd_\eb (\Bb^0 \bar p)\circ \circpare{~} + \pd_\eb \ul{\Hb}\bar{\ul{\vphi}}\circ \circpare{~},
\\
\ol{\Bb}(\bar\ub,\bar p) & = \Bb^0 + \ol{\pd}\Bb^0
+ \pd_p (\Bb^0 \bar p)\circ \circparp{~} -\pd_p (\Aop^0\eeb{\bar\ub })\circ \circparp{~} -\pd_p (\ul{\Hb}^0\bar{\ul{\vphi}} )\circ \circparp{~},
\\
\ol{\Db}(\bar\ub,\bar p,\ub^0,p^0) & = \Bb^0 + \ol{\pd}\Bb^0 
+ \pd_\eb \Bb^0 :(\eeb{\bar\ub }-\eeb{\ub^0})\circ \circpare{~} \\
& \quad + \pd_\eb M^0 (\bar p - p^0)\circ \circpare{~} - \pd_\eb \ul{Z}^0 (\bar{\ul{\vphi}} - \ul{\vphi}^0)\circ \circpare{~},\\
\ol{M}(\bar\ub,\bar p,\ub^0,p^0) & =  M^0 + \ol{\pd}M^0 + \pd_p  M^0 (\bar p - p^0)\circ \circparp{~} \\
 & \quad + \pd_p \Bb^0 :(\eeb{\bar\ub }-\eeb{\ub^0})\circ \circparp{~} - \pd_p \ul{Z}^0 (\bar{\ul{\vphi}} - \ul{\vphi}^0)\circ \circparp{~}\;,\\
\ol{\Kb}(\bar\ub,\bar p) & = \Kb^0 + \ol{\pd}\Kb^0\;,\\
\ol{\Gb} & = \pd_\eb \Kb^0 (\nabla \bar p - \fb)\circ \circpare{~}\;, \\
\ol{\Qb} & = \pd_p \Kb^0 (\nabla \bar p - \fb))\circ \circparp{~}.
\end{split}
\end{equation}


%% file: pz_flow_simulations.tex
\section{Numerical simulations}\label{sec:numex}

In this section, we present the results obtained by finite element (FE)
two-scale numerical simulations of flows in piezoelectric porous structures. In
order to validate the two-scale model of the homogenized medium presented in
Section~\ref{sec-homog}, as the ``reference model'' we consider the one
governing responses of the heterogeneous medium. Accordingly, ``the reference
solutions'' are computed using the \textit{direct numerical simulations} (DNS)
of the non-homogenized medium discretized on a ``full'' FE mesh which is
generated as a periodic lattice of the ``scaled'' periodic unit cells $\veps_0
Y$ with a given $\veps_0$. Both the two-scale simulations using the Homogenized
model (\hmodel) and the DNS of the Reference model (\rmodel) response have been
implemented in the Python based finite element package {\it SfePy}: Simple
Finite Elements in Python \cite{Sfepy_2019} which, in general, is well suited
for solving problems with multiscale and multiphysical features.
\revVa{The code for the calculations in this section is published
in the repository \cite{Lukes_code}.}{}

The evolutionary problems for both the \hmodel and the \rmodel are solved using
time stepping algorithms where the time derivatives are replaced by the
backward finite differences. For the \hmodel, Algorithm~\ref{sec:nonlin-macro}
is implemented which includes the homogenized coefficients updating. For DNS
with the \rmodel, the reference configuration is being updated to capture
the effect of nonlinear deformation on the flow.
The following spatial FE
approximations using Lagrangian elements are employed in all calculations: the
displacement, the electric potential, and the fluid pressure fields are
approximated using the P1 (piecewise linear) elements, while the fluid velocity
is approximated using the P2 (piecewise quadratic) ones.
\revVa{The finite element \rmodel consists of four-node tetrahedrons with
$4 \times (u_1, u_2, u_3, \varphi) = 16$ degrees of freedom (DOFs) per element
and of ten-node tetrahedrons with $10 \times (w_1, w_2, w_3) + 4 \times p = 34$
DOFs per element. The macroscopic \hmodel is approximated by eight-node
hexahedrons with $8 \times (u_1, u_2, u_3, p) = 32$ DOFs per element and the
microscopic \hmodel employs the same approximation as the \rmodel.}

\subsection{Two-scale computations with the \hmodel}

Due to the treatment of the deformation-dependent homogenized coefficients
$\tilde\Hop$ explained in Section~\ref{sec:appHC}, the micro-problems
\eq{eq-mip1}-\eq{eq-mip4} and \eq{eq-A10c} for characteristic responses of the
(non-deformed) microconfiguration $\Mcal(Y)$, as well as the associated
sensitivities \eq{eq-S40}-\eq{eq-S41} are solved independently of the
macro-problem \eq{eq-S25TD}.

\begin{figure}[h!]
    \centering
    \includegraphics[width=0.9\linewidth]{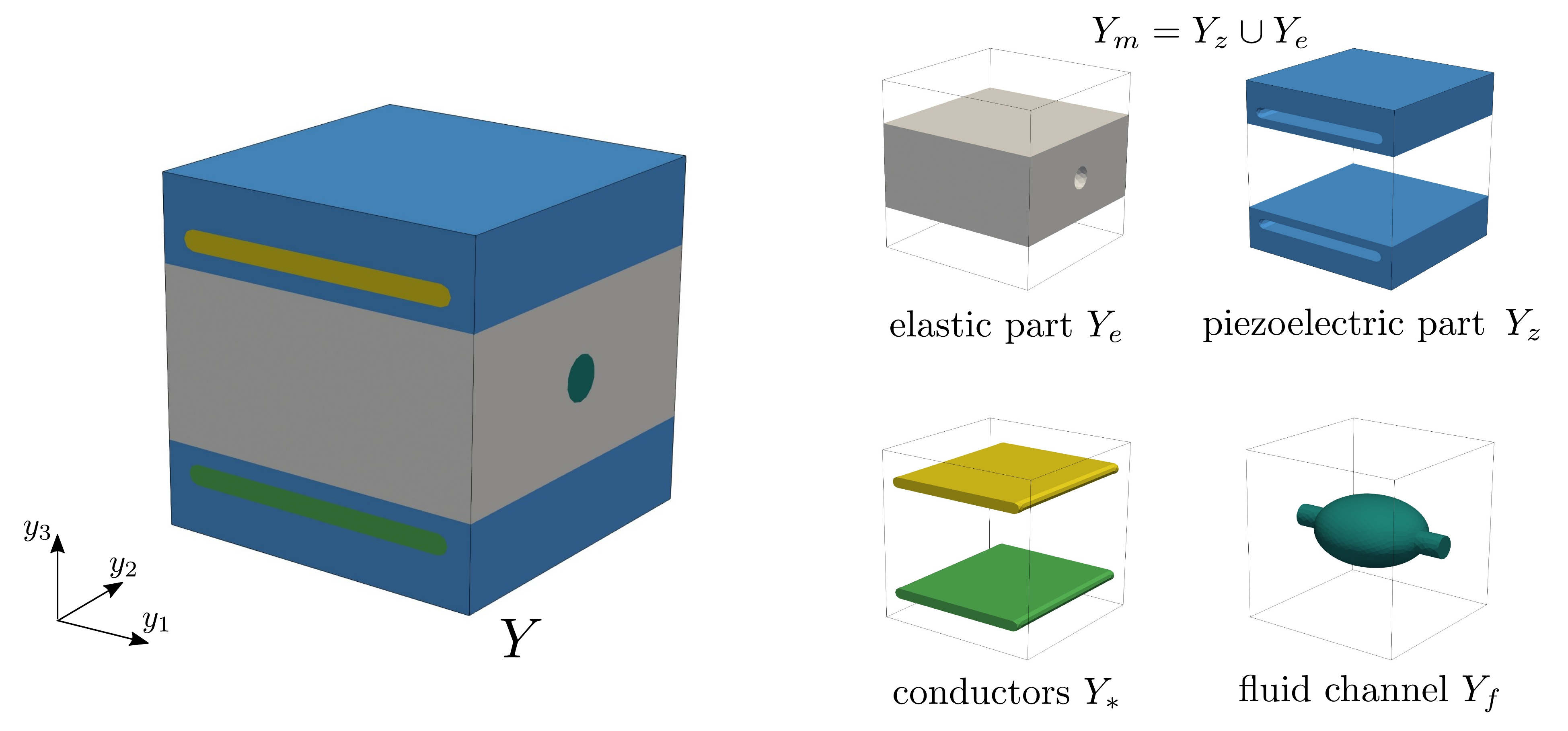}
    \caption{Decomposition of the microscopic periodic cell.}
    \label{fig:num-micro-geom}
\end{figure}

The microscopic domain $Y$ is decomposed into several parts according to
\eq{eq-6}: compliant elastic part $Y_e$, piezoelectric matrix $Y_z$ ($Y_e \cup
Y_z = Y_m$), two mutually disconnected conductors $Y_\ast^1$ and $Y_\ast^2$,
and fluid part $Y_f$, see Figure~\ref{fig:num-micro-geom}.
The material properties of the constituents are summarized in
Table~\ref{tab:num-material-properties}. Note that, by virtue of the scaling
introduced in Section~\ref{sec:mat}, the characteristic responses of the
microstructure are computed with the rescaled coefficients $\bar\parg, \bar\db$
and $\bar\eta$ obtained for a give microstructure size, \ie a given
$\veps:=\vepsz$.
\begin{table}
    \begin{center}
    \begin{tabular}{|l|c|}
        \hline
        {elastic part $Y_e$ (elastomer)} & E = 0.02\;GPa, $\nu$ = 0.49\\
        \hline
        \multirow{3}{*}{\shortstack{piezoelectric part $Y_m$ \\ (piezo-polymer, \cite{koutsawa_belouettar_makradi_nasser_2010})}} &
            $\Aop = \begin{bmatrix}
                     6.0& 3.72& 3.83& 0& 0& 0 \\
                     3.72& 6.0& 3.83& 0& 0& 0 \\
                     3.83& 3.83& 20.3& 0& 0& 0 \\
                     0& 0& 0& 1.23& 0& 0 \\
                     0& 0& 0& 0& 1.23& 0 \\
                     0& 0& 0& 0& 0& 1.23
                 \end{bmatrix} \cdot 10^7$\;Pa\\
         & $\parg^\veps = \begin{bmatrix}
                                0& 0& 0& 0& 0.01& 0\\0.34
                                0& 0& 0& 0& 0& 0.01\\
                                -0.09& -0.09& 5.91& 0& 0& 0
                            \end{bmatrix}$\;C/m$^2$\\
         & $\db^\veps = \begin{bmatrix}
                                18& 0& 0\\
                                0& 18& 0\\
                                0& 0& 255.3
                            \end{bmatrix} \cdot 8.854 \cdot 10^{-12}$\;C/Vm\\
        \hline
        {conductors $Y_*$ (steal alloy)} & E = 200\;GPa, $\nu$ = 0.25\\
        \hline
        {fluid $Y_f$ (water)} & $\gamma = 1.0 / (2.15 \cdot 10^9)$\;Pa, $\mu^\veps = 8.9 \cdot 10^{-4}$\;Pa\;s\\
        \hline
    \end{tabular}
    \end{center}
    \caption{Material properties of the piezoelectric skeleton, Voigth representation of the strain is employed for tensors $\Aop$ and $\parg^\veps$. Note that $\bar\parg = \bar\parg^\veps/\veps$ and $\bar\db = \db^\veps/\veps$ are computed for a given scale $\veps:=\veps_0$.}
    \label{tab:num-material-properties}
\end{table}

In order to compare the \hmodel solutions with the corresponding \rmodel
responses obtained by the DNS we consider a ``pseudo 1D problem'', such that,
although $\Om$ is a 3D block, all the macro-problem parameters, including the
boundary conditions make the response of the \hmodel to depend on the $x_1$
coordinate only.

\begin{figure}[h!]
    \centering
    \includegraphics[width=0.8\linewidth]{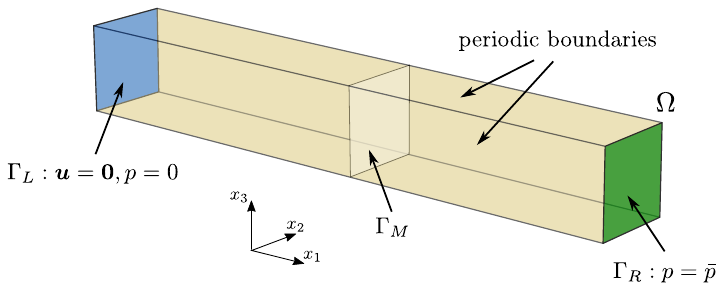}
    \caption{Macroscopic domain $\Omega = ]0,L[ \times ]0,a[^2$ and applied boundary conditions. The positions of sections $\Gamma_L$, $\Gamma_M$, and $\Gamma_R$ are $x_1 = 0, L/2$ and $L$, respectively.}
    \label{fig:num-macro-geom-bc}
\end{figure}

The macroscopic problem \eq{eq-M1}-\eq{eq-M2} recast for the semilinear
formulation \eq{eq-rif1}-\eq{eq-S25TD}, is imposed in $\Omega$ represented as
block with dimensions $L\times a \times a$ spanned on axes $x_1,x_2,$ and $x_3$
in the respective order. The boundary conditions, see
Fig.~\ref{fig:num-macro-geom-bc}, are considered according to \eq{eq-H*S2},
being modified on $\Gamma_0$ by periodic conditions. They are applied in the
directions $x_2$, $x_3$ on all four faces with the normals aligned with either
of these directions. In particular, on the faces $\Gamma_1 \equiv \Gamma_L$ and
$\Gamma_2 \equiv \Gamma_R$ perpendicular to $x_1$ axis, we prescribe:
$\ub=\bmi{0}$, $p=0$ on $\Gamma_L$, whereas $p= \bar p$ on $\Gamma_R$, with
$\bar p$ [Pa] a given pressure value; this defines also the surface traction
$\bbav = -\bar p \nb$. The boundary conditions are depicted in
Fig.~\ref{fig:num-macro-geom-bc}

While $\vphi^1 = 0$ is prescribed on the first electrode, the potential waves propagate due to the second electrode, whereby the electric potential  $\varphi^2$ to be a function depending on time $t$
and coordinates $x_1 \in ]0,L[$ and $x_2\in]0,a[$:
\begin{equation}\label{eq:num_electric_potential}
  \begin{split}
    \psi(x_1,x_2,t) & := b_1\,x_1 + b_2\,x_2 - c\,t + d \;,\\
    \varphi^2(x_1, x_2, t) & =
    \begin{cases}
        \frac{1}{2}\left[1 + \cos(\psi(x_1,x_2,t)+ \pi) \right]\epot,
            \quad \mbox{for}\quad \psi(x_1,x_2,t) < 0\\
        0 \quad \mbox{otherwise}
    \end{cases},
 \end{split}   
\end{equation}
where $b_1$, $b_2$, $c$, and $d$ are constants specified for the examples
reported below as ``1D'' and ``2D'' cases, and the voltage wave amplitude
$\epot$ is adjusted accordingly to the case, see
Tab.~\ref{tab:phi}. The above setting ensure the macroscopic
displacements $\ub$ and the pressure gradient $\nabla_x p$ induced by the
electric potential to be non-zero only in $x_1$ direction, thus, the ``pseudo
1D problem'' intentions are met. Due to the periodicity in the $x_2$ and $x_3$
directions, the length $a$ can be arbitrary, but for a better comparison with
the reference calculation due to the \rmodel, we choose $a := L/N$, where $N$
is the number of the real-sized periods $Y^{\veps_0} = \veps_0 Y$, see the next
subsection.

\begin{table}
    \begin{center}
      \begin{tabular}{|l|lllll|}
        \hline
        Structure / flow & $\epot$ [V]& $b_1$ &  $b_1$ & $c$ & $d$ \\
        \hline
        case 1D & $4\cdot10^5$ & $\pi / 0.03$ & $0$ & $10\pi$ & $0$ \\
        case 2D & $4\cdot10^5$ & $\pi / 0.03$ & $\pi/ 0.075$ & $10\pi$ & $0$ \\
        \hline
      \end{tabular}
      
      \caption{Parameters of the electric potential wave, see \eq{eq:num_electric_potential}.}\label{tab:phi}
    \end{center}
\end{table}

\subsection{Reference simulation with the \rmodel}

For a validation of the homogenized model and verification of its numerical
implementation the approach of the direct numerical simulations (DNS) of the
heterogeneous continuum response is employed. The reference response of the
\rmodel is obtained as the solutions of problem \eq{eq-pzfl2} with the material
parameters given Tab.~\ref{tab:num-material-properties}. In connection with the
\hmodel, we set $L = 0.1$ m and $a = L/N = 0.005$ m with $N = 20$, which yields
$\veps_0 = 0.005$. The general decomposition \eq{eq-pzm1} holds, where in the
solid phase $\Om_s^\vepsz = \Om_e^\vepsz\cup\Om_z^\vepsz\cup\Om_*^\vepsz$ we
distinguish an elastic-dielectric part $\Om_e^\vepsz$, the piezoelectric
segments $\Om_z^\vepsz$ and the electrodes
$\Om_*^{\kkk,\vepsz} \subset \Om_*^\vepsz$, $\kkk = 1,2$,
being generated by periodically repeated copies of the reference cell $Y$
components, see Fig.~\ref{fig:num-micro-geom}. Using
\eq{eq:num_electric_potential} with the constants given in Tab.~\ref{tab:phi},
the electric potential wave (the voltage) is imposed at $\Om_*^{2,\vepsz}$. To
adhere with finite size electrodes, the unfolding brackets are employed: $\hat
x_k = \vepsz\{x_k/\vepsz\}_Y$ and \eq{eq:num_electric_potential} is employed
with $\varphi^2(\hat x_1, \hat x_2, t)$, attaining piecewise constant values.


%
\begin{figure}[h!]
    \centering
    \includegraphics[width=0.9\linewidth]{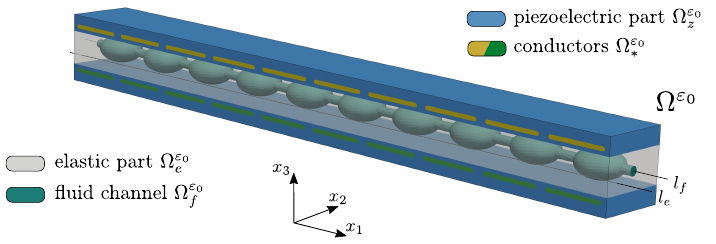}
    \caption{The geometry employed in the reference numerical simulation
        and two probe lines $l_f$, $l_e$ along which the values
        in Section~\ref{sec:num-validation} are displayed.}
    \label{fig:num-dir-geom}
\end{figure}

To capture the model nonlinearity induced by imposing the \rmodel equations in
the ``spatial'' configuration deforming in time, we chose updating the geometry
of the subdomains of $\Om^{\varepsilon_0}(t)$ directly. In this way, we avoid any additional
approximation as the one employed in the \hmodel, being based on the
sensitivity analysis of the \hmodel effective coefficients. Instead,
we implemented an iterative incremental computational algorithm with updates of
the FE mesh used to discretize the subdomains on $\Om{\varepsilon_0}(t)$. In each iteration
within a single time step, the nodes of the FE mesh are moved to new positions
defined by the actual displacements. Such an iterative updating loop is
repeated until the geometry updates are negligible, \ie the iterations
terminate when the norm of the displacement corrections is small enough.

\subsection{Validation test}\label{sec:num-validation}

We consider the ``pseudo 1D problem'' imposed in the block specimen $\Om$, with
the boundary conditions explained above, however, setting the fluid pressure to
be zero on both the ends $\Gamma_L$ and $\Gamma_R$ of the specimen, $\bar p =
0$. Responses of both the \reflabel- and \homlabel-models are computed for the actuation by
$\bar\vphi^2 = 4\cdot10^5$\,V according to \eq{eq:num_electric_potential} and
the constants adjusted for the ``1D'' case, see Tab.~\ref{tab:phi}, applied in
50 uniform time steps in the time interval $[0, 1]$\,s. Note that $b_2 = 0$ so
that $\varphi^2 = \varphi^2(x_1, t)$, whereby the above mentioned modification
due to the finite size electrodes applies in the \rmodel.





\begin{figure}[h]
    \centering
    \includegraphics[width=0.49\linewidth]{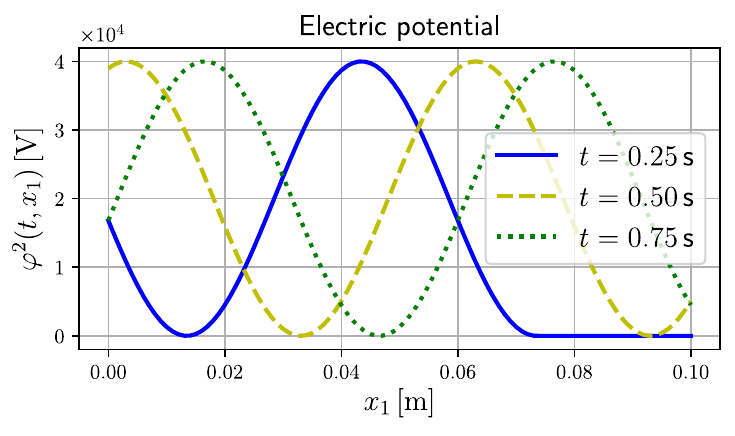}
    \caption{Spatial distribution of the control electric potential $\varphi^2$
    at selected times $t_k$; note the step-wise modification applies in the
    \rmodel.} \label{fig:num-electric_potential}
\end{figure}

\begin{figure}[h]
    \centering
    \includegraphics[width=0.99\linewidth]{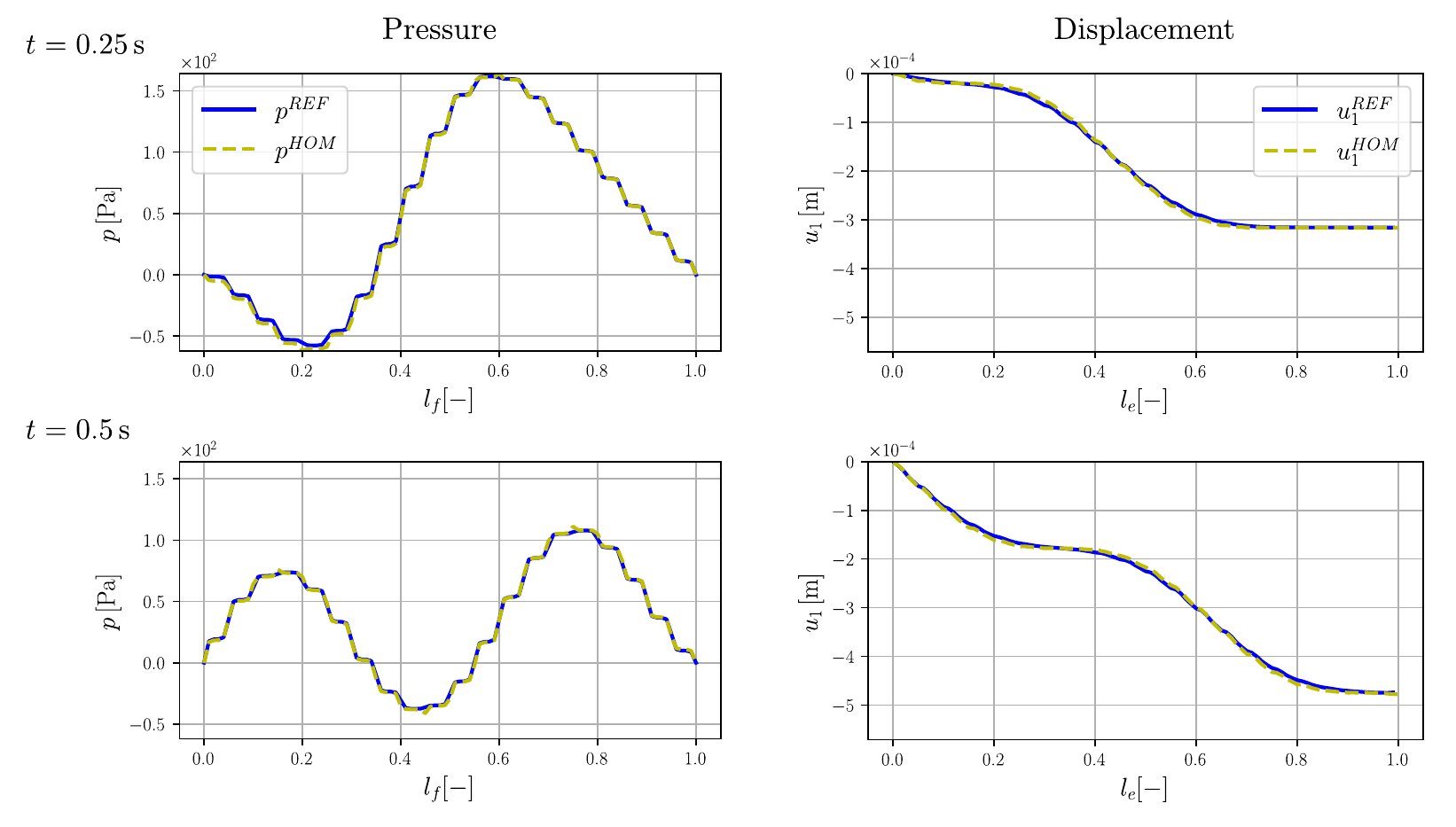}
    \caption{Comparison of homogenized and reference solutions -- linear model.}
    \label{fig:num-hom_x_ref-fields-lin}
\end{figure}

First we compare the linear \hmodel defined in \eq{eq-S25}, for which the homogenized coefficients
$\Hop^0$ are given by the non-perturbed configuration. Accordingly, the \rmodel
is used without updating the configuration, thus, solutions are computed on the
fixed FE mesh To provide consistent responses corresponding to the finite
heterogeneity captured by the \rmodel, the macroscopic solutions of the \hmodel 
must be reconstructed at the microscopic level for a given finite scale
$\varepsilon^0$. The two-scale field reconstruction (sometimes called the
``downscaling'') is based on the decomposition formulae \eq{eq-L2b}, details are
discussed in \cite{Rohan-Lukes-PZ-porel}. The reconstructed pressure field
$p^\Hveps$ and the displacement field $u_1^\Hveps$ together with the reference
solutions $p^{\reflabel,\vepsz}$, $u_1^{\reflabel,\vepsz}$ at selected time levels
$t_k = 0.25, 0.5\deleted{, 0.75}$\,s are compared in Fig.~\ref{fig:num-hom_x_ref-fields-lin}. Distributions
of these fields are displayed along lines $l_f$, $l_e$ parallel to the axis
$x_1$, see Fig.~\ref{fig:num-dir-geom}. The electric potential
$\varphi^2(x_1,t_k)$ at the three selected times $t_k$ is shown in
Fig.~\ref{fig:num-electric_potential}.
\revVa{It can be seen from Fig.~{\ref{fig:num-hom_x_ref-fields-lin}}
that the agreement between the reference and the homogenized model for the linear
case is very good.}{}

\begin{figure}[h!]
    \centering
    \includegraphics[width=0.99\linewidth]{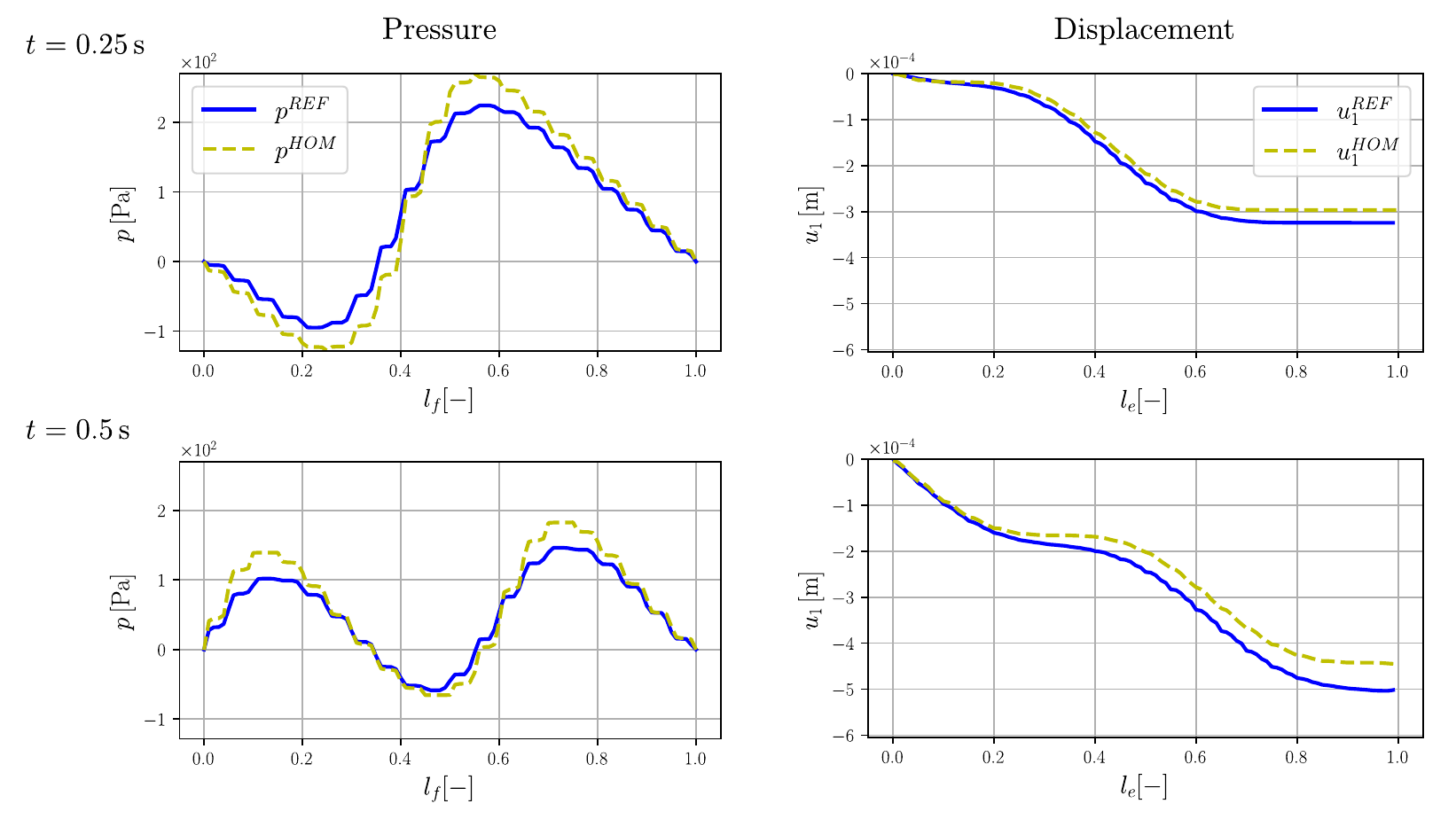}
    \caption{Comparison of homogenized and reference solutions -- nonlinear model.}
    \label{fig:num-hom_x_ref-fields-nlin}
\end{figure}

\begin{figure}[h!]
    \centering
    \includegraphics[width=0.99\linewidth]{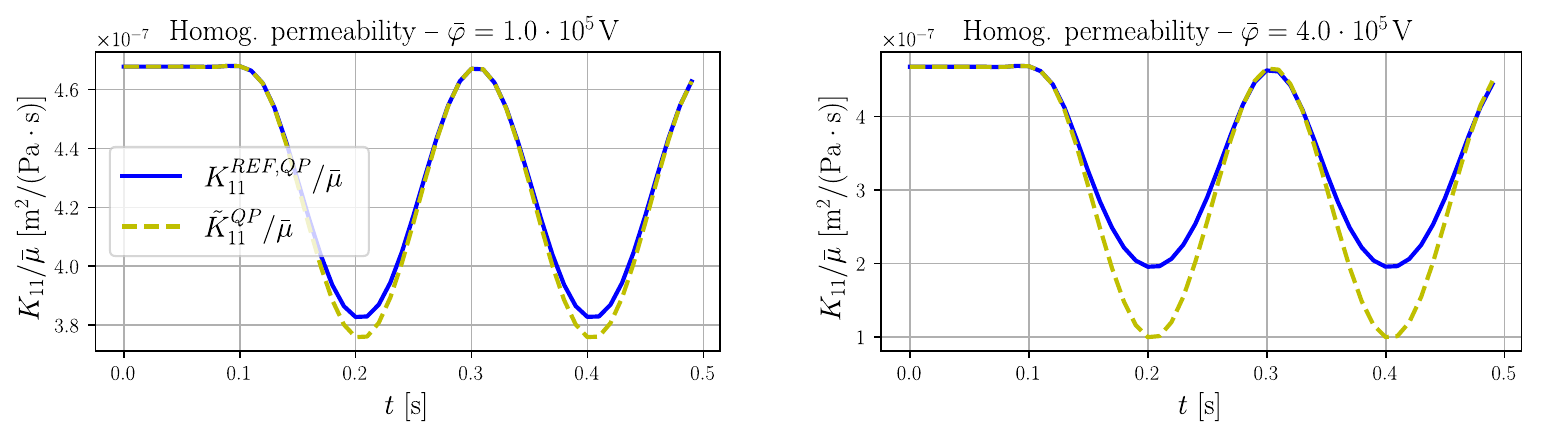}
    \caption{
    \revVa{
        Comparison of homogenized and reference permeability ${K}_{11}$
        at macroscopic quadrature point QP:
            $\bar\varphi = 10^5$\,V -- left,
            $\bar\varphi = 4\cdot 10^5$\,V -- right.
    }
    }\label{fig:num-hom_x_ref-nlin_Kx}

\end{figure}

\begin{figure}[h!]
    \centering
    \includegraphics[width=0.99\linewidth]{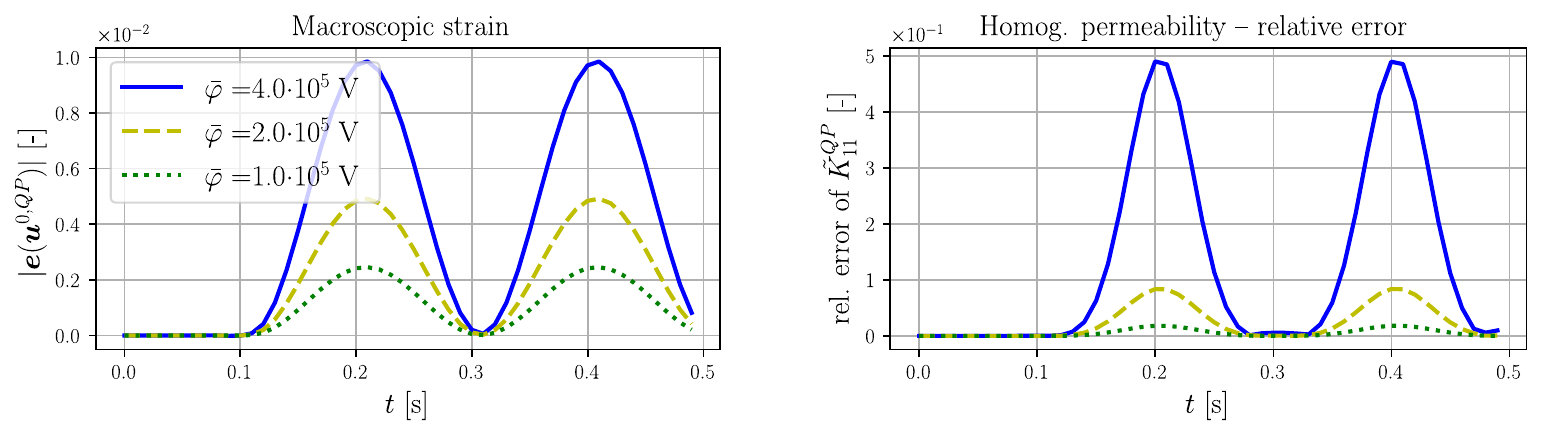}
    \caption{
    \revVa{
        Dependence of the quantities on the applied electric potential:
        strain magnitude $\vert \eb(\ub^{0, QP})\vert$ -- left,
        relative error of $\tilde{K}_{11}^{QP}$ -- right.
    }
    }\label{fig:num-hom_x_ref-nlin_K}
\end{figure}

In the second comparison, we consider nonlinear models so that the homogenized
coefficients featuring the \hmodel are updated according to \eq{eq-S24},
whereby the macroscopic response is given by the algorithm described is
Sec.~\ref{sec:nonlin-macro}. Accordingly, the \rmodel responses are computed
with the mesh updates in an iterative loop embedded in each time step, as
described above. Results of this second comparison test are presented as in the
first one; the reconstructed pressure and displacement fields with that
obtained by the reference model are reported in
Fig.~\ref{fig:num-hom_x_ref-fields-nlin}.
\revVa{The disagreement of the homogenized model with the reference model for the
nonlinear case is quite significant. This discrepancy is due
to the fact that the applied electric potential induces relatively large
strains on the order of $10^{-2}$ and the approximation of the homogenized
coefficients by sensitivity analysis is no longer accurate. The dependence of
the permeability on the applied electric potential and hence on the deformation
is shown in Fig~\ref{fig:num-hom_x_ref-nlin_Kx} and the relative error
of the permeability is plotted in Fig~\ref{fig:num-hom_x_ref-nlin_K}/right.
The homogenized permeability $\tilde\Kb^{\QP}$, which is updated by sensitivity analysis, 
is evaluated at a given macroscopic
quadrature point $\QP$ and the reference permeability tensor $\Kb^{REF, \QP}$ is
computed on the deforming microscopic cell $Y$. The cell deformation is defined
by
\begin{equation}\label{eq-num-micro-def}
    \delta u(\QP, y) = (\omegabf^{ij} + \Pibf^{ij})e^x_{ij}(\delta\ub^{0, \QP})
    - \omegabf^P \delta p^{0, \QP}
    + \sum_\alpha \delta\hat\omegabf^\alpha \hat\varphi^{\alpha, \QP},
\end{equation}
where $\delta u^{0, \QP}$, $\delta p^{0, \QP}$, and $\delta \hat\varphi^{\alpha, \QP}$ 
are the increments of the macroscopic variables evaluated at $\QP$.
}{}

\begin{figure}[h!]
    \centering
    \includegraphics[width=0.99\linewidth]{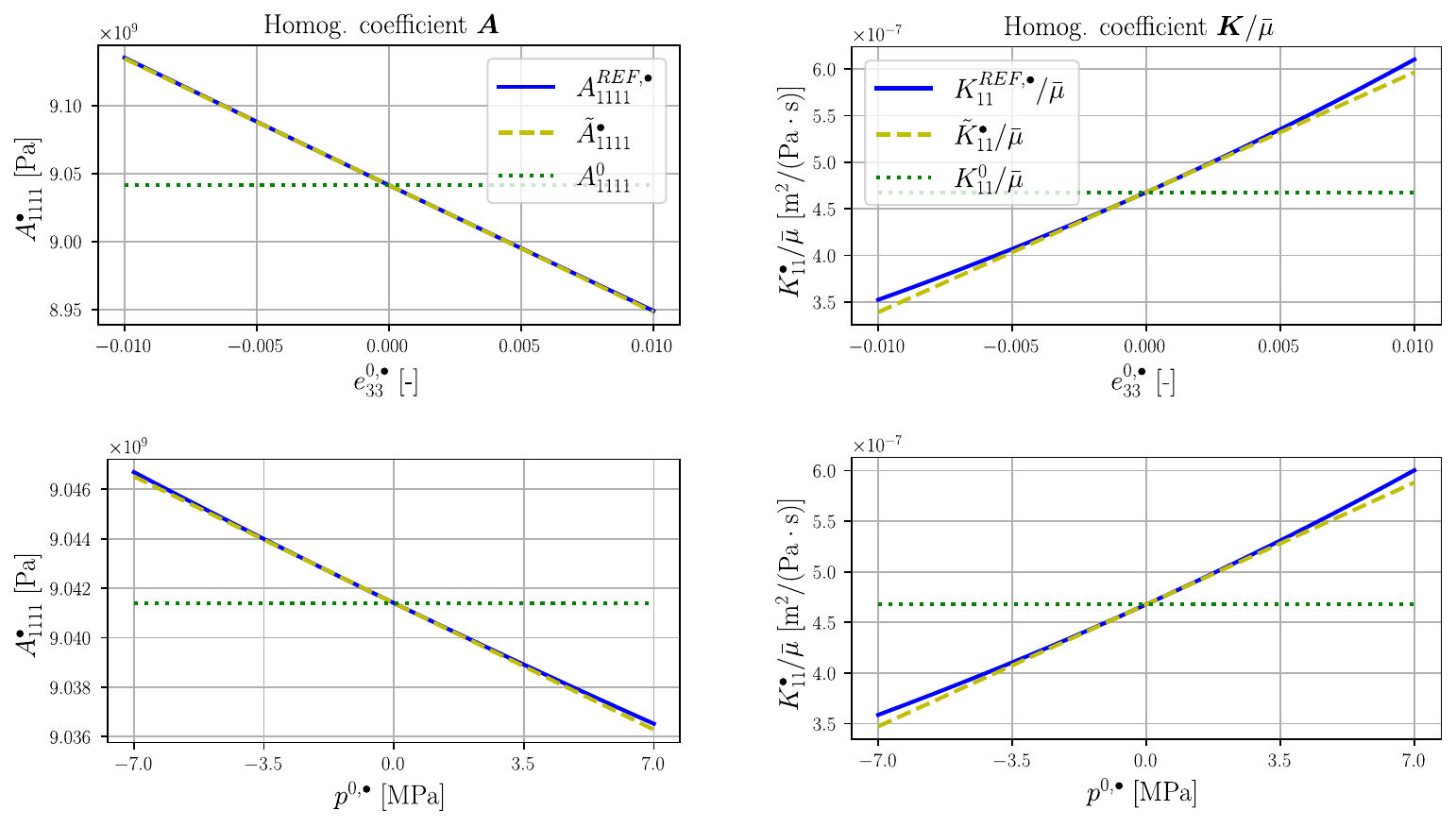}
    \caption{
    \revVa{
        Dependence of the homogenized coefficients $\tilde{A}_{1111}$ (left)
        and $\tilde{K}_{11}$ (right) on the macroscopic quantities.
        The unit cell $Y^\bullet$ is loaded either by the macroscopic strain
        $e^{0,\bullet}_{33} \in \left[-0.01, 0.01\right]$ (top)
        or by the macroscopic pressure $p^{0,\bullet} \in \left[-7, 7\right]$\,MPa (bottom).
        The coefficients approximated by sensitivity analysis are compared with the reference values.
    }
    }\label{fig:num-hom_x_ref-nlin_AK}
\end{figure}

\revVa{
To show the dependence of the homogenized coefficients on macroscopic
deformation and pressure, we perform a test on a single reference cell,
denoted $Y^\bullet$. The cell is deformed by the macroscopic strain
$e^{0,\bullet}_{33} \in \left[-0.01, 0.01\right]$ in the first test and by the macroscopic pressure
$p^{0,\bullet} \in \left[-7, 7\right]$\,MPa the second test according to equation \eq{eq-num-micro-def}.
While the elastic coefficient $\tilde{A}_{1111}^{\bullet}$ changes for the prescribed strain
only on the order of one percent, the homogenized permeability $\tilde{K}_{11}^{\bullet}$ changes by about 25 percent.
For the varying macroscopic pressure, the elastic coefficient remains almost unchanged,
whereas the change in the homogenized permeability is again in the order of tens of
percent. Fig.~\ref{fig:num-hom_x_ref-nlin_AK} also presents the comparison of the coefficients
obtained by the sensitivity analysis (denoted by $\tilde{}$\,) and those calculated on the
deforming microscopic cell (denoted by ${}^{REF}$).
}{}

\begin{figure}[h!]
    \centering
    \includegraphics[width=0.98\linewidth]{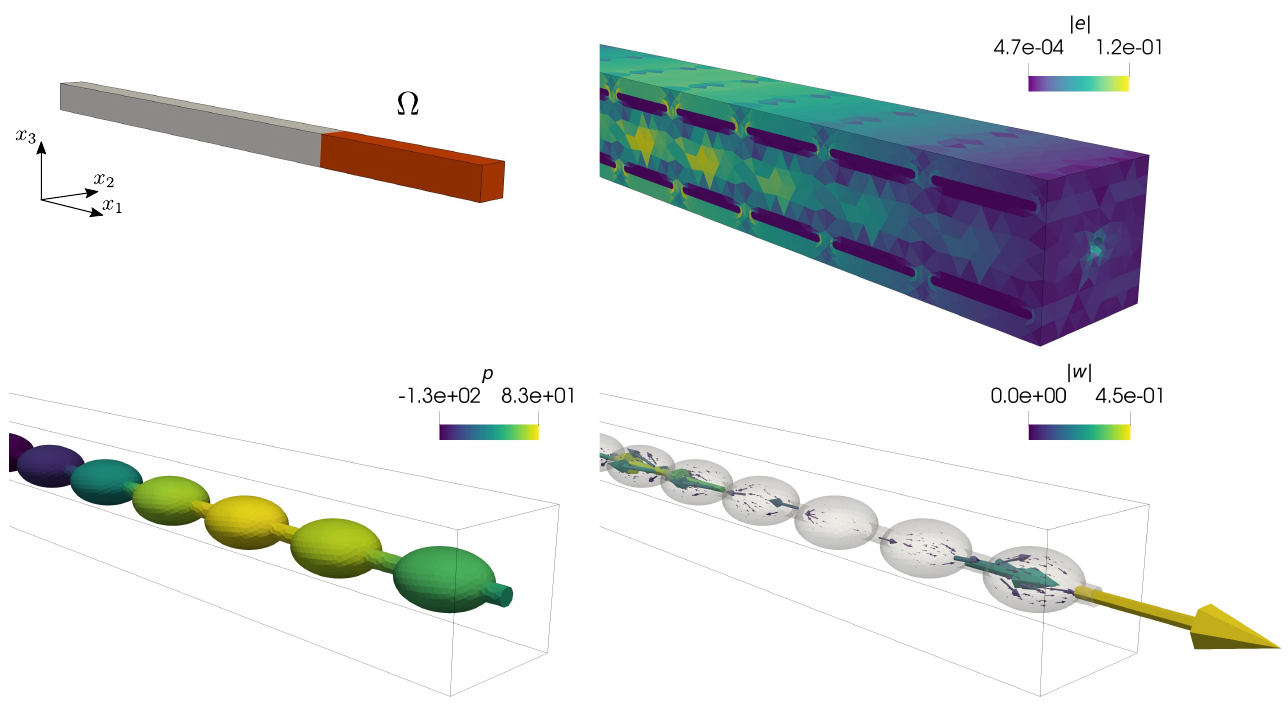}
    \caption{Reconstructed homogenized solution in the highlighted (red) part of $\Omega$ at time $t = 0.9$\,s:
        strain magnitude $|\eeb{\ub^\Hveps}|$ -- top right,  pressure $p^\Hveps$ -- bottom left, velocity -- bottom right.}
    \label{fig:num-hom_reconstructed_fields}
\end{figure}

The reconstructed microscopic pressure, fluid velocity, and strain magnitude
for a part of the specimen are depicted in
Fig.~\ref{fig:num-hom_reconstructed_fields}.

\begin{figure}[h!]
    \centering
    \includegraphics[width=0.98\linewidth]{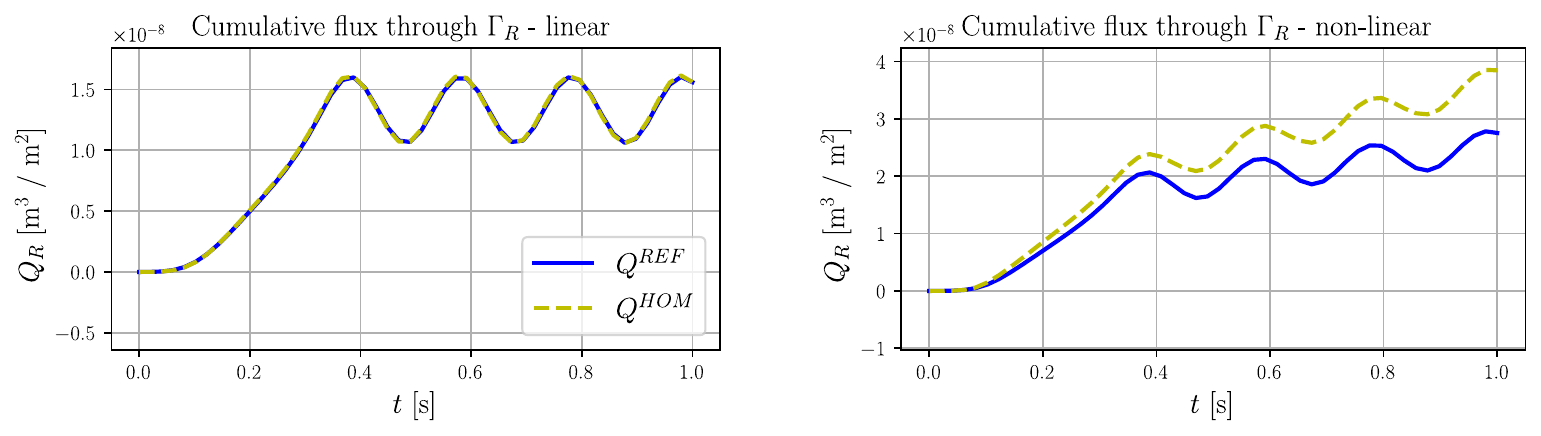}
    \caption{Cumulative flux through the \revVr{right faces $\Gamma_R$}{internal faces $\Gamma_M$}:
        linear model -- left, nonlinear model -- right.}
    \label{fig:num-hom_x_ref-flux}
\end{figure}

To illustrate the effect of the model nonlinearity in the simulations, the
cumulative fluxes through the \revVr{right face $\Gamma_R$}{middle face $\Gamma_M$} computed as in
\eq{eq-1d-8}, but with $w:= w_1$ (fluxes aligned with the $x_1$ axis), are
compared in Fig.~\ref{fig:num-hom_x_ref-flux} for the linear and nonlinear
cases, \ie without and with updating the homogenized coefficients using the
approximation \eq{eq-S24}. Positive flux values mean the flow from $\Gamma_L$
to $\Gamma_R$. While in linear case the flux oscillates around \revVr{a constant value}{zero} indicating
no effective flux, in the nonlinear case, the pumping effect is obvious since
the cumulative flux values increase with time, so that the fluid is transported
from $\Gamma_L$ to $\Gamma_R$.
\revVa{Note that the results are calculated for the actuation $\epot = 4\cdot
10^5$\,V, i.e. for relatively large strains, and therefore we observe a big difference 
in the cumulative flux for the nonlinear case. If we decrease $\epot$, the difference becomes smaller.}{}

\revVa{For the number of periods $N=20$, the \rmodel has 317265
tetrahedral elements, of which 22328 are for-node tetrahedrons and 83980 are
ten-node tetrahedrons, for a total of more than $3\cdot 10^6$ DOFs.}



\subsection{Peristaltic pumping against natural flow}

\begin{figure}[h!]
    \centering
    \includegraphics[width=0.98\linewidth]{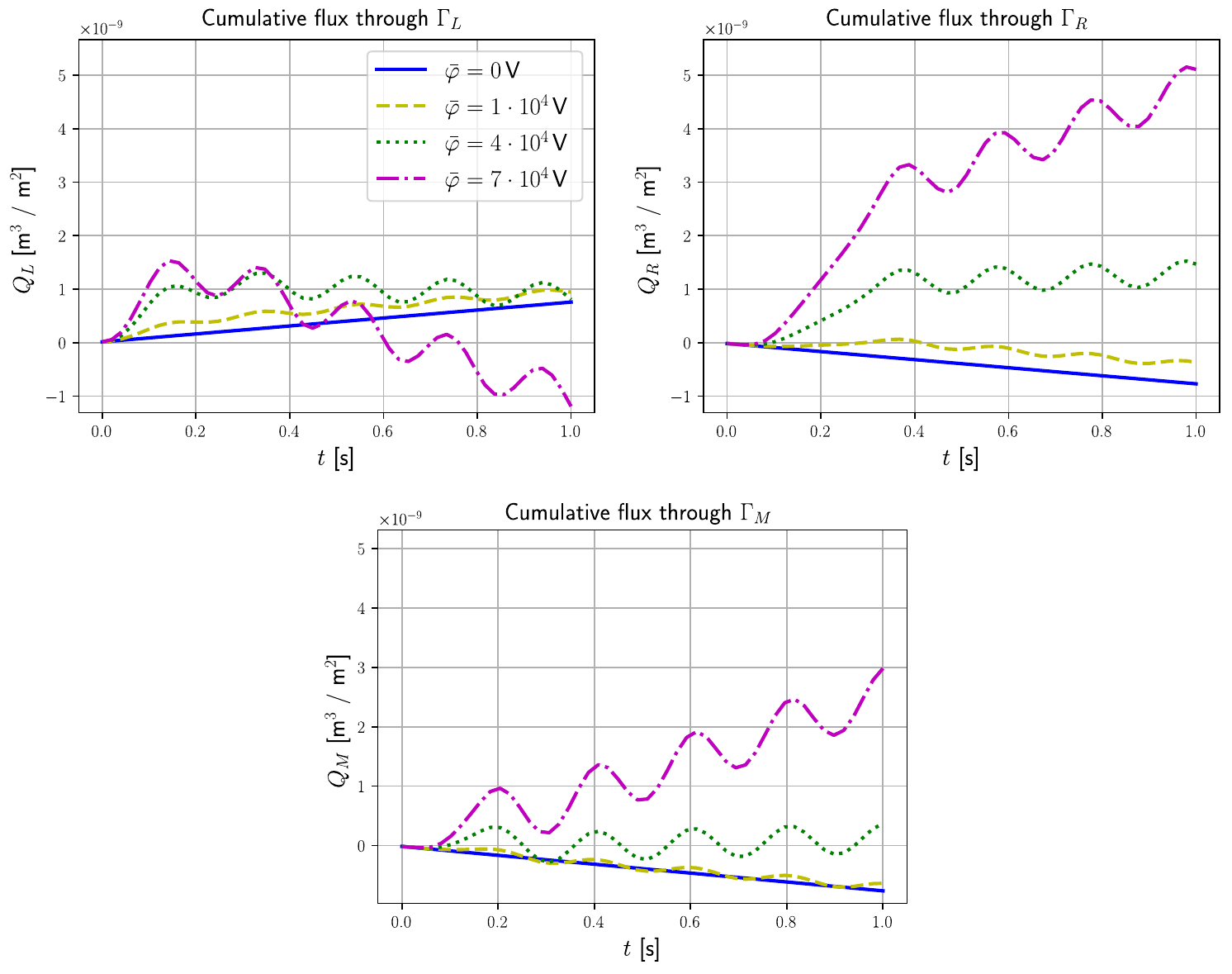}\hfill
    \caption{Cumulative flux through faces $\Gamma_L$, $\Gamma_R$, and $\Gamma_M$ for various $\epot$.}
    \label{fig:num-hom_flow_transport}
\end{figure}

Now we examine the peristaltic pumping effect in the context the linear and
nonlinear models. For this, modified data are considered. Given non-zero
pressure $\bar p$ on $\Gamma_R$, whereas $p = 0$ on $\Gamma_L$, the Darcy low
induces a natural flow towards $\Gamma_L$. The pumping actuated by the induced
peristaltic deformation wave running to the right should produce a net flux
towards $\Gamma_R$, just in the opposite direction than the one of the natural
flow. We choose $\bar p = 1$\,kPa, and now put $\veps^0 = 10^{-3}$, which
means different material properties than those considered in the above
validation tests. The aim is to assess the influence of
the electric voltage actuation $\bar\vphi^*$ on the amount of the fluid
transported to the right $\Gamma_R$, see Fig.~\ref{fig:num-hom_flow_transport}.
For the electric potential $\epot = 0$\,V, obviously, there is no
peristalsis driven flow, the natural flow towards $\Gamma_L$ due to the
pressure gradient is obtained. By increasing $\epot$, the reverse flow
can be achieved and for $\epot = 7 \cdot 10^4$\,V a stable fluid
transport due to the peristaltic pumping occurs.

\begin{figure}[h!]
    \centering
    \includegraphics[width=0.5\linewidth]{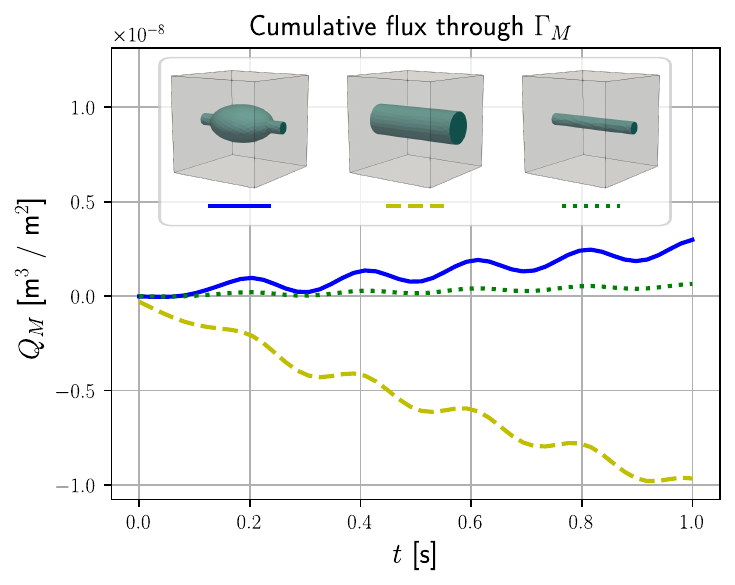}\hfill
    \caption{Cumulative flux through face $\Gamma_M$ for three different
             fluid channel shapes, $\epot = 7 \cdot 10^4$\,V.}
    \label{fig:num-hom_flow_transport_geom}
\end{figure}

The ability of the structure to transport fluid against pressure drop is
strongly affected by the geometrical arrangement of the microstructure,
especially by the shape of the fluid channel. This effect is illustrated in
Fig.~\ref{fig:num-hom_flow_transport_geom} where the cumulative fluxes for the
``balloon $\cap$ cylinder'' fluid channel geometry and the simple cylindrical
geometries are compared. As can be seen from the figure, purely cylindrical
fluid channels are not suitable for peristaltic fluid transport because
cylinders with large diameters cannot be sufficiently throttled to eliminate
reflux, while the smaller channels can only transport a limited amount of fluid
for a given pressure gradient and electric potential.

\subsection{Peristaltic flow in ``2D'' structure}

\begin{figure}[h!]
    \centering
    \includegraphics[width=0.98\linewidth]{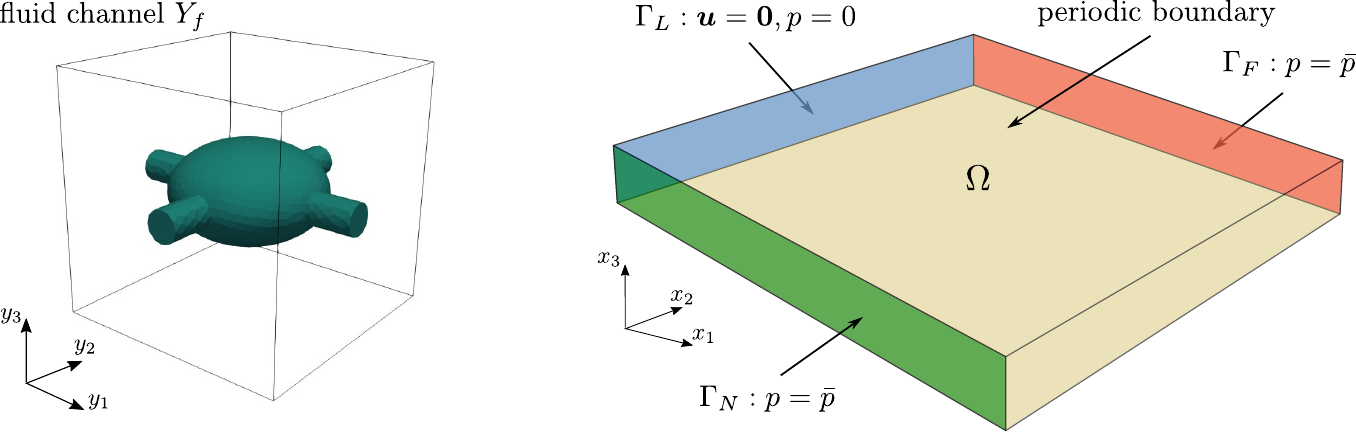}
    \caption{Microscopic periodic cell $Y$ (left) and macroscopic domain 
        $\Omega$ with an illustration of the applied boundary conditions (right).}
    \label{fig:num-macro2D-geom-bc}
\end{figure}

In this section, we employ a modified periodic unit cell which allows for flows
in the $x_1$ and $x_2$ directions, see Fig.~\ref{fig:num-macro2D-geom-bc}/left,
and the macroscopic sample of dimensions $0.1 \times 0.1 \times 0.0625$\,m to
which the boundary conditions depicted in
Fig.~\ref{fig:num-macro2D-geom-bc}/right are applied.

\begin{figure}[h!]
    \centering
    \includegraphics[width=0.98\linewidth]{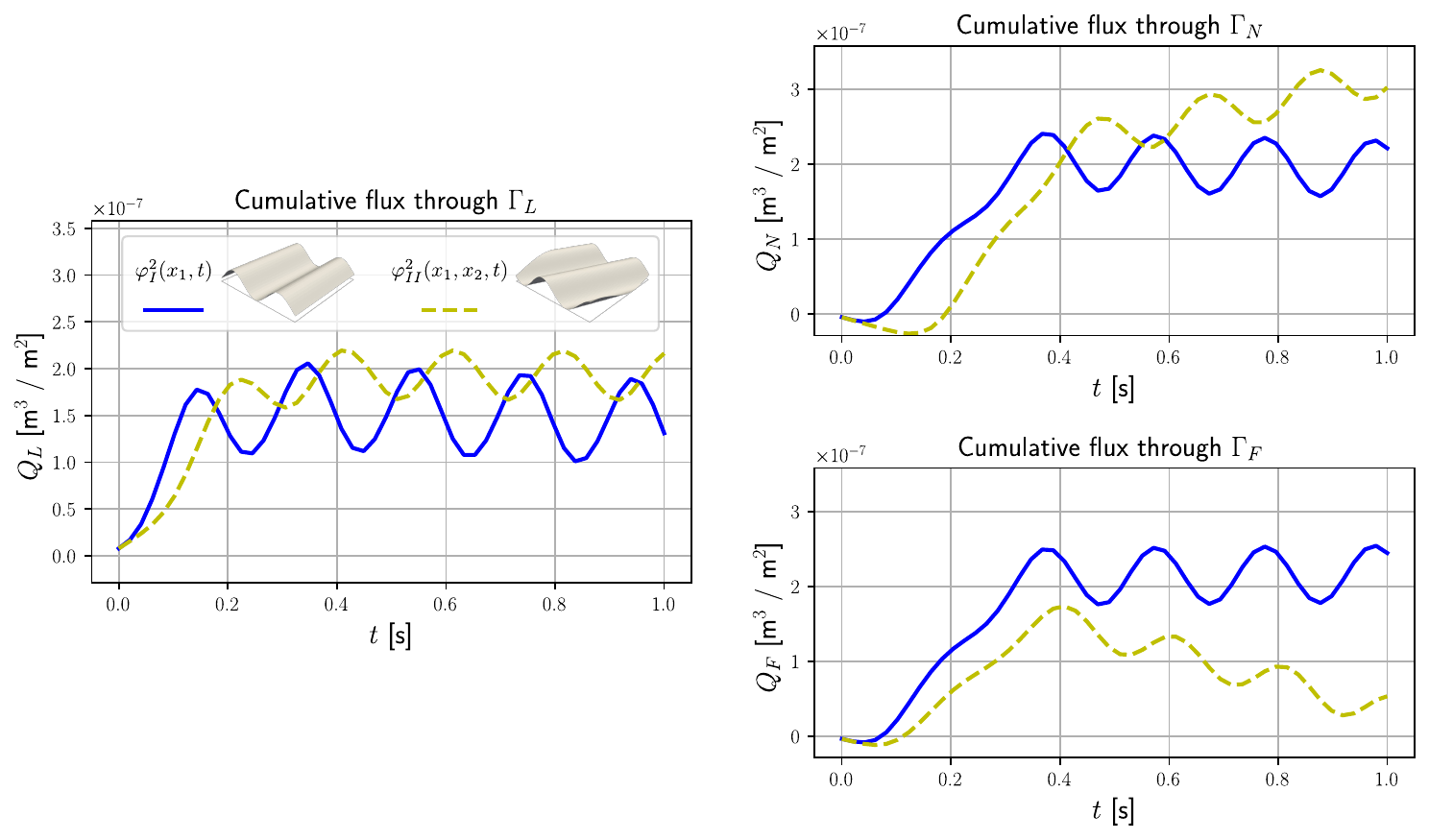}
    \caption{Cumulative fluxes through faces $\Gamma_L$, $\Gamma_N$, $\Gamma_F$
    for the electric potential given by $\varphi^2_{I}$ and $\varphi^2_{II}$.}
    \label{fig:num-macro2D-Q}
\end{figure}

In order to demonstrate how the fluid transport through the structure can be
controlled by the prescribed electric potential, we consider two distinct
functions $\varphi^2$. In the first simulation, we employ the function
$\varphi^2 = \varphi^2_{I}(x_1, t)$ identical to that used in
Sec.~\ref{sec:num-validation} while the second simulation is performed for
non-zero $d=\pi / 0.75$ and $b_2=\pi/ 0.075$, the ``case 2D'' in
Tab.~\ref{tab:phi}, so that $\varphi^2 = \varphi^2_{II}(x_1,x_2,t)$. The
comparison of the cumulative fluxes through the faces $\Gamma_N$, $\Gamma_F$,
$\Gamma_L$ for $\varphi^2_{I}$ and $\varphi^2_{II}$ is shown in
Fig.~\ref{fig:num-macro2D-Q}.
\revVa{The figures show that changing the parameters of the control
electric potential can have a major effect on the resulting fluid transport and
ability to pump in the desired direction.
While for $\varphi^2_{I}$ there is a uniform distribution of flux from
$\Gamma_L$ to $\Gamma_F$ and $\Gamma_N$, in the case of $\varphi^2_{II}$ we
altered the flux by changing the parameter $d$ so that fluid flows from
$\Gamma_L$ and $\Gamma_F$ to $\Gamma_N$.
Optimizing the control function parameters with respect to the
required flow is a challenge for further research.}


\begin{figure}[h!]
    \centering
    \includegraphics[width=0.98\linewidth]{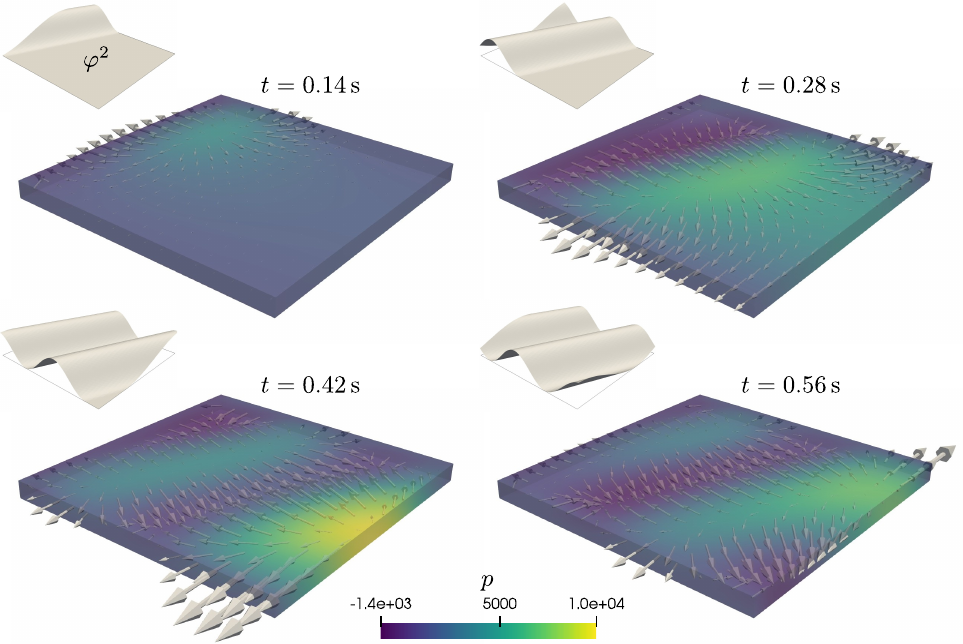}
    \caption{Macroscopic pressure and velocity fields at
        times $t_k=0.14, 0.28, 0.42, 0.56$\,s, $\varphi^2 = \varphi^2_{II}$.}
    \label{fig:num-macro2D-pressure-velocity}
\end{figure}

The time evolution of the macroscopic pressure and velocity fields at
simulation  time levels $t_k=0.14, 0.28, 0.42, 0.56$\,s for the prescribed electric
potential $\varphi^2_{I}$ is depicted in
Fig.~\ref{fig:num-macro2D-pressure-velocity}.

%% file: aux-Appendix-pz-porel-short.tex
\appendix
\section{\newE{Direct numerical simulations -- incremental formulation}}\label{sec-DNS}
We present the weak formulation of the linearized problem arising from \eq{eq-fsi1a} adapted for the fluid-structure interaction problem and the heterogeneous medium introduced in Section~\ref{sec:prob}, see Remark~\ref{rem-bi-fsi}.

We consider time levels $t_0 < t_1 < t_2,\dots$, such that backward approximations in time are used to define $\dot \ub(t_k,x) \approx (\ub^k - \ub^{k-1})/\Dlt t$. Within each time level, the ``loads'' are presented by $\dlt \bar p$ and by prescribed electric potential $\dlt \bar \vphi$. At each time level $t_k$, solutions $\ssb(t_k,x) \approx \ssb^k(x)$, where $\ssb^k(x^i) =\ssb^{k-1}(x^0) + \Dlt \ssb(x^0)$, where $x^0 \in \Om_{k-1}$, and $x^i$ is the spatial position of $x^0$ for the $i$-iteration, such that
\begin{equation}\label{eq-fsi2}
  \begin{split}
    \Om(t_{k-1}) = \Om_0, \quad \Om(t_k) \approx \Om_i,\quad i = 0,1,2,\dots \\
    \ssb(x^i) \approx \ssb^k(\hat x),\quad x^i \approx \hat x = x^0 + \Dlt x,
 \end{split}
\end{equation}
We shall now consider the following simplified problem where the ``iterations'' are merged with ``time steps''. The aim is to find a new equilibrium for perturbed loads $\dlt \bar p$ and electric potential $\dlt \bar \vphi$, as stated below. The corresponding solution perturbations $\dlt\ub,\dlt\wb,\dlt p$ and $\dlt\vphi$ must satisfy the following equations,
\begin{equation}\label{eq-fsi3}
  \begin{split}
    \int_{\Om_R^s}\left(\Dop\eeb{\dlt\ub}  - \ul{\gb}^T\cdot\nabla\dlt\vphi\right):\eb_R(\tilde \vb) \dV   &     \\
    + \int_{\Om_R^f}\left(\mu \eeb{\wb + \dlt\wb + \dot\ub + \dlt\dot\ub} - (p +\dlt p) \Ib\right):\eb_R(\tilde \vb) \dV & \\+ \int_{\Om_R^s}\sigmabf_R^s(\ub,\vphi):\eb_R(\tilde \vb)  \dV
   & =  \int_{\pd_\sigma \Om_R^f} (\bar p + \dlt \bar p)\nb\cdot \tilde \vb \dS\;,\\
    \int_{\Om_R^s}\left(\ul{\gb}:\eeb{\dlt\ub} + \db\nabla \dlt\vphi\right)\cdot\nabla\psi\dV & = - \int_{\Om_R^s}\ul{D}_R(\ub,\vphi)\cdot\nabla\psi  \dV\;,\\
    \int_{\Om_R^f}\left(\mu \eeb{\wb + \dlt\wb + \dot\ub + \dlt\dot\ub} - (p +\dlt p) \Ib\right):\eb_R(\vtheta)  \dV& =  \int_{\pd_\sigma \Om_R^f} (\bar p + \dlt \bar p) \nb \cdot \vtheta\dS\;,\\
     \int_{\Om_R^f}q \nabla\cdot(\wb + \dlt\wb + \dot\ub + \dlt\dot\ub)  \dV & = 0\;,
  \end{split}
\end{equation}
to hold for any $\tilde \vb \in V_0(\Om_R)$, $\psi \in \Psi_0(\Om_R)$, and $\vtheta \in W_0(\Om_R)$, $ q \in Q_0(\Om_R)$. Above the solid velocity $\dot\ub$ is considered as the extension in the fluid, whereby,
\begin{equation}\label{eq-fsi2a}
  \begin{split}
 \ub^k & \approx \bar\ub^k + \dlt \ub\;,\quad \dot\ub^k \approx (\bar\ub^k - \ub^{k-1})/\Dlt t + \dlt \ub/\Dlt t \;,
 \end{split}
\end{equation}
where $\bar\ub^k$ is the current approximation of the unknown $\ub^k$.

The incremental problem is defined, as follows. Assume the reference configuration $\Om_R$ is given, for
the load due to the pressure increment $\dlt\bar p$ acting at $\pd_\sigma \Om_R^f$, and the actual electric potential $\dlt \bar\vphi$ imposed at the electrodes, compute the solution increments $\dlt \ssb = (\dlt \ub,\dlt \vphi, \dlt p, \dlt \wb) \in \Wcal(\Om_R) = V(\Om_R)\times\Psi(\Om_R) \times Q(\Om_R) \times W(\Om_R)$, such that \eq{eq-fsi3} is satisfied.

The stress is updated for the ``next'' iteration / time level, as follows
\begin{equation}\label{eq-fsi4}
  \begin{split}
    \Om_R^* = \Om_R + \{\dlt \ub\}_{\Om_R}\;,\quad x:= x^* = x + \dlt \ub(x)\;, \\
    \taubf_R^*(x^*) = \taubf_R(x) + \sym\{\taubf_R\nabla \dlt \ub \}(x) + \dlt_R^\Lcal \taubf_R(x)\;, \\
    \dlt_R^\Lcal \taubf_R^s = \Dop\eeb{\dlt\ub}  - \ul{\gb}^T\cdot\nabla\dlt\vphi\;.
 \end{split}
\end{equation}
In the expression for $\taubf_R^*(x^*)$, the term $(\nabla \dlt \ub) \taubf_R$ can be neglected in a case of very small deformation.

\paragraph{Time increments and iterations}
The load $\dlt \bar p := \Dlt \bar p$ and $\dlt \bar \vphi = \Dlt \bar \vphi$ can be applied at the first iteration $i=1$ of each time step $\Dlt t =  t_k - t_{k-1}$, see \eq{eq-fsi2}, which produces the out-of-balance in \eq{eq-fsi3}.  For the iterations $i = 2,3, \dots$, $\dlt \bar p = 0$ and $\dlt \bar \vphi = 0$, such that $\dlt\ssb^i \rightarrow 0$  with increasing $i$.



\input{aux-Appx-nl-problem1.tex}

\section{Sensitivity analysis of piezo-poroelastic coefficients}\label{sec-sa}


The sensitivity analysis of the homogenized coefficients is based on the
so-called ``design velocity method'' of the shape sensitivity which is well
known in the Structural Optimization,
\cite{haslinger_makinen_2003:_introduction_shape_optimization,haug_choi_komkov_1986:_design_sensitivity_analysis,haslinger_neittaanmaki_1988:_finite_element_approximation}.
A ``flux'' of material points $y \in Y$ is given in terms of a differentiable
and Y-periodic vector field $\vec\Vcal(y)$, such that $z_i(y,\tau) = y_i +
\tau\Vcal_i(y)$, $y \in Y$, $i=1,2$, where $\tau$ is the ``time-like'' variable
describes the perturbed position of $y\in Y$. We use the notion of the
following derivatives: $\delta(\cdot)$ is the total ({\it material})
derivative, $\delta_\tau(\cdot)$ is the partial ({\it local}) derivative \wrt
$\tau$. These derivatives are computed as the directional derivatives in the
direction of $\vec\Vcal(y)$, $y \in Y$.

\subsection*{Preliminary results on the shape sensitivity}

We shall use some particular results of the domain method of the shape sensitivity analysis which relies on the convection velocity field $\vec\Vcal$.
  Auxiliary results. Using \eq{eq-mip0a},
\begin{equation}\label{eq-sa5}
  \begin{split}
    \dltsh|Y_d| = \int_{Y_d} \nabla_y\cdot\vec\Vcal \dVy = \int_{\pd Y_d} \vec\Vcal\cdot\nb\dVy\;,
    \quad \dltsh \phi_f = |Y|^{-1} (\dltsh|Y_c| - \phi_f\dltsh|Y|)\;,\\
  \dltsh\intY_{\Gamma_f} \vb\cdot\nb^\fx \dSy  =  - \dltsh\intY_{Y_m}\nabla_y\cdot\vb\dVy 
   = \intY_{Y_m} \pd_i^y\Vcal_k \pd_k^y v_i \dVy\\
  - \intY_{Y_m}\nabla_y\cdot\vb \nabla_y\cdot\vec\Vcal \dVy + |Y|^{-1}\dltsh|Y|\intY_{Y_m}\nabla_y\cdot\vb\dVy
  \;.
\end{split}
\end{equation}
Note that $\dltsh|Y| = 0$ when $\vec\Vcal$ is $Y$-periodic.

For the bilinear forms involved in the local problems, the following sensitivity expressions hold:
 \begin{equation}\label{eq-sa6}
\begin{split}
\dltsh \aYms{\ub}{\vb} &= \intY_{Y_m\cup Y_*} 
D_{irks}\left(\delta_{rj}\delta_{sl} \nabla_y\cdot \vec\Vcal - \delta_{jr}\pd_s^y \Vcal_l - \delta_{ls}\pd_r^y \Vcal_j\right)
e_{kl}^y(\ub) e_{ij}^y(\vb)\dVy\\
& \quad - \aYms{\ub}{\vb}\intY_{Y} \nabla_y\cdot \vec\Vcal\dVy\;,\\
\end{split}
 \end{equation}
\begin{equation}\label{eq-sa7}
\begin{split}
\dltsh \dYm{\vphi}{\psi} &= \intY_{Y_m} 
d_{rs}\left(\delta_{ri}\delta_{sk} \nabla_y\cdot \vec\Vcal - \delta_{ks}\pd_r^y \Vcal_i - \delta_{ir}\pd_s^y \Vcal_k\right)
\pd_k^y \vphi \pd_i^y \psi\dVy \\
& \quad - \dYm{\vphi}{\psi}\intY_{Y} \nabla_y\cdot \vec\Vcal\dVy\;,\\
\end{split}
 \end{equation}
\begin{equation}\label{eq-sa8}
\begin{split}
\dltsh \gYm{\ub}{\psi} &= \intY_{Y_m} 
g_{jks}\left(\delta_{ij}\delta_{sl} \nabla_y\cdot \vec\Vcal - \delta_{ij}\pd_s^y \Vcal_l  - \delta_{sl}\pd_j^y \Vcal_i\right)
e_{kl}^y(\ub)\pd_i^y \psi\dVy \\
& \quad - \gYm{\ub}{\psi}\intY_{Y} \nabla_y\cdot \vec\Vcal\dVy\;,\\
\end{split}
\end{equation}
 
We shall use the following notation:
\begin{equation}\label{eq-mip0}
\begin{split}
\Xibf^{ij} = \omegabf^{ij} + \Pibf^{ij}\;, \mbox{ so that }  \dlt\Xibf^{ij} = \dlt\omegabf^{ij} + \dltsh\Pibf^{ij}\;.
\end{split}
\end{equation}

\subsection*{Differentiation of the microproblems}

Using the identities of the local characteristic problems
\eq{eq-mip1}-\eq{eq-mip4}, and \eq{eq-A10c}, and their differentiated forms
derived above, the sensitivity expressions of $\dlt\Hop\circ\vec{\Vcal}$ for
all the homogenized coefficients \eq{eq-HC1} can be derived; note that the
sensitivity expressions of the poroelastic coefficients treated in
\cite{Rohan-Lukes-nlBiot2015} cannot be used since, in our model, all theses
are modified by the skeleton piezoelectricity.

\paragraph{Elasticity}
\begin{equation}\label{eq-HC4a}
\begin{split}
 \dlt A_{klij}^H &  = 
 \dltsh\aYms{\Xibf^{ij}}{\Xibf^{kl}} - \dltsh\dYm{\hat\eta^{kl}}{\hat\eta^{ij}}
  + \aYms{\Xibf^{ij}}{\dltsh\Pibf^{kl}} + \aYms{\dltsh\Pibf^{ij}}{\Xibf^{kl}}\\
  & \quad - \left(\dltsh\gYm{\Xibf^{kl}}{\hat\eta^{ij}}  + \dltsh\gYm{\Xibf^{ij}}{\hat\eta^{kl}} + \gYm{\dltsh\Pibf^{kl}}{\hat\eta^{ij}}   + \gYm{\dltsh\Pibf^{ij}}{\hat\eta^{kl}}  \right),
\end{split}
\end{equation}

\paragraph{Biot modulus}
\begin{equation}\label{eq-HC7}
\begin{split}
  \dlt M^H &  = 
  \gamma \dltsh \phi_f  - 2\dltsh \intY_{\Gamma_f} \omegabf^P\cdot\nb^\fx\dSy\\
& \quad + \dltsh\dYm{\hat\eta^P}{\hat\eta^P} + 2\dltsh\gYm{\omegabf^P}{\hat\eta^P}-\dltsh\aYms{\omegabf^P}{\omegabf^P}\;.
\end{split}
\end{equation}

\paragraph{Biot coupling coefficients}
\begin{equation}\label{eq-HC13}
\begin{split}
  \dlt B_{ij}^H &  = 
  \dltsh \phi_f \delta_{ij} - \dltsh\intY_{Y_m} \nabla_y\cdot \omegabf^{ij}\dVy
  + \aYms{\omegabf^P}{\dltsh\Pibf^{ij}} + \dltsh\aYms{\omegabf^P}{\Xibf^{ij}}\\
 & \quad  -\left(\dltsh\gYm{\omegabf^P}{\hat\eta^{ij}}
   + \dltsh\gYm{\Xibf^{ij}}{\hat\eta^P}   + \gYm{\dltsh\Pibf^{ij}}{\hat\eta^P} + \dltsh\dYm{\hat\eta^P}{\hat\eta^{ij}}\right)\;.
\end{split}
\end{equation}

\paragraph{Stress-voltage coupling coefficients}
\begin{equation}\label{eq-HC14a}
\begin{split}
  \dlt H_{ij}^\kkk &  =  \aYms{\hat\omegabf^k}{\dltsh\Pibf^{ij}} - \gYm{\dltsh\Pibf^{ij}}{\hat\vphi^\kkk}
  - \dltsh\gYm{\hat\omegabf^\kkk}{\hat\eta^{ij}} - \dltsh\dYm{\hat\vphi^\kkk}{\hat\eta^{ij}}\\
& \quad + \dltsh\aYms{\hat\omegabf^\kkk}{\omegabf^{ij}+\Pibf^{ij}} - \dltsh\gYm{\omegabf^{ij}+\Pibf^{ij}}{\hat\vphi^\kkk}\;.
\end{split}
\end{equation} 

\paragraph{Pressure-voltage coupling coefficients}
\begin{equation}\label{eq-HC20}
\begin{split}
  \dlt Z^\kkk &  = \dltsh\gYm{\omegabf^P}{\hat\vphi^\kkk} - \dltsh\aYms{\hat\omegabf^k}{\omegabf^P} \\
  & \quad \dltsh\gYm{\hat\omegabf^\kkk}{\hat\eta^P} + \dltsh\dYm{\hat\vphi^\kkk}{\hat\eta^P}
  - \dltsh\intY_{\Gamma_f} \hat\omegabf^\kkk\cdot \nb^\fx\dSy\;.
\end{split}
\end{equation}

\begin{equation}\label{eq-S11}
\begin{split}
\delta K_{ij} & = \intY_{Y_c}\left( \ww_i^j + \ww_j^i- \nabla_y\wb^i:\nabla_y\wb^j
+\pi^i \nabla_y \cdot \wb^j +\pi^j \nabla_y \cdot \wb^i\right)\nabla_y\cdot \vec\Vcal\dVy \\
& \quad - \intY_Y\nabla_y\cdot \vec\Vcal\dVy
\intY_{Y_c}\nabla_y\wb^i:\nabla_y\wb^j\dVy
 \\
& \quad + \intY_{Y_c} \left( 
\pd_l^y\Vcal_r \pd_r^y\ww_k^i\pd_l^y\ww_k^j + \pd_l^y\Vcal_r \pd_r^y\ww_k^j\pd_l^y\ww_k^i - 
\pi^i \pd_k^y\Vcal_r \pd_r^y\ww_k^j
-  \pi^j \pd_k^y\Vcal_r \pd_r^y\ww_k^i 
\right)\dVy\;.
\end{split}
\end{equation}
Although $\delta K_{ij}$ depends formally on $\vec\Vcal$, this field can
be constructed arbitrarily in the interior of $Y_c$ without any
influence on $\delta K_{ij}$.

Note that the sensitivities of the boundary integrals in \eq{eq-HC7} and
\eq{eq-HC20} are evaluated according to \eq{eq-sa5}.

%% file: aux-Appx-nl-problem1.tex
\section{\newE{Homogenization of poroelastic medium described in deforming configuration}}\label{app-incr}
We consider simplified quasistatic problem for porous linear elastic medium with periodic occluded (disconnected) porosities saturated by (imobile) fluid. The aim is to justify the sensitivity analysis approach employed to introduce the ``reduced order model'' with approximated poroelastic homogenized coefficients. By virtue of the equilibrium \eq{eq-fsi1a}, the micromodel of the heterogeneous structure

We now consider a quasi-periodic structure in the deformed configuration with the unfolded domain $\Om_R\times Z(x)$, for any $x \in \OmR$, the microscopic configuration $Z(x)$ is defined by the local deformation of the initial reference cell $Y$.  In what follows, we write simply $Z$ instead of $Z(x)$.
The decomposition of $Y$ into the solid, $Y_S$, and the fluid, $Y_F$ applies, hence also $Z = Z_S \cup Z_F \cup \Gamma_Z$, where $\Gamma_Z$ is the solid-fluid interface.
In the limit, $\veps\rightarrow 0$, the homogenization in the deformed configuration leads to the two-scale equation
\begin{equation}\label{eq-E4}
  \begin{split}
    \int_\OmR \intY_{Z_S}\hat\Dop [\eebx{\dlt\ub^0} + \eebz{\dlt\hat\ub^1}]:[\eebx{\vb^0} + \eebz{\vb^1}]\dVxz\\
    - \int_\OmR \dlt p^0 \intY_{Z_F} (\nabla_x\cdot\vb^0 + \nabla_z\cdot\vb^1)\dVxz = \\
    \int_\OmR p^0 \intY_{Z_F} (\nabla_x\cdot\vb^0 + \nabla_z\cdot\vb^1)\dVxz -
    \int_\OmR \intY_{Z_S} \sigmabf_S^0:[\eebx{\vb^0} + \eebz{\vb^1}]\dVxz\;,  \\
    \int_\OmR q^0 \intY_{Z_F} (\nabla_x\cdot\dlt \ub^0 + \nabla_z\cdot \dlt \hat\ub^1)\dVxz = 0\;,
\end{split}
\end{equation}
for all test displacements $\vb^0 \in V_0(\Om)$, $\vb^1 \in \Hpdb(Z)$ and all $q^0 \in L_2(\Om)$.

\subsection{Local problems transformed in the initial configuration}
The local problem in the deformed configuration $Z(x)$ is obtained in analogy with the linear model considered in the fixed (initial) configuration; it reads, as follows: Find $\dlt\ub^1\in \Hpdb(Z)$, such that
\begin{equation}\label{eq-E4loc}
  \begin{split}
  \intY_{Z_S}\hat\Dop [\eebx{\dlt\ub^0} + \eebz{\dlt\hat\ub^1}]:\eebz{\vb^1}\dVz 
     + \dlt p^0 \intY_{Z_S} \nabla_z\cdot\vb^1 \dVz& = - r_Z(x,\vb^1)\;,\\
\end{split}
\end{equation}
holds for a.a. $x\in \Omega$, for all test displacements  $\vb^1 \in \Hpdb(Z)$. Note that $\hat\Dop(z) = \Dop(y)$, whereby $z = \psi(y, \cdot)$, the spatial position of $y \in Y$, is defined below, being dependent on the macroscopic position $x$.
Assuming the local equilibrium holds for a current approximation of the stress $\sigmabf_S^0$ and fluid pressure $p^0$, the out-of-balance vanishes,
\begin{equation}\label{eq-E4a}
  \begin{split}
 r_Z(x,\vb^1) =  \intY_{Z_S} \sigmabf_S^0:\eebz{\vb^1}\dVz + 
  p^0 \intY_{Z_R}  \nabla_z\cdot\vb^1\dVz \approx 0\;,\quad \forall \vb^1 \in \Hpdb(Z) \;,
 \end{split}
\end{equation}
so that the ``out of balance'' $R_Z(\vb^0)$ deduced from the terms in \eq{eq-E4}, involving the averaged stress $\bar\sigmabf_S$ and $p^0$ being tested by $\vb^0$  is respected only on the macroscopic level, namely in \eq{eq-fsi-Mac0} and \eq{eq-fsi-Mac1}.
Due to the linearity, local characteristic responses $\hat\chibf^{ij}$ and $\hat\chibf^P$ can be introduced, such that 
\begin{equation}\label{eq-E6}
  \begin{split}
 \dlt\hat\ub^1 = \hat\chibf^{ij}e_{ij}^x(\dlt\ub^0) + \hat\chibf^P \dlt p^0\;.
 \end{split}
\end{equation}
Upon substituting in \eq{eq-E4loc}, the local characteristic problems are obtained which have to be solved at a.e. $x \in \Om$: 
Find $\hat\chibf^{ij} \in \Hpdb(Z)$ and $\hat\chibf^P \in \Hpdb(Z)$, such that (note $\hat\Pibf^{ij} = \Pibf^{ij}(z)$, $\Pi_k^{ij}(z) = z_j\dlt_{ik}$)
\begin{equation}\label{eq-E7}
  \begin{split}
    \intY_{Z_S}\hat\Dop \eebz{\hat\chibf^{ij} + \hat\Pibf^{ij}}:\eebz{\vb^1}\dVz & = 0\;,\quad \forall \vb^1\in \Hpdb(Z)\;, \\
     \intY_{Z_S}\hat\Dop \eebz{\hat\chibf^P}:\eebz{\vb^1}\dVz & = -\intY_{Z_S} \nabla_z\cdot\vb^1\dVz\;,\quad \forall \vb^1\in \Hpdb(Z)\;.
\end{split}
\end{equation}
Homogenized coefficients $\hat\Aop,\hat\Bb$ and $\hat M$ analogous to those involved in the linear model, see \eq{eq-HC1}, can be computed at a.e. $x \in \Om$ to constitute the macroscopic problem, as explained below. The spatial (deformed) micro-configurations are updated using 
\begin{equation}\label{eq-fsi7}
  \begin{split}
       \dlt\hat\ub^\mic & = \dlt\hat\ub^1 + \Pibf^{ij}(z)e_{ij}^x(\dlt\ub^0(x))\;.\quad \\
\end{split}
\end{equation}
%
Hence, the micro-displacement transforming $Y$ to the local micro-configuration $Z(x)$ is given by 
\begin{equation}\label{eq-fsi8}
\begin{split}
\hat\ub^\mic(z,\eb^x,p^0) =  \int_0^{\eb^x}(\hat\chibf^{ij}(\tilde z)+\Pibf^{ij}(\tilde z))\dlt\tilde e_{ij}(x) +  \int_0^{p^0} \hat\chibf^P(\tilde z) \dlt \tilde p^0 
\end{split}
\end{equation}
where $z =  \psi(y,\eb^x,p^0)$, $y \in Y$ 
(and $\tilde z=  \psi(y,\tilde\eb^x,\tilde p^0)$) depends on both the macroscopic variables. 
 Obviously, mapping $\psi:y\mapsto z$ is parameterized by $\xibf^x:=(\eb^x,p^0)$, however, we shall write $z = \psi(y)$ for the sake of brevity.
Using $\ub^\mic(y) = \hat\ub^\mic(\psi(y))$ and mapping $\psi$, the deformation gradient $\Fb(y) = \Ib + \nabla_y\ub^\mic$ and its inverse $\hat\Fb = \Fb^{-1} = \Ib - \nabla_z\hat\ub^\mic$ can be established for  $y \in Y$ and $z = \psi(y)$. By virtue of the small deformation assumption, approximation $\hat\Fb \approx  \Ib - \nabla_y\ub^\mic$ can be used for $y\in Y$. Now, the transformed elasticity and the divergence operator can be defined: 
\begin{equation}\label{eq-fsi9}
  \begin{split}
    \tilde D_{irks}(y,\xibf^x) & = D_{ijkl}(\delta_{ls}-\pd_l^z \hat u^\mic_s)(\delta_{jr} - \pd_j^z  \hat u^\mic_r) \det(\Ib + \nabla_y \ub^\mic))\\
    & = D_{ijkl}(y)\hat F_{ls}(\psi(y))\hat F_{jr}(\psi(y))\,\det(\Fb(y))\;,\\
    \tilde \Dvg_y & = \det(\Fb)\Ib: (\hat\Fb\nabla_y)\;, \quad \tilde \pd_i^y =  \det(\Fb) \hat F_{ik}\pd_k^y\;.
\end{split}
\end{equation}


The local problem \eq{eq-E4loc} imposed in the deformed configuration can be transformed in the initial one using the approximations of $\Fb$ and $\hat\Fb$ described above. For a given local macroscopic response $\xibf^x = (\eb^x(\ub^0), p^0, \sigmabf_S^0)$ at  a.a. $x\in \Omega$, find $\dlt\ub^1\in \Hpdb(Y)$, such that
\begin{equation}\label{eq-fsi-loc0}
  \begin{split}
  \intY_{Y_S}\tilde\Dop(y,\xibf^x) \eeby{\tilde\Pibf^{kl} \dlt e_{kl}^x + \dlt\ub^1}:\eeby{\vb^1}\dVy & \\
    + \dlt p^0 \intY_{Y_S}\Ib: (\hat\Fb\nabla_y)\vb^1 \det(\Fb)\dVy& = - r_Y(x,\vb^1)\;,
\end{split}
\end{equation}
holds for all test displacements  $\vb^1 \in \Hpdb(Y)$; the \rhs functional vanishes due to \eq{eq-E4a}. Above $\tilde\Pibf^{kl}(y) = \Pibf^{kl}(z)$ is defined using $\hat\Pibf^{kl}$ introduced in \eq{eq-E6}, thus, $\tilde\Pi_j^{kl}(y) = (y_l + u_l^\mic(y))\dlt_{jk}$.

It is possible to consider  \eq{eq-fsi-loc0} and solve $\dlt\ub^1$ for a given ``pre-deformed'' configuration with $\tilde\Dop(y,\eb^x)$, $\Fb$ and $\hat\Fb$ being given.
In analogy with \eq{eq-E6} and using \eq{eq-fsi7} transformed in $Y_S$, we get
\begin{equation}\label{eq-E6Y}
  \begin{split}
    \dlt\ub^1 = \tilde\chibf^{ij}e_{ij}^x(\dlt\ub^0) + \tilde\chibf^P \dlt p^0\;,\\
    \dlt\ub^\mic = \tilde\Xibf^{ij}e_{ij}^x(\dlt\ub^0) + \tilde\chibf^P \dlt p^0\;,\mbox{ where } \tilde\Xibf^{ij} := \tilde\chibf^{ij} +  \tilde\Pibf^{ij}\;,
 \end{split}
\end{equation}
at any $y \in Y_S$, where $\tilde\chibf^{ij} \in \Hpdb(Y)$ and $\tilde\chibf^P \in \Hpdb(Y)$ satisfy (at any $x \in \Om_R$)
\begin{equation}\label{eq-locA}
  \begin{split}
    \intY_{Y_S}\tilde\Dop \eeby{\tilde\chibf^{ij} + \tilde\Pibf^{ij}}:\eeby{\vb^1} \dVy & = 0\;,\quad \forall \vb^1\in \Hpdb(Y)\;, \\
    \intY_{Y_S}\tilde\Dop \eeby{\tilde\chibf^P}:\eeby{\vb^1}\dVy & = -\intY_{Y_S}\tilde \Dvg_y\cdot\vb^1\dVy \;,\quad \forall \vb^1\in \Hpdb(Y)\;.
\end{split}
\end{equation}

\subsection{The macroscopic incremental problem}
For the macroscopic problem we use the ``spatial'' reference configuration $\Om(t) = \Om_R$, the domain $\Om_R$ is being updated using $\dlt\ub^0$, see \eq{eq-fsi4}. From \eq{eq-E4} we get
\begin{equation}\label{eq-fsi-Mac0}
  \begin{split}
\int_{\Om_R}  \intY_{Y_S}\tilde\Dop(y,\xibf^x) \eeby{\tilde\Pibf^{kl} e_{kl}^x(\dlt\ub^0) + \dlt\ub^1}:\eeby{\tilde\Pibf^{ij}}e_{ij}^x(\vb^0)\dVxy & \\
     - \int_{\Om_R} \dlt p^0 \MeanYs{\det(\Fb)}{Y_F} \nabla_x \cdot \vb^0 \dVx& = -R_Z(\vb^0)\;,\\
     \int_{\Om_R}  q\left(\MeanYs{\det(\Fb)}{Y_F} \nabla_x\cdot\dlt\ub^0 - \intY_{Y_S} \tilde \Dvg_y\cdot\dlt \ub^1\dVy\right) \dVx & = 0\;,
\end{split}
\end{equation}
for all $\vb^0 \in V_0(\Om_R)$ and $q \in L^2(\Om_R)$. Using the characteristic solutions, the poroelastic coefficients $\tilde\Aop = (\tilde A_{ijkl})$, $\tilde\Bb = (\tilde B_{ij})$ and $\tilde M$ can be identified,
\begin{equation}\label{eq-fsi-Mac0a}
  \begin{split}
    \tilde A_{ijkl} & = \intY_{Y_S}\tilde\Dop \eeby{\tilde\chibf^{kl} + \tilde\Pibf^{kl}}: \eeby{\tilde\chibf^{ij} + \tilde\Pibf^{ij}}\dVy\;,\\
    \tilde C_{ij} & = -\intY_{Y_S}\tilde\Dop \eeby{\tilde\chibf^P}:\eeby{\tilde\Pibf^{ij}}\dVy =  \tilde C_{ij}' = -\intY_{Y_S}\tilde\Dvg_y\tilde\chibf^{ij}\dVy\;,\\
    \tilde M & = - \intY_{Y_S}\tilde\Dvg_y\tilde\chibf^P \dVy = \intY_{Y_S}\tilde\Dop \eeby{\tilde\chibf^P}:\eeby{\tilde\chibf^P}\dVy\;,\\
    \tilde\Bb & = \tilde\phi_F \Ib + \tilde\Cb\;,\quad \tilde\phi_F = \MeanYs{\det(\Fb)}{Y_F}\;.
    \end{split}
\end{equation}
Recall that these $\tilde{~}$ marked coefficients depend on the macroscopic state represented by $\xibf^x = (\eebx{\ub^0}, p^0)$.
Now the macroscopic incremental  problem can be written in terms of the homogenized coefficients: Find $\dlt\ub^0 \in V(\Om_R)$ and $\dlt p^0 \in L^2(\Om_R)$, such that 
\begin{equation}\label{eq-fsi-Mac1}
  \begin{split}
    \int_{\Om_R}\left(\tilde\Aop \eebx{\dlt\ub^0} - (\Ib \tilde \phi_F + \tilde\Cb)\dlt p^0\right):\eebx{\vb^0} \dVx
     & =-R_Z(\vb^0)\;,\\
      \int_{\Om_R}q\left((\Ib \tilde \phi_F + \tilde\Cb):\eebx{\dlt\ub^0} + \tilde M \dlt p^0\right) \dVx & = 0\;,
 \end{split}
\end{equation}
for all $\vb^0 \in V_0(\Om_R)$ and $q \in L^2(\Om_R)$. To conclude, the outlined iterative computing scheme is based on commuting the macroscopic step \eq{eq-fsi-Mac1} imposed in updated $\Om_R$ according to \eq{eq-fsi4}, and solving the local problems \eq{eq-locA} with updated $\tilde\Dop$ and $\tilde\Dvg_y$ according to \eq{eq-fsi9}, involving $\ub^\mic$, \eq{eq-fsi8}. 

\subsection{Reduced order modelling based on the sensitivity analysis}
Using an approximation of the local characteristic responses, specific to $x \in \Om_R$, the computational complexity of the ``macro-micro'' commuting steps  can be avoided. The key is to approximate the  local deformation at the heterogeneity level. We employ the following abstract notation: 
\begin{equation}\label{eq-fsi16a}
  \begin{split}
\dlt\ub^\mic(y;\tilde\xibf^x) = \tilde\Xibf^*(\tilde\xibf^x) \dlt\xibf_*^x = (\Xibf^*(\zerobf) + \dlt_\xibf \Xibf^{*}(\zerobf)\circ\tilde\xibf^x)\dlt\xibf_*^x\;,
\end{split}
\end{equation}
where summation applies over ``index'' $*$ which stands for $I = {ij}$, or $P$, such that $\dlt\xibf^x = (\dlt\eb^x,\dlt p^0)$ and $\tilde\Xibf^* = (\tilde\Xibf^I,\tilde\chibf^P)$, recalling $\tilde\Xibf^I$ introduced in \eq{eq-E6Y}; the differential $\dlt_{\xibf} \Xibf^{*}\circ\xibf^x = \pd_{\xibf}^\odot\Xibf^{*}\xibf_\odot^x$.

Using the approximation  $\tilde\Xibf^*\approx \Xibf^*(\zerobf) + \dlt_{\xibf} \Xibf^{*}\circ\xibf^x$, \ie using the sensitivity determined in the initial configuration $\xibf^x = \zerobf$, 
\begin{equation}\label{eq-fsi16}
  \begin{split}
 \ub^\mic(y,\xibf^x)   & \approx \int_0^{\xibf^x} \left( \Xibf^{*}(y,\zerobf) +   \pd_{\xibf}^\odot\Xibf^{*}(y,\zerobf)\xibf_\odot^x \right)\dlt \tilde \xibf_* \\
   & \approx  \Xibf^*(y,\zerobf) \xibf_*^x + \frac{1}{2} \pd_{\xibf}^\odot\Xibf^{*}(y,\zerobf)  \xibf_\odot^x\otimes\xibf_*^x\;.
\end{split}
\end{equation}

For the 1st order approximation of $\ub^\mic$ 
\wrt the macroscopic response (as here represented by $\xibf^x$), we consider
\begin{equation}\label{eq-fsi17}
  \begin{split}
    \ub^\mic(x,y) \approx \vthetabf(1,x,y)\;,\quad \vthetabf(\tau,x,y) = \tau \hat\vthetabf(x,y)\;,\quad  \tau \in [0,1]\;,\\
    \mbox{ where } \quad    \hat\vthetabf(x,y) := \Xibf^*(y,\zerobf) \hat\xibf_*^x\;.
\end{split}
\end{equation}
Besides the approximation of $\tilde\Dop(y,\eb^x)$ established in \eq{eq-fsi9} using $\ub^\mic$, we employ 
\begin{equation}\label{eq-fsi15}
  \begin{split}
    \Ib: (\hat\Fb\nabla_y)\vb^1 \det(\Fb)& \approx
    \Ib: ((\Ib - \nabla_y\ub^\mic)\nabla_y)\vb^1 \det(\Ib + \nabla_y\ub^\mic) \\
    & \approx \Ib: ((\Ib - \nabla_y\ub^\mic)\nabla_y)\vb^1 (1 + \nabla_y\cdot\ub^\mic) \\
    & = \Ib: \left(\nabla_y\vb^1 - (\nabla_y\ub^\mic\nabla_y)\vb^1 + \nabla_y\vb^1 (\nabla_y\cdot\ub^\mic)\right) \\
   & \quad -\Ib: (\nabla_y\ub^\mic\nabla_y)\vb^1 (\nabla_y\cdot\ub^\mic)\;,
\end{split}
\end{equation}
where the last term can be neglected. 
In \eq{eq-fsi-loc0}, the following approximations relevant up to order $o(|\xibf^x|)$ defined in $Y$ are needed:
\begin{equation}\label{eq-fsi18}
  \begin{split}
    \hat\Fb& \approx \Ib -\nabla_y\vthetabf\;, \\  \Fb& \approx \Ib +\nabla_y\vthetabf\;,\\
    \Ib: (\hat\Fb\nabla_y)\vb^1 \det(\Fb)& \approx\Ib: \left(\nabla_y\vb^1 - (\nabla_y\vthetabf\nabla_y)\vb^1 + \nabla_y\vb^1 (\nabla_y\cdot\vthetabf)\right) \;,\\
    \MeanYs{\det(\Fb)}{Y_F}  & \approx \phi_F + \MeanYs{\nabla_y\cdot\vthetabf}{Y_F}\;,\\
    \tilde\Dop(\cdot,\xibf^x)& \approx \Dop + \left(\Dop \nabla_y\cdot\vthetabf - (\Dop\nabla_y\vthetabf)^T - \Dop\nabla_y\vthetabf\right)\;,\\
    \dlt_\vthetabf\tilde\Dop(\cdot,\xibf^x)\circ\vthetabf & \approx
    \Dop \nabla_y\cdot\vthetabf - (\Dop\nabla_y\vthetabf)^T - \Dop\nabla_y\vthetabf\;.
\end{split}
\end{equation}
Using this approximation, we can rewrite \eq{eq-fsi-loc0}, whereby 
the residual $r_Y(x,\vb^1) \approx 0$. Find $\dlt\ub^1\in \Hpdb(Y)$, such that
for a given macroscopic response perturbation $\dlt\xibf^x = (\dlt\eb^x,\dlt p^0)$ and given perturbation micro-displacement field $\vthetabf$, 
\begin{equation}\label{eq-fsi-loc1}
  \begin{split}
    \intY_{Y_S} \Dop\eeby{\tilde\Pibf^{kl} \dlt e_{kl}^x + \dlt\ub^1}:\eeby{\vb^1}\dVy \\
  +  \intY_{Y_S} \left(\Dop \nabla_y\cdot\vthetabf - (\Dop\nabla_y\vthetabf)^T - \Dop\nabla_y\vthetabf\right)\eeby{\tilde\Pibf^{kl} \dlt e_{kl}^x + \dlt\ub^1}:\eeby{\vb^1}\dVy \\
    + \dlt p^0 \intY_{Y_S} \Ib: \left(\nabla_y\vb^1 - (\nabla_y\vthetabf\nabla_y)\vb^1 + \nabla_y\vb^1 (\nabla_y\cdot\vthetabf)\right)\dVy = 0\;,\\
\end{split}
\end{equation}   
holds for all test displacements  $\vb^1 \in \Hpdb(Y)$. Solution to this problem can be approximated using the sensitivity analysis of the characteristic solutions computed in the initial configuration. As will be shown, solutions to \eq{eq-fsi-loc1} are not needed due to the approximation based on the sensitivity analysis.

For this, we use the local problems 
\eq{eq-locA} defined in the initial configuration, \ie for undeformed continuum, in analogy with the treatment in Section~\ref{sec-charpb}; $\chibf^{ij},\chibf^P\in\Hpdb(Y_S)$ satisfy
\begin{equation}\label{eq-E7Ya}
  \begin{split}
  \aYs{\chibf^{ij} + \Pibf^{ij}}{\vb} & = 0\;,\quad \forall \vb\in \Hpdb(Y_S)\;, \\
  \aYs{\chibf^P}{\vb} + \bYs{1}{\vb} & = 0\;,\quad \forall \vb\in \Hpdb(Y_S)\;,
\end{split}
\end{equation}
where $\Pi_k^{ij} = y_j \dlt_{ik}$ and 
\begin{equation}\label{eq-E7Y-bil}
  \begin{split}
  \aYs{\ub}{\vb} & = \intY_{Y_S}[\Dop \eeby{\ub}]:\eeby{\vb}\dVy\;,\\
  \bYs{q}{\vb} 
  & = \intY_{Y_S} q \Ib:\eeby{\vb}\dVy = \intY_{Y_S} q \nabla_y{\vb}\dVy\;.
\end{split}
\end{equation}
The homogenized coefficient associated with the undeformed configuration are given by (note, that $\Aop_{IK}$ combines $e_K^x(\dlt \ub^0)$ and $e_I^x(\vb^0)$, $\Cb$ combines $\dlt p^0$ and $e_I^x(\vb^0)$,
\begin{equation}\label{eq-HC1b}
  \begin{split}
    \Aop_{IK} = \intY_{Y_S}[\Dop \eeby{\Xibf^K}]:\eeby{\Pibf^I}\dVy = \aYs{\Xibf^K}{\Pibf^I}\;,\\
    \Cb_I = \intY_{Y_S}[\Dop \eeby{\chibf^P}]:\eeby{\Pibf^I}\dVy =  \aYs{\chibf^P}{\Pibf^I}\;,\\
    \Bb_I = \Cb_I + \phi_F \Ib\;.
\end{split}
\end{equation}

We shall now pursue the linearization procedure. Using \eq{eq-fsi18}, two-scale (macroscopic) problem \eq{eq-fsi-Mac0} is transformed in the following form:
\begin{equation}\label{eq-fsi-Mac2}
  \begin{split}
  \int_{\Om_R}  \intY_{Y_S} \Dop\eeby{\tilde\Pibf^{K} e_{K}^x(\dlt\ub^0) + \dlt\ub^1}:\eeby{\tilde\Pibf^{I}} e_{K}^x(\vb^0)\dVxy \\
    + \int_{\Om_R} \intY_{Y_S} \Dop_{\vthetabf}\eeby{\tilde\Pibf^{K} e_{K}^x(\dlt\ub^0) + \dlt\ub^1}:\eeby{\tilde\Pibf^{I}} e_{K}^x(\vb^0) \dVxy\\
    - \int_{\Om_R} \dlt p^0 (\phi_F + \MeanYs{\nabla_y\cdot\vthetabf}{Y_F}) \nabla_x \cdot \vb^0\dVx  = -R_Z(\vb^0)\;,\\
 \left(\phi_F + \MeanYs{\nabla_y\cdot\vthetabf}{Y_F}\right) \nabla_x\cdot\dlt\ub^0 - \intY_{Y_S}\Ib: \left(\nabla_y\dlt\ub^1 \right)\dVy \\
    -\intY_{Y_S}\Ib: \left( (\nabla_y\vthetabf\nabla_y)\dlt\ub^1 + \nabla_y\dlt\ub^1 (\nabla_y\cdot\vthetabf)\right) \dVy = 0\;,
\end{split}
\end{equation}
for all  $\vb^0 \in V_0(\Om_R)$ and $q \in L^2(\Om_R)$, where
 $\Dop_\vthetabf = \Dop \nabla_y\cdot\vthetabf - (\Dop\nabla_y\vthetabf)^T - \Dop\nabla_y\vthetabf$.
We proceed by substituting there the decomposition of $\dlt\ub^1$. Due to the perturbation by the micro-displacements $\vthetabf$, the following expressions hold,
\begin{equation}\label{eq-fsi-Mac3a}
  \begin{split}
   & \dlt\ub^\mic = \tilde\Pibf^{I}\dlt e_I^x +  \dlt\ub^1
    = \tilde\Xibf^*\dlt\xibf_*^x = (\Pibf^I+\pd_\xibf^*\Pibf^I \xibf_*^x)\dlt e_I^x  +  \dlt\ub^1   \;,\\
  &   \dlt\ub^1 = (\chibf^I + \dlt_\vthetabf\chibf^I\circ \vthetabf)\dlt e_I^x + (\chibf^P + \dlt_\vthetabf\chibf^P\circ \vthetabf)\dlt p^0\;,
\end{split}
\end{equation}
where the 2nd expression of $\dlt\ub^1$ uses $\vthetabf \approx \Xibf^*(\zerobf)\xibf_*$ with $\Xibf^* := (\Xibf^I,\chibf^P)$ and $\xibf = (e_I^x, p^0)$, see \eq{eq-fsi16a}.
Note that $\pd_\xibf^*\Xibf^{I,P}\xibf_*^x$ yields a linear expression in terms of $\xibf^x = (\eb^x(\ub^0),p^0)$ in our case. 

Further we express also the perturbation displacement $\vthetabf = \Xibf^*\xibf_*^x$ so that, see \eq{eq-fsi18},
\begin{equation}\label{eq-fsi-Mac4}
  \begin{split}
     \tilde\Dop(\cdot,\xibf^x)& \approx \Dop + \left(\Dop \nabla_y\cdot\vthetabf - (\Dop\nabla_y\vthetabf)^T - \Dop\nabla_y\vthetabf\right)\\
    & \approx  \Dop + \dlt_\xibf\tilde\Dop(\cdot,\zerobf)\circ\xibf^x \equiv \Dop + \pd_\xibf^*\Dop\xibf_*^x\;,
\end{split}
\end{equation}
where the last expression is introduced for the sake of brevity. Note $\MeanYs{\nabla_y\cdot\vthetabf}{Y_F} = \MeanYs{\nabla_y\cdot\Xibf^*}{Y_F}\xibf_*^x = \MeanYs{\nabla_y\cdot\Xibf^*\xibf_*^x}{Y_F}$. Substituting in \eq{eq-fsi-Mac2}, we get
\begin{equation}\label{eq-fsi-Mac5-I}
  \begin{split}
    \int_{\Om_R}  \intY_{Y_S} \left[(\Dop + \pd_\xibf^*\Dop\xibf_*^x)
      \eeby{\Xibf^K+\dlt_\xibf^*\Xibf^K \xibf_*^x}\right]:\eeby{\Pibf^I +\pd_\xibf^*\Pibf^I \xibf_*^x} e_K^x(\dlt\ub^0) e_I^x(\vb^0)\dVxy \\
    +   \int_{\Om_R}  \intY_{Y_S} \left[(\Dop + \pd_\xibf^*\Dop\xibf_*^x)
      \eeby{\chibf^P+\pd_\xibf^*\chibf^P \xibf_*^x}\right]:\eeby{\Pibf^I +\pd_\xibf^*\Pibf^I \xibf_*^x} \dlt p^0 e_I^x(\vb^0)\dVxy \\
        -\int_{\Om_R}\dlt p^0\left(\phi_F + \MeanYs{\nabla_y\cdot\Xibf^*\xibf_*^x}{Y_F}\right)\nabla_x \cdot \vb^0 \dVx = -R_Z(\vb^0)\;,
    \end{split}
\end{equation}
\begin{equation}\label{eq-fsi-Mac5-II}
  \begin{split}
      \int_{\Om_R}q\left(\phi_F + \MeanYs{\nabla_y\cdot\Xibf^*\xibf_*^x}{Y_F}\right) \nabla_x\cdot\dlt\ub^0 \dVx\\
       -\int_{\Om_R} q\intY_{Y_S}\nabla_y\cdot \left( \chibf^K+\pd_\xibf^*\chibf^K \xibf_*^x\right) e_K^x( \dlt\ub^0)\dVy \dVx\\
     - \int_{\Om_R}q\intY_{Y_S} \Ib: \left(\nabla_y(\Xibf^*\xibf_*^x)\nabla_y( \chibf^K+\pd_\xibf^*\chibf^K \xibf_*^x)\right) e_K^x( \dlt\ub^0)\dVy \dVx\\
    - \int_{\Om_R}q\intY_{Y_S}\nabla_y \cdot(\chibf^K+\pd_\xibf^*\chibf^K \xibf_*^x)\nabla_y\cdot(\Xibf^*\xibf_*^x) e_K^x( \dlt\ub^0)\dVy \dVx \\
  -\int_{\Om_R} q\intY_{Y_S}\nabla_y\cdot \left( \chibf^P+\pd_\xibf^*\chibf^P \xibf_*^x\right) \dlt p^0\dVy \dVx\\
       - \int_{\Om_R}q\intY_{Y_S} \Ib: \left(\nabla_y(\Xibf^*\xibf_*^x)\nabla_y( \chibf^P+\pd_\xibf^*\chibf^P \xibf_*^x)\right) \dlt p^0\dVy \dVx\\
       - \int_{\Om_R}q\intY_{Y_S}\nabla_y \cdot(\chibf^P+\pd_\xibf^*\chibf^P \xibf_*^x)\nabla_y\cdot(\Xibf^*\xibf_*^x) \dlt p^0\dVy \dVx
    = 0\;.
\end{split}
\end{equation}

The aim is now to show the relationship between the above linearized macroscopic problem and the formulation \eq{eq-fsi-Mac1}, where the perturbed  homogenized coefficients are approximated using the sensitivity analysis, as employed in the main part of this paper. For this, we rewrite \eq{eq-fsi-Mac5-I} and \eq{eq-fsi-Mac5-II}  in terms of the bilinear forms \eq{eq-E7Y-bil} neglecting higher order terms \wrt the macroscopic ``perturbation'' $\xibf_*^x$. The differential  $\dlt_\xibf()\circ\xibf^x$ is abbreviated simply by ``$\dlt$''  in the context of \eq{eq-fsi-Mac3a} (see also \eq{eq-HC2} below), whereas $\Dlt\ub^0:=\dlt\ub^0$ and $\Dlt p^0:=\dlt p^0$ to avoid confusion only in the following equation:
\begin{equation}\label{eq-fsi-Mac5a}
  \begin{split}
     \int_{\Om_R}\left[\aYs{\Xibf^K}{\Pibf^I} + \dlt\aYs{\Xibf^K}{\Pibf^I}\right] e_K^x(\Dlt\ub^0) e_I(\vb^0)\dVx\\
     \int_{\Om_R}\left[\aYs{\dlt\Xibf^K}{\Pibf^I} + \aYs{\Xibf^K}{\dlt\Pibf^I}\right] e_K^x(\Dlt\ub^0) e_I(\vb^0)\dVx\\
    + \int_{\Om_R}\left[\aYs{\chibf^P}{\Pibf^I} + \dlt\aYs{\chibf^P}{\Pibf^I}\right] \Dlt p^0 e_I(\vb^0)\dVx \\
    + \int_{\Om_R}\left[\aYs{\dlt\chibf^P}{\Pibf^I} + \aYs{\chibf^P}{\dlt\Pibf^I}\right] \Dlt p^0 e_I(\vb^0)\dVx \\
    - \int_{\Om_R}\dlt p^0 (\phi_F + \dlt\phi_F)\nabla\cdot\vb^0\dVx =  -R_Z(\vb^0)\;,\\
    \int_{\Om_R}q (\phi_F + \dlt\phi_F)\nabla\cdot\dlt\ub^0 \dVx\\
     -\int_{\Om_R}q\left[
 \bYs{1}{\chibf^K} + \bYs{1}{\dlt\chibf^K} + \dlt\bYs{1}{\chibf^K} 
 \right]e_K^x(\Dlt\ub^0)\dVx \\
    -\int_{\Om_R}q\left[
 \bYs{1}{\chibf^P} + \bYs{1}{\dlt\chibf^P} + \dlt\bYs{1}{\chibf^P} 
 \right]\Dlt p^0\dVx = 0\;.
\end{split}
\end{equation}
    
The sensitivity of the bilinear forms \wrt the perturbation by $\vthetabf$ (see \eq{eq-fsi17},
recalling  the multiindex $ij = I$, $K = kl$
) yields
\begin{equation}\label{eq-E7Yb}
  \begin{split}
  \dlt_\tau\aYs{\chibf^I + \Pibf^I}{\vb} +\aYs{\dlt\chibf^I + \dlt_\tau\Pibf^I}{\vb} & = 0\;,\\
  \dlt_\tau\aYs{\chibf^P}{\vb} + \dlt_\tau\bYs{1}{\vb} +\aYs{\dlt\chibf^P}{\vb} & = 0\;,\\
\end{split}
\end{equation}
where the differential $\dlt_\tau$ means $\pd/\pd \tau|_{\tau = 0}$ in the context of \eq{eq-fsi17},
\begin{equation}\label{eq-E7Yc}
  \begin{split}
  \dlt_\tau\aYs{\chibf^I}{\vb} =  \intY_{Y_S} \left(\Dop \nabla_y\cdot\vthetabf - (\Dop\nabla_y\vthetabf)^T - \Dop\nabla_y\vthetabf\right)\eeby{\chibf^I}:\eeby{\vb}\dVy\;,\\
  \dlt_\tau\bYs{1}{\vb} = \intY_{Y_S}\Ib: \left( -(\nabla_y\vthetabf\nabla_y)\vb + \nabla_y\vb(\nabla_y\cdot\vthetabf)\right)\dVy\;,
\end{split}
\end{equation}
From \eq{eq-E7Yb} and \eq{eq-E7Yb}, using $\Xibf^I = \chibf^I + \Pibf^I$, we get
\begin{equation}\label{eq-E7Yd}
  \begin{split}
    \aYs{\dlt\Xibf^K}{\Pibf^I} & = \dlt_\tau\aYs{\Xibf^K}{\chibf^I} + \aYs{\dlt\Pibf^K}{\Xibf^I}\;,\\
    \aYs{\dlt\chibf^P}{\chibf^K} & = - \dlt_\tau\aYs{\chibf^P}{\chibf^K} - \dlt_\tau\bYs{1}{\chibf^K}\;,\\
    \aYs{\dlt\chibf^P}{\Pibf^I} & = -\aYs{\dlt\chibf^P}{\chibf^I} =
    \dlt_\tau\aYs{\chibf^P}{\chibf^I} +  \dlt_\tau\bYs{1}{\chibf^I}\;,\\
    \aYs{\chibf^P}{\dlt\chibf^I} & = - \aYs{\chibf^P}{\dlt\Pibf^I} - \dlt_\tau\aYs{\chibf^P}{\chibf^I}\;.
\end{split}
\end{equation}
Due to \eq{eq-fsi17} introducing $\vthetabf$, the following identity holds,
\begin{equation}\label{eq-HC2}
  \begin{split}
 \dlt_\xibf^*\Xibf\circ \xibf_*^x = \pd_\vthetabf\Xibf\circ \dlt\vthetabf = \pd_\vthetabf\Xibf\circ\pd_\xibf^* \vthetabf \circ \xibf_*^x =
    (\dlt_\vthetabf\Xibf\circ\Xibf^*)\circ \xibf_*^x \;.
\end{split}
\end{equation}
To express $\dlt_\vthetabf \Aop$ and $\dlt_\vthetabf \Bb$, in \eq{eq-fsi-Mac5-I} and \eq{eq-fsi-Mac5-II}, we collect the corresponding terms which yield expressions up to order $o(|\vthetabf|)$, recalling  \eq{eq-HC2},
using \eq{eq-E7Yd}, we get from \eq{eq-fsi-Mac5a}$_1$ for the effective elasticity $\Aop_{IK}$,
\begin{equation}\label{eq-HC4A}
  \begin{split}
    \dlt_\vthetabf\Aop_{IK} & =  \dlt_\tau\aYs{\Xibf^K}{\Pibf^I} +\aYs{\dlt\Xibf^K}{\Pibf^I}+\aYs{\Xibf^K}{\dlt\Pibf^I} \\
    & = \dlt_\tau\aYs{\Xibf^K}{\Pibf^I} +  \dlt_\tau\aYs{\Xibf^K}{\chibf^I} + \aYs{\dlt\Pibf^K}{\Xibf^I} + \aYs{\Xibf^K}{\dlt\Pibf^I} \\
    & = \dlt_\tau\aYs{\Xibf^K}{\Xibf^I}  + \aYs{\dlt\Pibf^K}{\Xibf^I} + \aYs{\Xibf^K}{\dlt\Pibf^I}
    \;,
\end{split}
\end{equation}
for the Biot coupling coefficients $\Bb_I = \Cb_I + \phi_F \Ib$, where
\begin{equation}\label{eq-HC4B}
  \begin{split}
    \dlt_\vthetabf\Cb_I & =  -\left(\dlt_\tau\aYs{\chibf^P}{\Pibf^I} +\aYs{\dlt\chibf^P}{\Pibf^I}+\aYs{\chibf^P}{\dlt\Pibf^I} \right)\\
    & = -\left(\dlt_\tau\aYs{\chibf^P}{\Pibf^I} + \dlt_\tau\aYs{\chibf^P}{\chibf^I} +  \dlt_\tau\bYs{1}{\chibf^I} +\aYs{\chibf^P}{\dlt\Pibf^I}\right) \\
    & =  \-\left(\dlt_\tau\aYs{\chibf^P}{\Xibf^I}+ \dlt_\tau\bYs{1}{\chibf^I} +\aYs{\chibf^P}{\dlt\Pibf^I}\right)\;,
\end{split}
\end{equation}
and, thus
\begin{equation}\label{eq-HC4C}
  \begin{split}
    \pd_\vthetabf\Bb_I & = \pd_\vthetabf\Cb_I + \pd_\vthetabf\phi_F \Ib_I\;,\\
    \pd_\vthetabf\phi_F & =  \pd_\vthetabf\MeanYs{\nabla_y\cdot\vthetabf}{Y_F} = \MeanYs{\nabla_y\cdot\circ}{Y_F}\;.
\end{split}
\end{equation}
From \eq{eq-fsi-Mac5a}$_2$, in analogy
\begin{equation}\label{eq-HC5}
  \begin{split}
    \dlt_\vthetabf\Cb_I' & = -\dlt_\tau\bYs{1}{\chibf^I} -\bYs{1}{\dlt\chibf^I} \\
    & = 
   -\dlt_\tau\bYs{1}{\chibf^I} + \aYs{\chibf^P}{\dlt\Pibf^I} + \dlt_\tau\aYs{\chibf^P}{\chibf^I}\;,
\end{split}
\end{equation}
whereby $\dlt_\vthetabf\Cb_I' = \dlt_\vthetabf\Cb_I$, as expected. This can be seen using \eq{eq-E7Yd}$_4$ substituted in the final expression \eq{eq-HC4B}, which yields \eq{eq-HC5}. For the Biot compressibility, we obtain
\begin{equation}\label{eq-HC6}
  \begin{split}
    \dlt_\vthetabf M & = -\dlt_\tau\bYs{1}{\chibf^P} -\bYs{1}{\dlt\chibf^P} \\
    & = 
   -2\dlt_\tau\bYs{1}{\chibf^P} - \dlt_\tau\aYs{\chibf^P}{\chibf^P}\;.
\end{split}
\end{equation}
With these differentials and knowing $\dlt_\vthetabf(~) = \pd_\xibf^* (~)\xibf_*^x = \pd_\eb (~) \eb_x(\ub^0) + \pd_p (~) p^0$, the following coefficients, depending on the actual solution, can be introduced 
\begin{equation}\label{eq-HC7B}
  \begin{split}
    \ol{\Aop}(\xibf^x) := \Aop + \pd_\eb\Aop \eb_x(\ub^0) + \pd_p\Aop p^0 \;, \\
    \ol{\Bb}(\xibf^x) := \Bb + \pd_\eb\Bb \eb_x(\ub^0) + \pd_p\Bb p^0\;, \\
    \ol{M}(\xibf^x):= \pd_\eb M \eb_x(\ub^0) + \pd_p M p^0\;,
\end{split}
\end{equation}
where $\xibf^x = \bar\xibf^x + \dlt\xibf^x$ with $\dlt\xibf^x = (\eebx{\dlt\ub^0}, \dlt p^0)$ and $\bar\xibf^x = (\eebx{\bar\ub^0}, \bar p^0)$.
Now, the nonlinear problem arising from incremental formulation \eq{eq-fsi-Mac0} attains the following form:
\begin{equation}\label{eq-fsi-Mac3}
  \begin{split}
 \int_{\Om_R}\left(\ol{\Aop}(\xibf^x)\eebx{\dlt\ub^0} - \ol{\Bb}(\xibf^x)\dlt p^0\right):\eebx{\vb^0}\dVx & =-R_Z(\vb^0)\;,\\
      \int_{\Om_R}q\left(\ol{\Bb}(\xibf^x):\eebx{\dlt\ub^0} + \ol{M}(\xibf^x)\dlt p^0\right)\dVx & = 0\;,
 \end{split}
\end{equation}
for all $\vb^0 \in V_0(\Om_R)$ and $q \in L^2(\Om_R)$.

\begin{remark}\label{rem-SAHC}
The incremental, \ie a linearized formulation of the weakly nonlinear problem for the heterogeneous poroelastic medium introduced at the pore level in terms of the FSI problem in the deformed configuration was transformed into formulation \eq{eq-fsi-Mac3} using the homogenization in the deformed reference configuration and the approximation of the characteristic responses based on the response in the initial configuration. We have observed that this treatment leads to a consistent linearization of the problem where the involved homogenized coefficients are linear functions of the local macroscopic response $\xibf^x$ which can be obtained using the sensitivity analysis, as explained in a more detail in \ref{sec-sa}. The sensitivity analysis applies with the Y-periodic vector field $\hat\vthetabf = \vec\Vcal$, see \eq{eq-fsi17}. The formulae obtained for $\dlt_\vthetabf \Aop$, $\dlt_\vthetabf \Bb$, and $\dlt_\vthetabf M$ using \eq{eq-HC4A}-\eq{eq-HC6} are equivalent to those obtained for the ``full problem'' with piezoelectric material, formulae \eq{eq-HC4a}-\eq{eq-HC13} for vanishing piezoelectric properties (\ie neglecting $g_{ikl}$), whereby $\omegabf^{P} = - \chibf^P$.
  \end{remark}

%% file: aux-1Dmodel-R.tex
\section{1D example of the homogenized metamaterial model}\label{sec-1D}

Using a 1D model we can illustrate the modelling approach based on the
semi-linear extension of homogenized piezo-poroelastic medium controllable by
electric potential wave. \chE{Here we consider the nonlinearity affecting the
flow only, to demonstrate that the permeability dependence on the deformed
configuration is fundamental to capture the peristalsis driven flow.}

Suitable symmetric boundary conditions applied to the 3D model derived in subsequent sections lead to a reduced 1D model of the macroscopically homogeneous medium involving displacement  $u = u_1$, the fluid seepage $w=w_1$, the strain $e = e_{11} = u'$ and the electric potential $\vphi$ which are functions of $(x,t)$. We consider the 1D continuum represented by interval $X = ]x_-,x_+[$; upon denoting the spatial derivative by $\pd_x() \equiv ()'$,
\begin{equation}\label{eq-1d-1}
\begin{split}
  e \equiv e_{11}\;,\quad e_{ij} = 0\quad  \mbox{ for } i,j \not = 1\;, \\
  \nabla p = (p',0,0)\;,\quad \mbox{ and } \nabla \vphi = (\vphi',0,0)\;,\\
  \ub = (u,0,0)\;,\quad \wb = (w,0,0)\;.
\end{split}
\end{equation}
Since the inertia effects are neglected in the present paper, only the static permeability $K$ is required. Thus, we avoid the convolution integral arising in the dynamic Darcy flow model involving the dynamic permeability. Moreover, in a simple way, this quasistatic restriction enables to reduce the homogenized model represented by the following equations,
\begin{equation}\label{eq-1d-2}
  \begin{split}
   & -\sigma' = f\;, \mbox{ where }\sigma = A e(u) - p B + H \vphi\;, \\
  &  M \dot p +  w' + B e(\dot u) = Z \dot \vphi\;,\mbox{ where }
    w = - \tilde K p'\;. 
\end{split}
\end{equation}
Above the ``dot'' denotes the time derivative.
The stress $\sigma$ is given by the pressure $p$ and the potential $\vphi$ due to the extended poroelastic law involving coefficients $A,B$ and $H$.
Further we omit the volume forces $f$. 
The permeability depends on the response due to the linear expansion formula,
\begin{equation}\label{eq-1d-2b}
  \begin{split}
    \tilde K(e,p,\vphi) = K_0 + \pd_e K_0 e + \pd_p K_0 p + \pd_\vphi K_0 \vphi\;,
\end{split}
\end{equation}
which can be obtained using the sensitivity analysis of the microstructure response; $\pd_p K_0$ is the derivative of $K_0$ \wrt the pressure and, in analogy, the differentiation \wrt $e$ and $\vphi$ is denoted by $\pd_e K_0$ and $\pd_\vphi K_0$, respectively.
For other coefficients such expansions can be considered in analogy. 
The mass conservation involves the homogenized coefficients $M,B$ and $Z$. 

\begin{figure}
\centering 
\begin{tabular}{cc}
  \includegraphics[width=0.39\linewidth]{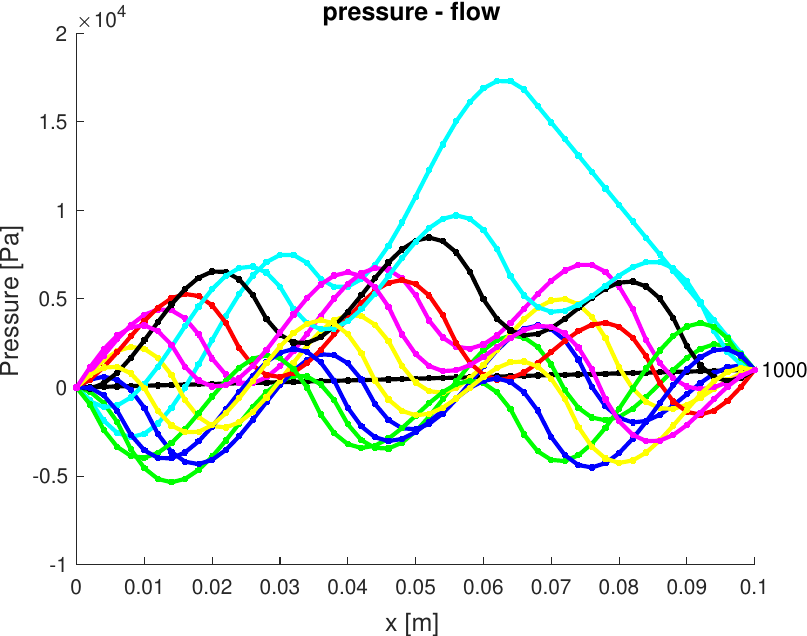} &
  \includegraphics[width=0.39\linewidth]{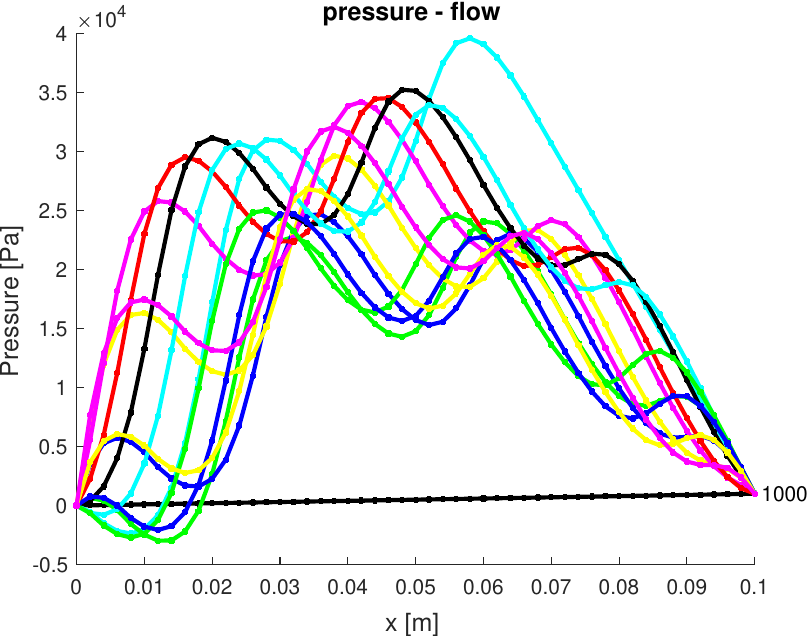} \\
  \includegraphics[width=0.39\linewidth]{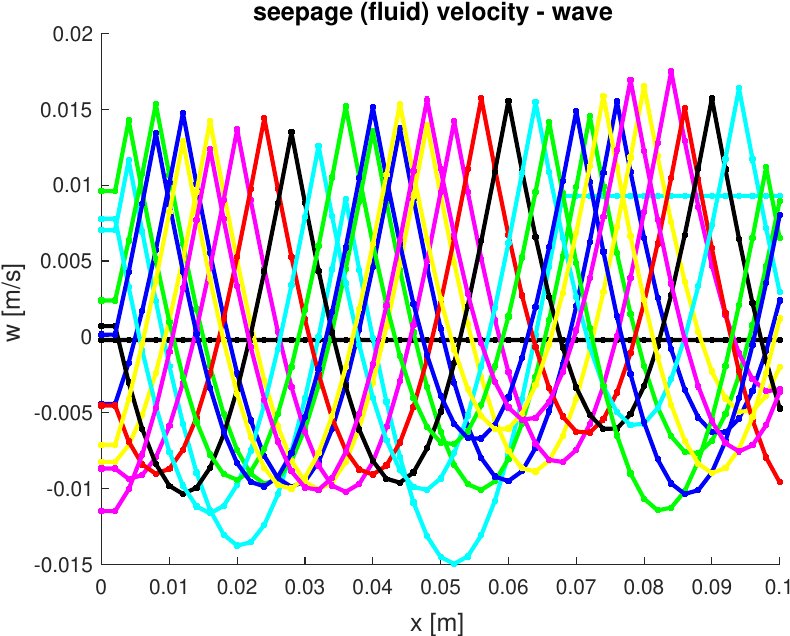} &
  \includegraphics[width=0.39\linewidth]{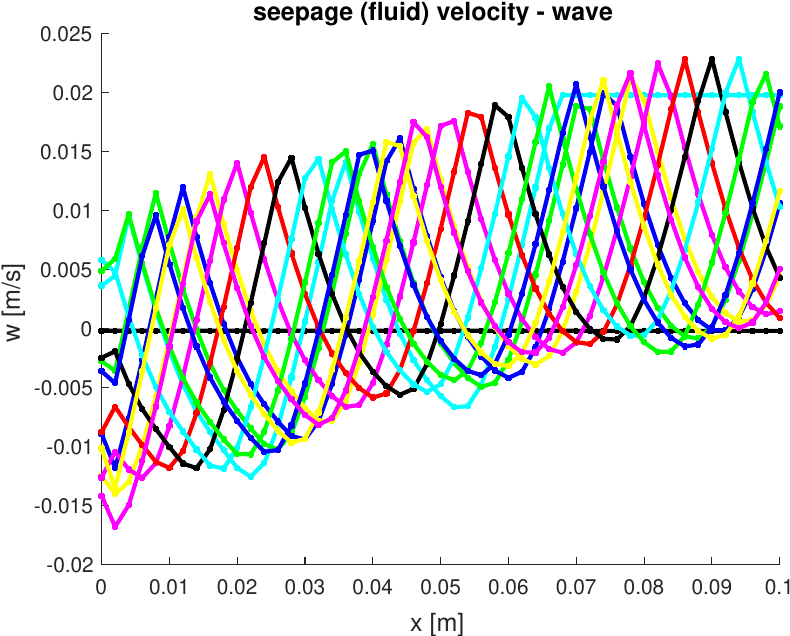}\\
  linear & nonlinear
\end{tabular}

\caption{Spatial distributions of pressure $p(\cdot,t)$ and seepage $w(\cdot,t)$ for various time levels (the time levels specific colours): linear (left) and nonlinear (right) models.}
  \label{fig-1D-pw}
\end{figure}

\begin{figure}
\centering 

\begin{tabular}{cc}
  \includegraphics[width=0.49\linewidth]{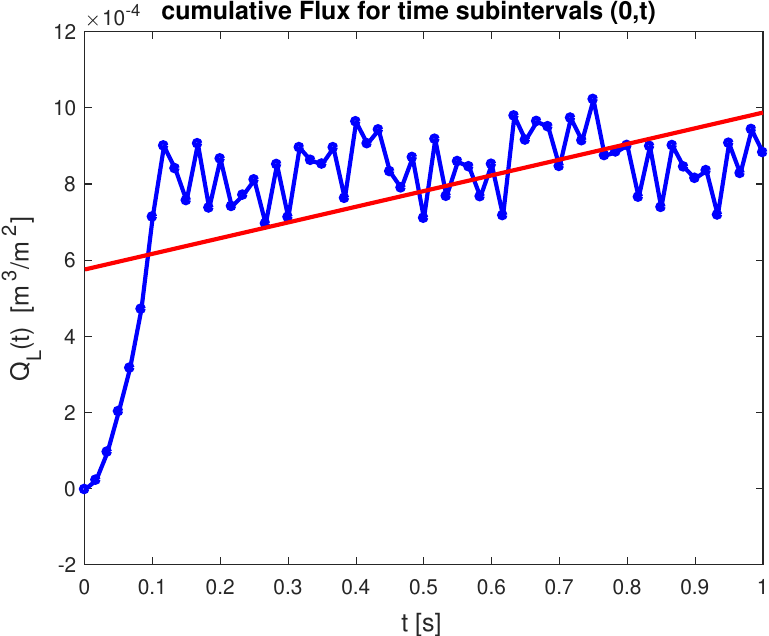} &
  \includegraphics[width=0.49\linewidth]{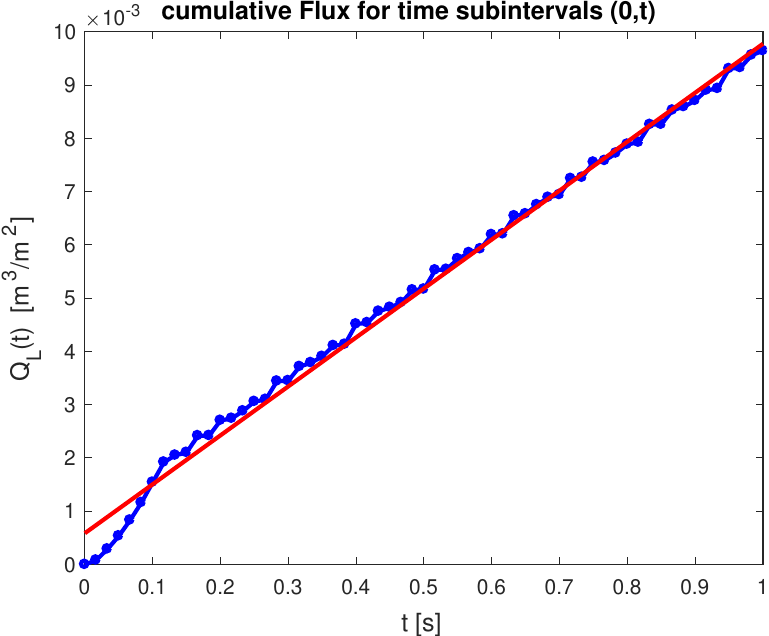} \\
   linear & nonlinear
\end{tabular}

 \caption{The cumulative fluxes $Q_+(t)$ [m$^3$/m$^2$] for the  linear and nonlinear models. The regression lines show the pumping capability.}
  \label{fig-1D-Q}
\end{figure}

Let us define
\begin{equation}\label{eq-1d-3}
\begin{split}
  C = M + BA^{-1} B\;, \quad F = Z + BA^{-1} H\;,\\
  K_p =  \pd_p K_0+\pd_e K_0A^{-1} B \;, \quad K_\vphi = \pd_\vphi K_0-\pd_e K_0A^{-1} H\;.
\end{split}
\end{equation}
When no volume forces and no inertia effects are considered, thus, $f = 0$ in \eq{eq-1d-2}$_1$, thus  $\bar \sigma' = 0$, the stress $\sigma(x,t) = \bar\sigma(t)$ in the 1D continuum. The strain can be eliminated from the mass conservation 
\begin{equation}\label{eq-1d-4}
\begin{split}
e = A^{-1}(B p - H \vphi + \bar\sigma)\;.
\end{split}
\end{equation}
Consequently, this enables to rewrite the expansion of permeability $\tilde K$ and its spatial derivative $\tilde K'$ in terms of $p$, $\bar\sigma$ and $\vphi$ only,
\begin{equation}\label{eq-1d-5}
  \tilde K = K_0 + \pd_e K_0A^{-1}\bar\sigma + K_p p + K_\vphi \vphi\;,\quad
  \tilde K' = K_p p' + K_\vphi \vphi'\;,
\end{equation}
noting $\bar \sigma' = 0$.
Using the Darcy law,  the linear and a non-linear flow models are obtained from \eq{eq-1d-2}$_2$,
\begin{equation}\label{eq-1d-6}
\begin{split}
  C \dot p - K_0 p'' & = F \dot \vphi\;, \\
  C \dot p - (\tilde K p'' +  \tilde K' p') & = F \dot \vphi\;, 
\end{split}
\end{equation}
where the nonlinearity arises due to the dependence of $\tilde K$ and $\tilde K'$ on $p$, as follows from \eq{eq-1d-5}.
Upon substituting $\tilde K$ and $\tilde K'$ in the nonlinear flow model \eq{eq-1d-6}$_2$ using \eq{eq-1d-5}, we get
\begin{equation}\label{eq-1d-7}
  C \dot p - K_0p'' - (K_p p' + K_\vphi \vphi') p'
  -( \pd_e K_0A^{-1}\bar\sigma + K_p p + K_\vphi \vphi ) p'' = F \dot \vphi\;.
\end{equation}
This illustrates, how the nonlinearity due to the ``deformation dependent
homogenized coefficients'' emerges. We shall apply the expansion $\tilde
K(e,p,\vphi)$ to all poroelastic coefficients, as detailed below in
Section~\ref{sec-homog}. The numerical results confirm the necessity to account
for this kind of nonlinearity to mimic the pumping effect which enables for the
fluid flow against the pressure slope. In Figs.~\ref{fig-1D-pw} and
\ref{fig-1D-Q}, this statement is supported by a numerical example with
prescribed pressures $p(x_-,t) = 0$ and $p(x_+,t) = \bar p = 0.1$ kPa, whereby
the controlling potential $\vphi = \bar\vphi_0 |\sin(\om t - k x)|$ with
$\bar\vphi_0 = -100$kV, $c = \om/k = 0.8$ m/s. Solutions of \eq{eq-1d-7} for
the linear and nonlinear cases are illustrated in Fig.~\ref{fig-1D-pw}, in
terms of the pressure $p(x,t)$ and the seepage $w(x,t)$. In
Fig.~\ref{fig-1D-Q}, we depict the cumulative flux $Q_+(t)$ (at the right end)
which is expressed at the endpoints $x_-$ and $x_+$ for times $t>0$ by
\begin{equation}\label{eq-1d-8}
\begin{split}
  Q_\pm(t) = \pm\int_0^t w(x_\pm,\tau) \dd \tau\;,
\end{split}
\end{equation}
noting that $Q_+(t) \approx Q_-(t)$ for times $t > t_0$, whereby $t_0$ is ``large enough'' to suppress effects of the initial conditions.
It is readily seen that $Q(t)$ is oscillating about a constant in the linear case, \revE{ \eq{eq-1d-6}$_1$}, indicating vanishing transport, \revE{whereas, in the non-linear case, \eq{eq-1d-7},} $Q(t)$ shows a linear average growth, hence the pumping effect.


